\RequirePackage{rotating}

\documentclass[10pt]{amsart}
\usepackage{amssymb,amsmath,amsthm}
\usepackage{bm}
\usepackage{fullpage}
\usepackage{enumerate}
\usepackage{verbatim}
\usepackage{graphicx,float,rotating}
\usepackage[mathscr]{eucal}

\usepackage[hidelinks]{hyperref}

\usepackage{amsrefs}

\newtheorem{theorem}[equation]{Theorem}
\newtheorem{lemma}[equation]{Lemma}
\newtheorem{prop}[equation]{Proposition}
\newtheorem{cor}[equation]{Corollary}
\newtheorem{definition}[equation]{Definition}
\newtheorem{remark}[equation]{Remark}
\newtheorem{notation}[equation]{Notation}

\numberwithin{equation}{section}

\newcommand{\Grp}{\mathscr{G}}
\newcommand{\GrpK}{\Grp_K}
\newcommand{\GrpT}[1]{\Grp_{T_{#1}}}

\newcommand{\lhat}{\widehat{L}}

\newcommand{\R}{\mathbb{R}}
\newcommand{\Sph}{\mathbb{S}}
\newcommand{\T}{\mathbb{T}}
\newcommand{\C}{\mathbb{C}}
\newcommand{\Z}{\mathbb{Z}}

\newcommand{\That}{\widehat{\T}}
\newcommand{\Tbar}{\overline{\T}}

\newcommand{\geuc}{g_{_E}}
\newcommand{\gsph}{g_{_S}}

\newcommand{\rr}{\mathrm{r}}
\newcommand{\xx}{\mathrm{x}}
\newcommand{\yy}{\mathrm{y}}
\newcommand{\zz}{\mathrm{z}}

\newcommand{\uj}[1]{u_{#1}}

\newcommand{\zzj}[1]{\zz_{#1}}
\newcommand{\rrj}[1]{\rr_{#1}}

\newcommand{\tauj}[1]{\tau_{#1}}
\newcommand{\zetaj}[1]{\zeta_{#1}}
\newcommand{\xij}[1]{\xi_{#1}}
\newcommand{\cj}[1]{b_{#1}}
\newcommand{\djj}[1]{d_{#1}}

\newcommand{\taubar}{\underline{\tau}}
\newcommand{\abar}{\underline{a}}
\newcommand{\taubarj}[1]{\taubar_{#1}}

\newcommand{\Hcal}{\mathcal{H}}
\newcommand{\Lcal}{\mathcal{L}}
\newcommand{\Rcal}{\mathcal{R}}
\newcommand{\Qcal}{\mathcal{Q}}

\newcommand{\Acal}{\mathcal{A}}

\newcommand{\Lchi}{\mathcal{L}_\chi}

\newcommand{\LhatK}{\widehat{\mathcal{L}}_K}
\newcommand{\RhatK}{\widehat{\mathcal{R}}_K}
\newcommand{\RhatT}[1]{\widehat{\mathcal{R}}_{T_{#1}}}
\newcommand{\RtildeK}{\widetilde{\mathcal{R}}_K}
\newcommand{\RtildeT}[1]{\widetilde{\mathcal{R}}_{\torr[#1]}}
\newcommand{\Ltilde}{\widetilde{\Lcal}}
\newcommand{\Rtilde}{\widetilde{\Rcal}}

\newcommand{\Pcal}{\mathcal{P}}

\newcommand{\Fcal}{\mathcal{F}}
\newcommand{\Ftilde}{\widetilde{\Fcal}}
\newcommand{\Dcal}{\mathcal{D}}

\newcommand{\Jcal}{\mathcal{J}}

\newcommand{\mubar}{\underline{\mu}}
\newcommand{\lambdabar}{\underline{\lambda}}
\newcommand{\wbar}{\underline{w}}

\newcommand{\cbar}{\underline{c}}

\newcommand{\sgn}{\operatorname{sgn}}
\newcommand{\sech}{\operatorname{sech}}
\newcommand{\tr}{\operatorname{tr}}
\newcommand{\arcosh}{\operatorname{arcosh}}
\newcommand{\spt}{\operatorname{spt}}

\newcommand{\cyl}{\mathbb{K}}
\newcommand{\cate}{\widehat{\kappa}}
\newcommand{\cat}{\kappa}

\newcommand{\catr}{\mathcal{K}}
\newcommand{\torr}{\mathcal{T}}

\newcommand{\rhohat}{\widehat{\rho}}

\newcommand{\chihat}{\widehat{\chi}}
\newcommand{\chihatK}{\chihat_{_K}}

\newcommand{\rot}{\mathsf{R}}
\newcommand{\xbar}{\underline{\mathsf{X}}}
\newcommand{\ybar}{\underline{\mathsf{Y}}}
\newcommand{\zbar}{\underline{\mathsf{Z}}}

\newcommand{\transe}{\widehat{\mathsf{T}}}
\newcommand{\xbarhat}{\widehat{\xbar}}
\newcommand{\ybarhat}{\widehat{\ybar}}
\newcommand{\zbarhat}{\widehat{\zbar}}

\newcommand{\abs}[1]{\left\lvert#1\right\rvert}
\newcommand{\norm}[1]{\left\|#1\right\|}
\newcommand{\cutoff}[2]{\psi\left[ #1,#2 \right]}

\title{Minimal surfaces in the $3$-sphere by stacking Clifford tori}
\author{David Wiygul}
\address{Department of Mathematics and Statistics, California State University, Long Beach, CA 90840}
\email{david.wiygul@csulb.edu}

\begin{document}

\begin{abstract}
Extending work of Kapouleas and Yang,
for any integers $N \geq 2$, $k, \ell \geq 1$,
and $m$ sufficiently large,
we apply gluing methods to construct in the round $3$-sphere
a closed embedded minimal surface
that has genus $k\ell m^2(N-1)+1$
and is invariant under a $D_{km} \times D_{\ell m}$ subgroup of $O(4)$,
where $D_n$ is the dihedral group of order $2n$.
Each such surface resembles the union of $N$ nested topological tori,
all small perturbations of a single Clifford torus $\mathbb{T}$,
that have been connected by
$k\ell m^2 (N-1)$ small catenoidal tunnels,
with $k \ell m^2$ tunnels joining each pair of neighboring tori.
In the large-$m$ limit for fixed $N$, $k$, and $\ell$,
the corresponding surfaces converge to $\mathbb{T}$ counted with multiplicity $N$.
\end{abstract}

\maketitle

\section{Introduction}
\label{intro}

In \cite{KY} Kapouleas and Yang constructed a sequence of embedded minimal
surfaces in the round $3$-sphere $\Sph^3$
converging to a Clifford torus $\T$ counted with multiplicity $2$;
each surface consists of two small perturbations of $\T$
connected by many catenoidal annuli
taking their centers at the sites of a square lattice on $\T$.
Accordingly they called their surfaces \emph{doublings} of the Clifford torus.
Kapouleas announced these in \cite{KapClay} as the first examples
of a general class of gluing constructions to double given minimal surfaces,
subsequently discussed further in \cite{KapRS}.
More recently in \cite{KapSphDbl}
he has doubled
the equatorial $2$-sphere in $\Sph^3$,
and now additional such doublings with different configurations of catenoidal tunnels
have been carried out by Kapouleas and McGrath \cite{KapMcG}.
Min-max methods have also been used to double minimal surfaces in $\Sph^3$.
Pitts and Rubinstein proposed a variety of such constructions in \cite{PR}.
One was completed by Ketover, Marques, and Neves in \cite{KMN},
where they too double the torus over square lattices,
conjecturally producing the same surfaces as \cite{KY}
when the lattice spacing becomes small,
and in \cite{Ket} Ketover has performed
more min-max constructions, including doublings,
previously described in \cite{PR}.
Doublings appear in the free-boundary setting as well.
Using variational rather than gluing methods,
for each integer $n \geq 3$ Fraser and Schoen (\cite{FS})
have constructed orientable free-boundary minimal embeddings in the unit ball
with genus $0$ and $n$ boundary components;
for large $n$ these surfaces look like doublings of the equatorial disc.
Later, in \cite{FPZ}, Folha, Pacard, and Zolotareva applied gluing techniques to
double the equatorial disc,
producing free-boundary examples with genus $0$
(possibly the same as those in \cite{FS})
or $1$ and a large number $n$ of boundary components.

Returning to \cite{KY}, the surfaces of Kapouleas and Yang are highly symmetric,
admitting many \emph{horizontal} symmetries, which preserve
as sets each of the two sides of the doubled Clifford torus and permute
the lattice sites, as well as \emph{vertical} symmetries,
each of which exchanges the two sides of the doubled torus
but fixes as a set a catenoidal tunnel
(and in fact every catenoidal tunnel centered on a certain great circle on $\T$).
All these symmetries are enforced throughout the construction
and exploited to simplify its execution.
The present article undertakes less
vertically symmetric doublings of the torus, with the symmetry broken in two ways.
First, we allow the catenoidal tunnels to be arranged on rectangular
rather than strictly square lattices. Any isometry of $\Sph^3$ exchanging
the two sides of the doubled Clifford torus will fail to preserve such
a lattice, unless it is square.
Second, we interpret \emph{doubling}
in a generalized sense, realizing also triplings, quadruplings,
and in fact embedded minimal surfaces resembling any prescribed finite number of
slightly perturbed copies of $\T$
connected to one another by many small catenoidal tunnels.
Whenever at least three copies are incorporated,
even if these tunnels are centered on square lattices,
the symmetry group will not act transitively on the collection of copies.

These new constructions add to the list of known closed minimal embeddings in $\Sph^3$,
so far comprising those found in
\cite{Law}, \cite{KPS}, \cite{KY},
\cite{CS}, \cite{KapSphDbl}, \cite{MNposric},
\cite{KMN}, \cite{Ket}, \cite{KapMcG}, \cite{KWtordesing},
and \cite{BWW}.
The survey article \cite{BreS3}
contains an outline of a few of the constructions just mentioned.
The constructions at hand should be of interest not only
as providing new examples in $\Sph^3$
but also as a basis for further doublings with
asymmetric sides in a variety of settings.
A program toward doubling constructions of increasing generality,
including potential applications, is described in \cite{KapRS}.
The present work naturally emulates, with a few departures,
the approach of \cite{KY}
and draws extensively from the general gluing technology developed
by Kapouleas, much of which can be found summarized
in \cite{KapClay}
and was itself inspired by techniques applied by Schoen in \cite{SchPC}.
Although the current article can be read without reference to \cite{KY}
or any other gluing constructions, for the rest of this introduction
we will make use, without detailed explanation, of terminology standardized
by Kapouleas, so that the reader already acquainted with it may easily
appreciate the principal differences between this construction and \cite{KY}.

We now outline our procedure in rough terms.
We first fix a Clifford torus $\T$, which by definition
is the locus of points in $\Sph^3 \subset \R^4$
at distance $\frac{\pi}{4}$ from some great circle $C_1$.
More generally, for each $r \in \left(0,\frac{\pi}{2}\right)$
the locus of points at distance $r$ from $C_1$
is a torus of constant mean curvature,
which is equivalently the locus of points
at distance $\frac{\pi}{2}-r$ from the great circle $C_2$
defined as the intersection with $\Sph^3$ of the orthogonal complement in $\R^4$
of the plane containing $C_1$.
We will refer to $C_1$ and $C_2$ as the \emph{axes} of each such torus;
directions tangential to any of these tori we will call \emph{horizontal},
while the orthogonal direction we will call \emph{vertical}.
As basic data for the construction we take integers $k,\ell,m\geq1$ and $N \geq 2$.
Corresponding to a choice of such data
an \emph{initial surface} is built as follows.

We start with $N$ constant-mean-curvature tori coaxial with $\T$,
labeled $\T[1],\T[2],\ldots,\T[N]$
so that the index increases with distance from $C_2$.
The precise placements of the tori (that is their signed distances from $\T$)
cannot be freely prescribed but
will be determined, as a function of the data, in the course of the construction.
For now we mention only that as $m \to \infty$
every $\T[j]$ tends to $\T$.
Next we will connect this collection of tori
by first from each one excising discs centered on certain lattices
(to be described momentarily)
and then gluing in truncated approximate catenoids,
which shrink to points on $\T$ as $m \to \infty$.
Like the arrangement of the tori,
the precise sizes of the catenoids
and the heights (signed distances from $\T$) of their centers
are variables whose values will be set
by conditions---to be described later---necessary for the completion of the construction.

On the other hand,
we impose enough horizontal symmetries
(isometries of $\Sph^3$ preserving each side of $\T$)
that the horizontal positioning of the catenoids
and equivalently the locations of the discs deleted from the tori
are directly determined by the data already listed.
Specifically, we mark $km$ equally spaced points on $C_1$
and $\ell m$ equally spaced points on $C_2$,
and we write $\Grp[k,\ell,m]$ for the subgroup of $O(4)$
preserving each of these sets of marked points.
Thus $\Grp[k,\ell,m]$ is isomorphic to $D_{km} \times D_{\ell m}$,
where $D_n$ is the dihedral group of order $2n$.
Note that $\Grp[k,\ell,m]$ is equally the subgroup of $O(4)$
preserving the union of the sets of marked points
except when $k=\ell$,
in which case $\Grp[k,\ell,m]$ is strictly contained in this last group,
which admits also reflections through certain great circles on $\T$,
exchanging $C_1$ and $C_2$.

We choose orientations on $C_1$ and $C_2$
(together selecting an orientation on $\R^4$),
and for each $i \in \{0,1\}$
we write $\rot_{C_i}^\theta$
for the element of $O(4)$
fixing $C_i$ pointwise
and rotating the great circle in the orthogonally complementary plane
through angle $\theta$ (according to its orientation).
The group
$\Grp[k,\ell,m]$ is then generated by
(i) $\rot_{C_2}^{2\pi/km}$,
(ii) $\rot_{C_1}^{2\pi/\ell m}$,
(iii) reflection through any great sphere having equator $C_2$
and one pole a marked point on $C_1$,
and (iv) reflection through any great sphere having equator $C_1$
and one pole a marked point on $C_2$.
(Of course the point antipodal to a given marked point on $C_1$ (or $C_2$)
is itself marked if and only if $km$ (or $\ell m$) is even.)
Obviously $\Grp[k,\ell,m]$ preserves $\T$ and each $\T[j]$.
We will design the initial surfaces
(and the final minimal surfaces they approximate)
to be likewise invariant under $\Grp[k,\ell,m]$.

Next we pick a marked point on $C_1$, a marked point on $C_2$,
and the minimizing geodesic segment (quarter great circle) joining them.
This segment intersects $\T$ at a single point,
whose orbit under $\Grp[k,\ell,m]$
is a $km \times \ell m$ rectangular lattice on $\T$, which we call $L_{0,0}$.
We will use $L_{0,0}$ to fix the horizontal positions of the catenoidal
annuli connecting the $N$ tori in our configuration.
In fact there are precisely four $km \times \ell m$ rectangular lattices on $\T$
preserved by $\Grp[k,\ell,m]$.
It would be possible to carry out the constructions in this paper
using any of these lattices (to locate the catenoidal annuli)
without introducing any additional technical difficulties,
but to simplify the presentation we will make use of only
$L_{0,0}$
and $L_{1,1}:=\rot_{C_2}^{\frac{\pi}{km}}\rot_{C_1}^{\frac{\pi}{\ell m}}L_{0,0}$.
(Of course there are also finer lattices preserved by the same group.
Permitting such lattices in the construction of the initial surface
would allow for different numbers of catenoidal tunnels connecting
different pairs of tori while maintaining the high horizontal symmetry
but would demand a more complicated approach.)

Now for every $j \in \{1,2,\ldots N-1\}$
and for every point in $L_{(j+1) \bmod 2, (j+1) \bmod 2}$
we locate the closest point on $\T[j]$
and the closest point on $\T[j+1]$,
and we excise from the two tori two small discs
having these nearest points as their respective centers.
Then for each such pair of points,
using local coordinates for $\Sph^3$ adapted to the tori,
we smoothly glue the boundary circles of the deleted discs
to the boundary circles of a catenoidal annulus
centered on the geodesic segment connecting the two points.
The radii of the deleted discs are chosen comparable to the lattice spacing
but small enough so that all the discs are pairwise disjoint,
and the annuli are shaped so that
the resulting connected surface is invariant under $\Grp[k,\ell,m]$.
As already noted, additional information is needed to specify
the precise sizes and heights of these annuli,
but right now we mention that, when suitably scaled, each tends with large $m$
to a complete standard catenoid.

Thus we have produced a connected closed surface, the initial surface,
which is preserved by $\Grp[k,\ell,m]$,
is easily seen to have genus
$N + (N-1)k\ell m^2 - (N-1)=k\ell m^2(N-1)+1$
(since $N-1$ of the $(N-1)k\ell m^2$ catenoidal annuli
are spent to connect the $N$ tori, each of genus $1$,
while the remaining ones contribute genus),
and, for large $m$, is approximately minimal in a certain sense.
The construction will be completed by perturbing the
surface to exact minimality.
Two mechanisms of perturbation are applied in tandem.
One sort of perturbation is realized by considering graphs of small
functions over the initial surface.
To select the right function is then to solve
the elliptic quasilinear partial
differential equation prescribing zero mean curvature for the corresponding graph.
This equation can be studied by comparing the linearization of the operator
governing the mean curvature of graphs to certain large-$m$ limit operators on
the limit catenoids and limit torus. In the simplest scenario imaginable one
could solve the linearized equation on the toral and catenoidal components separately,
combine these solutions through an iterative procedure, and finally invoke
an inverse function theorem to solve the original nonlinear equation.
However, the presence of nontrivial kernel to the limit operators
gives rise to \emph{approximate kernel} that obstructs the approach just described.

The space of admissible perturbing functions is constrained to respect the
symmetries enjoyed by the initial surface, and so their imposition has
the effect of reducing the dimension of the approximate kernel.
Each torus turns out to carry one-dimensional approximate kernel of its own,
but in \cite{KY} the two tori can be exchanged
by reflections through certain great circles, and so together the tori
contribute just one dimension to the approximate kernel in \cite{KY} versus
$N$ dimensions more generally. Furthermore, in \cite{KY} these reflections
through circles render trivial the approximate kernel on the catenoidal tunnels.
Following the approach of \cite{KY} in the absence of these symmetries,
each tunnel would feature one-dimensional approximate kernel,
but we bypass this kernel altogether by altering,
as compared to \cite{KY}, the initial data
at the tunnel's waist for the rotationally invariant mode of the solution.

To overcome the obstruction posed by the approximate kernel
\cite{KY} introduces \emph{substitute kernel}, spanned by
a single function supported on the tori away from the circles where they
attach to the tunnels. By adding multiples of this function to the source
term of the linearized equation, the so modified source can be made
orthogonal to the approximate kernel, enabling the success of the above
scheme, but at the cost of solving the original equation only modulo
substitute kernel. For the same purpose the current construction introduces
$N$-dimensional substitute kernel, spanned by functions each of which
is supported on a single torus away from the tunnels.
(Actually, in this construction we never explicitly identify the approximate kernel,
nor do we invoke the $h$ metric employed in \cite{KY} for its analysis,
but our application of substitute kernel is morally identical.)

A further difficulty concerns the vast disparity in scale
between the waist radii of the
catenoidal tunnels on the one hand and the much greater
spacing between the tori on
the other.
The norm of the initial surface's second fundamental form
grows toward the waists of the catenoids
to a value diverging with $m$
from a value bounded uniformly in $m$ on the tori, and the embeddedness
of graphical perturbations is most precarious near the waists.
For these reasons,
as well as to ensure convergence of the iteratively defined global solution,
it is necessary to arrange for solutions on the tori
to decay toward the catenoidal waists.

All of the catenoids attaching to each of the two outermost
tori---the only type of torus appearing in \cite{KY}---are equivalent modulo the symmetries,
and adjustment of the source term
by the substitute kernel suffices to achieve such decay on these catenoids.
(Again, our actual approach deviates somewhat from this description,
applicable to \cite{KY}, but just superficially.)
However, each of the intermediate tori, $N-2$ in number,
attaches to catenoids of precisely two equivalence classes under $\Grp[k,\ell,m]$,
and so the appropriate decay of solutions requires the introduction
of another $N-2$ functions, linear combinations of which
are added to the source term to arrange decay,
a device originating in \cite{KapWenT} but unneeded in \cite{KY}.
In total we arrive at a $(2N-2)$-dimensional \emph{extended substitute kernel},
the sum of the substitute kernel and the span of these additional functions,
modulo which subspace we can, for large $m$, invert the linearized operator.

Thus an infinite-dimensional problem is reduced to a finite-dimensional one.
The resolution of this latter problem requires the second type
of perturbation and is best understood in terms of a correspondence,
which Kapouleas (\cite{KapWenT}) calls the \emph{geometric principle}, between
the initial geometry and the analytic obstructions that the
extended substitute kernel represents. In a few words, elements of the extended
substitute kernel can be generated, as components
of the initial surface's mean curvature,
by certain motions of its building blocks---here
catenoids and tori---relative to one another.
In accordance with this principle the other type of perturbation is realized by incorporating
parameters, one for each dimension of extended substitute kernel,
into the definition of the initial surface, whose
variation repositions the component tori and catenoids.
Thus for each choice of $k$, $\ell$, $m$, and $N$
we define not one initial surface but a $(2N-2)$-parameter family of initial surfaces.
More specifically, two parameters may be associated with each of the
$N-1$ classes of catenoids joining pairs of adjacent tori.
One set of parameters, $\{\zetaj{i}\}_{i=1}^{N-1}$, controls the waist radii,
while the other set, $\{\xij{i}\}_{i=1}^{N-1}$, adjusts the heights of the centers.
A degree of rigidity,
in the form of \emph{matching conditions},
is maintained to reposition the tori in response to the parameters,
and the surface is smoothed using cutoff functions as needed.
A single parameter $\zeta$ works for \cite{KY}, since
there $N=2$ and the symmetry between the sides of $\T$ forces $\xi=0$.

In the course of the construction
it is necessary to solve for the proper parameter values
along with the perturbing function.
The parameter dependence of the ``extended'' components
of the extended substitute kernel can be directly estimated with
accuracy adequate for our purposes. It turns out that these components
are primarily generated by \emph{dislocations} resulting from antisymmetric
variation in pairs of $\xi$ parameters associated to catenoids adjoining a common torus.
The parameter dependence of the
substitute kernel itself is more conveniently monitored indirectly,
as in \cite{KY}, via \emph{forces}.
On each torus
the elements of the approximate kernel, and so of the substitute kernel,
may be identified with approximate translations of the torus
relative to $\T$.
In fact $\Sph^3$ admits an exact Killing field which,
though it does not exactly generate this variation of the torus, does approximate
it in the vicinity of a great circle orthogonally intersecting $\T$.
The force in the direction of this Killing field through certain neighborhoods
of a given torus then serves as an estimate of the projection of the mean curvature
onto the approximate kernel and thereby as a proxy for the corresponding
component of substitute kernel itself.
The \emph{balancing} equations and
the analysis of the parameter dependence of the forces here are substantially more
complicated than those of \cite{KY} but no different in principle.

Finally, estimates for the initial geometry, the linearized equation,
the nonlinear terms, and the parameter dependence of the forces and
dislocations are applied in conjunction with the Schauder fixed-point theorem
to prove our main result,
which we state informally now,
a more refined version appearing as Theorem \ref{mainthm},
which makes use of notation developed throughout in the paper.

\begin{theorem}[Informal statement of the main theorem]
Let $k, \ell \geq 1$, and $N \geq 2$ be given integers.
For sufficiently large $m$ there exist
both a choice of parameters and a smooth, appropriately symmetric
perturbing function such that the resulting
surface (as described above)
is minimal, invariant under $\Grp[k,\ell,m]$,
and a small perturbation of the corresponding initial surface,
so in particular embedded and of genus $k\ell m^2(N-1)+1$.
\end{theorem}

\begin{remark}[The full symmetry group]
\label{verticalsymremark}
As described earlier,
the catenoidal tunnels joining a pair of adjacent tori in a given initial surface
take their centers on geodesic arcs intersecting $\T$ at the sites of
a $km \times \ell m$ rectangular lattice invariant under $\Grp[k,\ell,m]$.
Because the minimal surfaces produced by the construction
are small perturbations of the initial surfaces,
the symmetry group of each resulting minimal surface---that is the subgroup of $O(4)$ preserving
it as a set---must preserve each of these lattices.
When $k \neq \ell$, we can therefore conclude that this symmetry group
does not merely contain $\Grp[k,\ell,m]$ but coincides with it.
On the other hand, when $k=\ell$,
there are additional isometries of $\Sph^3$, not belonging to $\Grp[k,\ell,m]$
(namely vertical ones, exchanging the sides of $\T$),
that preserve each lattice,
which one could easily enforce in the construction to obtain minimal surfaces
enjoying these extra symmetries as well.
To avoid complicating the presentation, however,
we do not carry out these details.
Without making such modifications
it is not immediately clear
whether or not the minimal surfaces
resulting from our construction in the $k=\ell$ case
necessarily possess vertical symmetries;
to prove they do it would suffice to establish that the fixed point
in the proof of Theorem \ref{mainthm} is unique.
\end{remark}

\begin{remark}[Choice of lattices]
\label{latticeremark}
We have also already mentioned
that there are in fact four $km \times \ell m$ rectangular lattices on $\T$
invariant under $\Grp[k,\ell,m]$,
but the constructions in the present article
utilize only two (and of course just one in the special case that $N=2$)
in alternating fashion
to distribute the catenoidal tunnels joining each pair of adjacent tori.
It would be possible,
without incurring any technical difficulties that we do not already confront,
to avail ourselves of any of the four lattices
when fixing the horizontal locations of the catenoidal annuli
(subject only to the obvious constraint that two lattices corresponding
to consecutive pairs of tori be distinct).
By doing so, for $N>2$, we could construct a variety of examples
not congruent to one another but having the same genus and symmetry group.
This same flexibility would also allow us to construct
examples with $k=\ell$ and $N>2$ which indubitably do not enjoy any vertical symmetries.
(See Remark \ref{verticalsymremark}, just above.)
To avoid complicating the definitions in this article any further
we do not present our construction in this generality.
In a more ambitious construction one could even attempt to allow
a different (large) number of catenoidal annuli at each layer,
but this modification would require a genuinely more elaborate approach.
\end{remark}

\subsection*{Outline of the presentation}
In Section 2 we define the initial surfaces.
In Section 3 we analyze the dependence on the $\zeta$ and $\xi$ parameters of the dislocations
and vertical forces through various regions.
In Section 4 we obtain estimates for the geometry of the initial surfaces.
In Section 5 we study the linearized problem.
In Section 6 we solve the nonlinear problem, proving the main theorem.

\subsection*{Acknowledgments}
This article presents the results of my PhD thesis,
supervised by Nicos Kapouleas, who suggested the problems studied therein;
I am deeply grateful for his guidance.
Additional thanks are owed to Scott Field
for discussions that further stimulated my interest in $N>2$ stacking,
to Christina Danton for help preparing the figures,
to Hung Tran for much appreciated feedback on a portion of the paper,
and to Rick Schoen for essential advice
throughout the process of preparing the article.

\section{Initial surfaces}
In this section we must make a number of preliminary definitions
before defining the initial surfaces themselves.
We also try to offer some motivation for these definitions.
The eager reader may wish to look ahead to \eqref{initsurfdef}
and the references immediately preceding it,
consulting the intervening material only as needed.

We realize $\Sph^3$
as the unit sphere $\{(z_1,z_2) \in \C^2 \; : \; \abs{z_1}^2+\abs{z_2}^2=1\}$
in $\C^2$ and set
  \begin{equation}
  \label{Cliffordtorusdef}
    \T := \left\{ (z_1,z_2) : \abs{z_1}=\abs{z_2} = \frac{1}{\sqrt{2}} \right\},
  \end{equation}
the Clifford torus whose axes (as defined in Section 1)
are simply
the coordinate unit circles $C_1 := \{ z_2 = 0 \}$ and $C_2 := \{ z_1 = 0 \}$.
We define the covering map
  \begin{equation}
  \label{Phidef}
    \begin{aligned}
      &\Phi:
      \R \times \R \times \left(-\frac{\pi}{4},\frac{\pi}{4}\right)
        \to \Sph^3 \backslash (C_1 \cup C_2) \mbox{ by} \\
      &\Phi(\xx,\yy,\zz)
      :=
      \left( e^{i\sqrt{2}\xx} \sin \left(\zz + \frac{\pi}{4}\right),
        e^{i\sqrt{2}\yy} \cos \left(\zz + \frac{\pi}{4}\right)\right),
    \end{aligned}
  \end{equation}
which maps (i) horizontal planes to constant-mean-curvature tori
having axes $C_1$ and $C_2$, with $\Phi(\{\zz=0\})=\T$ in particular,
(ii) vertical lines to quarter great circles orthogonal to $C_1$, $C_2$, and $\T$,
(iii) vertical planes of constant $\xx$ to great hemispheres with equator $C_2$,
(iv) vertical planes of constant $\yy$ to great hemispheres with equator $C_1$,
and
(v) vertical planes of constant $\xx \pm \yy$ to half Clifford tori
through $C_1$ and $C_2$, orthogonally intersecting $\T$ along great circles.
Writing $\gsph$ for the standard round metric on $\Sph^3$
and $\geuc$ for the standard flat metric on $\R^3$,
we find
  \begin{equation}
  \label{Phullback}
    \Phi^*\gsph = \geuc + (\sin 2\zz)\left(d\xx^2 - d\yy^2\right).
  \end{equation}
The initial surfaces will be assembled by applying $\Phi$ to a stack
of horizontal planes connected by staggered catenoidal columns.

  \begin{figure}
    \includegraphics[width=4.5in]{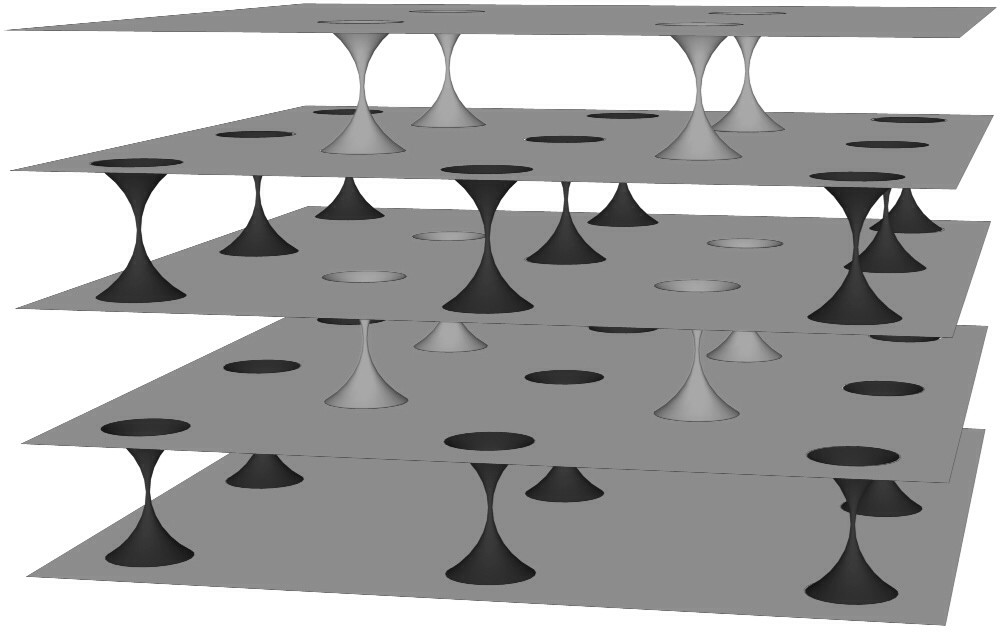}
      \caption{An $N=5$ initial surface, cut by four hemispheres
      and viewed in the $(\xx,\yy,\zz)$ coordinate system; for clarity
      the picture grossly exaggerates the ratio,
      which in fact tends to $0$ with large $m$, of the vertical
      spacing between the tori
      to the horizontal spacing between the catenoids.}
  \end{figure}

\subsection*{Half-catenoids bent to planes}
We shall make frequent use of cutoff functions
throughout the construction,
so we fix now a smooth, nondecreasing $\Psi: \R \to [0,1]$ with $\Psi$  identically $0$
on $(-\infty,-1]$, identically $1$ on $[1,\infty)$,
and such that $\Psi - \frac{1}{2}$ is odd.
We then define, for any $a,b \in \R$, the function $\cutoff{a}{b}: \R \to [0,1]$ by
  \begin{equation}
  \label{cutoff}
    \cutoff{a}{b} := \Psi \circ L_{a,b},
  \end{equation}
where $L_{a,b}: \R \to \R$ is the linear function satisfying $L(a)=-3$ and $L(b)=3$.

We write $B[(x,y),r]$ for the open Euclidean disc in $\R^2$ with radius $r$ and center $(x,y)$.
Given $X,Y,\tau>0$ with $\tau<\min \{X,Y\}$,
we set
  \begin{equation}
  \label{Tone}
    \T_{X,Y,\tau} := \left([-X,X] \times [-Y,Y]\right) \backslash B[(0,0),\tau],
  \end{equation}
a solid rectangle with a disc removed from its center.
Given also $\zzj{K},\zzj{T} \in \R$ and $R>\tau$ with $2R<\min \{X,Y\}$,
we define the function
$\phi[\zzj{K},\zzj{T},R,X,Y,\tau]: \T_{X,Y,\tau} \to \R$
by
  \begin{equation}
    \begin{aligned}
      \phi[\zzj{K},\zzj{T},R,X,Y,\tau](\xx,\yy) :=
      & \zzj{K} + (\zzj{T} -\zzj{K}) \cutoff{R}{2R}\left(\sqrt{\xx^2+\yy^2}\right) \\
      & + \sgn(\zzj{T} - \zzj{K}) \left(\tau \arcosh \frac{\sqrt{\xx^2+\yy^2}}{\tau}\right)
        \cutoff{2R}{R}\left(\sqrt{\xx^2+\yy^2}\right),
    \end{aligned}
  \end{equation}
where 
$\arcosh: [1,\infty) \to [0,\infty)$
is the inverse of
the restriction to $[0,\infty)$ of the hyperbolic cosine function
and the sign function $\sgn: \R \to \R$
takes the value $1$ when its argument is nonnegative and the value $-1$ otherwise.
Thus, inside the cylinder of radius $R$ about the $\zz$-axis
the graph of $\phi$ coincides with the portion
between the planes $\zz=\zzj{K}$ and $\zz=\zzj{T}$
of the catenoid
with vertical axis, center $\left(0,0,\zzj{K}\right)$, and waist radius $\tau$,
while outside the cylinder of radius $2R$ about the $\zz$-axis
it coincides with the solid rectangle $[-X,X] \times [-Y,Y] \times \left\{\zz=\zzj{T}\right\}$.
In between these two cylinders the cutoff function is used to bend the end of the half-catenoid
to become exactly horizontal.
With the additional data $(\xx_{0},\yy_{0}) \in \R^2$
we define the embedding
  \begin{equation}
  \label{Text}
    \begin{aligned}
      &T_{ext}[(\xx_{0},\yy_{0}),\zzj{K},\zzj{T},R,X,Y,\tau]
      : \T_{X,Y,\tau} \to \R^3 \mbox{ by} \\
      &T_{ext}(\xx,\yy)
      :=
      \left(\xx_{0} + \xx, \yy_{0} + \yy, \phi[\zzj{K},\zzj{T},R,X,Y,\tau](\xx,\yy)\right),
    \end{aligned}
  \end{equation}
whose image is the graph of $\phi$ translated by $(\xx_{0},\yy_{0},0)$.

Given also $\zzj{K}' \in \R$ along with $\tau'>0$,
now assuming $\max \{\tau,\tau'\}<R<\frac{1}{4}\min \{X,Y\}$,
we set
  \begin{equation}
  \label{Ttwo}
    \T_{X,Y,\tau,\tau'} := ([-X,X] \times [-Y,Y]) \left\backslash
      \left( B\left[\left(-\frac{X}{2},-\frac{Y}{2}\right),\tau\right]
        \cup  B\left[\left(\frac{X}{2},\frac{Y}{2}\right),\tau'\right]
        \right)\right.,
  \end{equation}
and we define the embedding
  \begin{equation}
  \label{Tint}
    \begin{aligned}
      &T_{int}[(\xx_{0},\yy{_0}),\zzj{K},\zzj{K}',\zzj{T},R,X,Y,\tau,\tau']
      :\T_{X,Y,\tau,\tau'} \to \R^3 \mbox{ by} \\
      &T_{int}(\xx,\yy) := \\
      &\;\;  \begin{cases}
          \left(\xx_{0}+\xx,\yy_{0}+\yy,
            \phi\left[\zzj{K},\zzj{T},R,\frac{X}{2},\frac{Y}{2},\tau\right]
            \left(\xx+\frac{X}{2},\yy+\frac{Y}{2}\right)
          \right)
            \mbox{ for } (\xx,\yy) \in [-X,0] \times [-Y,0] \\
          \left(\xx_{0}+\xx,\yy_{0}+\yy,
            \phi\left[\zzj{K}',\zzj{T},R,\frac{X}{2},\frac{Y}{2},\tau'\right]
            \left(\xx-\frac{X}{2},\yy-\frac{Y}{2}\right)
          \right)
            \mbox{ for } (\xx,\yy) \in [0,X] \times [0,Y] \\
          (\xx_{0}+\xx,\yy_{0}+\yy,\zzj{T}) \mbox{ everywhere else. }
        \end{cases}
    \end{aligned}
  \end{equation}
Thus the image of $T_{int}$ looks like a solid $2X \times 2Y$ rectangle in the $\zz=\zzj{T}$ plane
with two discs centered at
$\left(\xx_{0},\yy_{0},\zzj{T}\right) \pm \left(\frac{X}{2},\frac{Y}{2},0\right)$
replaced by catenoidal annuli terminating on waist circles at heights $\zzj{K}'$
and $\zzj{K}$.

The initial surfaces will be built from various applications of $T_{ext}$,
for extreme or outermost tori and adjoining half-catenoids, and of $T_{int}$,
for the intermediate tori and pairs of adjoining half-catenoids.
The horizontal positions of the catenoids
(the values of $(\xx_{0},\yy_{0})$ in the parametrizations above)
and the dimensions ($X$ and $Y$) of the parametrizing solid rectangles
(or equivalently the lattice edges)
will be set directly by the data $k$, $\ell$, and $m$.
The radii of the annuli of transition (determined by $R$)
will be chosen on the order of $\min\{X,Y\}$
(but smaller than it by a wide enough margin that the images
of $T_{ext}$ and $T_{int}$ are horizontal near their boundaries as assumed above).
There remain $N-1$ selections of waist radii ($\tau$ and $\tau'$)
and $N-1$ more of waist heights ($\zzj{K}$ and $\zzj{K}'$).
\emph{Balancing conditions} studied in the next section allow
for the estimation of the values these $2N-2$ unknowns must assume
for the construction to succeed, but their precise specification is made
by the $\zeta$ and $\xi$ parameters described in Section 1.
\emph{Matching conditions} will then fix the height $\zzj{T}$ of each torus,
by requiring these heights to agree with the heights of the adjoining catenoids where they
meet the transition annuli, in the case of the extreme tori, or to agree
with the average of the heights of the upper and lower catenoids
at the transition circles they adjoin, in the case of the intermediate tori.

\subsection*{The hierarchy of data}
For any positive integers $k$, $\ell$, and $N \geq 2$
the first phase of the construction produces a sequence,
indexed by $m$, of $(2N-2)$-parameter families of initial surfaces.
In order to obtain adequate estimates for the initial mean curvature,
the linearized operator, and the nonlinear terms,
we will routinely make the assumption that $m$
is as large as needed in terms of $k$, $\ell$, $N$, and all of the parameters.
Of course we expect the ultimate parameter choices themselves to depend on $m$,
so, in taking $m$ large as just described, it is necessary to assume
that the parameters are all bounded in absolute value
by a constant $c>0$ independent of $m$.
Naturally we do not yet know
what range is needed for the parameters, but eventually we will
be able to pick $c$ in terms of $k$, $\ell$, and $N$
so that for every sufficiently large $m$
we will be able to find parameters bounded by $c$
so that the corresponding initial surface can be perturbed to minimality.

To continue with the definition of the intial surfaces
we assume we are given
integers $k,\ell,m \geq 1$ and $N \geq 2$
as well as
a constant $c>0$ and parameters $\zeta, \xi \in [-c,c]^{N-1}$.
For notational simplicity we assume
  \begin{equation}
    k \leq \ell
  \end{equation}
and we write $n=n[N]$ for the greatest integer no greater than $N/2$, so that
  \begin{equation}
  \label{ndef}
    N=\begin{cases}2n \mbox{ when $N$ is even} \\ 2n+1 \mbox{ when $N$ is odd.}\end{cases}
  \end{equation}
We acknowledge a certain redundancy in the minimal surfaces ultimately exhibited,
one which is easily removed by taking $k$ and $\ell$ relatively prime.

\subsection*{Catenoidal radii and vertical specifications}
Modulo the symmetries that we will impose
when defining our initial surfaces
we have $N-1$ catenoidal waist radii
$\tauj{1},\ldots,\tauj{N-1}$ to prescribe,
one for each pair of adjacent tori to be joined.
Their selection is critical to the success of the construction,
and the next section (Section 3) is devoted in part to making a viable choice.
Specifically, in Lemma \ref{cilemma}, when $N \geq 4$,
we will determine a collection
$\{\cj{j}\}_{j=2}^n$ of $n-1$ positive real numbers
(recalling \eqref{ndef} just above),
which in turn we will use in conjunction with the $\zeta$ parameters
to set the waist radii of our catenoidal tunnels.
Each $\cj{j}$ will be a function of $k$, $\ell$, $N$, and $m$
but (for large $m$) will have an upper bound and a positive lower bound depending only on $N$;
in particular $\cj{2}$ will have an upper bound depending on $N$ but always strictly less than $2$.
Of course in the simpler cases $N=2$ and $N=3$ we have $n-1=0$,
but it is nevertheless notationally convenient to set
$\cj{2}:=0$ when $N=2$ and $\cj{2}:=1$ when $N=3$.
We emphasize that each $\cj{j}=\cj{j}[N,k,\ell,m]$ depends at least on $N$
and generally on the data $k$, $\ell$, and $m$ as well,
but for brevity in our notation
we will frequently suppress the expression of this dependence,
as we do for many other quantities of interest.
To summarize:
  \begin{equation}
  \label{cisummary}
    \begin{aligned}
      &\cj{2}[N=2,k,\ell,m]:=0,
      \qquad \cj{2}[N=3,k,\ell,m]:=1,
      \qquad \mbox{and for $N \geq 4$} \\
      &\{\cj{i}[N,k,\ell,m]\}_{i=2}^n
        \mbox{ is determined in Lemma \ref{cilemma} (recalling \eqref{ndef})}, \\
      &\mbox{so there exists a constant $C[N,k,\ell]>0$ such that for $N \geq 4$} \\
      &1<\max \{\cj{i}\}_{i=2}^n<C[N,k,\ell]
        \quad \mbox{and} \quad \cj{2} \leq 2-1/C[N,k,\ell] \\
      &\mbox{whenever $m$ is large enough in terms of $N$}.
    \end{aligned}
  \end{equation}

Having identified these numbers, we first define the collection
$\{\taubarj{i}\}_{i=1}^{N-1}$
of waist radii when $\zeta=0$ by
  \begin{equation}
  \label{taubardef}
    \taubarj{i}
    =
    \taubarj{i}[N,k,\ell,m]
    :=
    \begin{cases}
      \frac{1}{10 \ell m}e^{-\frac{k \ell m^2}{4\pi}\left(1-\frac{1}{2}\cj{2}[N,k,\ell,m]\right)}
        \mbox{ for } i=1 \\
      \cj{i}[N,k,\ell,m]\taubarj{1}[N,k,\ell,m] \mbox{ for } 2 \leq i \leq n \\
      \taubarj{N-i}[N,k,\ell,m] \mbox{ for } n+1 \leq i \leq N-1
    \end{cases}
  \end{equation}
(recalling \eqref{ndef})
and then for general $\zeta$ we define the radii
  \begin{equation}
  \label{taudef}
    \tauj{i}
    =
    \tauj{i}[N,k,\ell,m,\zeta]
    :=
      \begin{cases}
        e^{\zetaj{1}}\taubarj{1}[N,k,\ell,m] \mbox{ for } i=1 \\
        e^{\zetaj{1} + k^{-1}\ell^{-1} m^{-2}\zetaj{i}}\taubarj{i}[N,k,\ell,m] \mbox{ for } 1 < i < N.
      \end{cases}
  \end{equation}

We next define the $N$ heights $\zzj{i}=\zzj{i}[N,k,\ell,m,\zeta,\xi]$ of the tori
by
  \begin{equation}
  \label{torheights}
    \begin{aligned}
    &\zzj{1}
    :=
        \tauj{1} \xij{1} - 2^{N \bmod 2} \tauj{n} \ln \frac{1}{10 \ell m \tauj{n}}
          - 2 \displaystyle{\sum_{j=1}^{n-1}} \tauj{j} \ln \frac{1}{10 \ell m \tauj{j}}, \\
    &\zzj{N} :=
        \tauj{N-1} \xij{N-1} - 2^{N \bmod 2} \tauj{n} \ln \frac{1}{10 \ell m \tauj{n}}
          + 2 \displaystyle{\sum_{j=n}^{N-1}}
            \tauj{j} \ln \frac{1}{10 \ell m \tauj{j}}, \mbox{ and for $1<i<N$} \\
    &\zzj{i} :=
      \frac{\tauj{i-1}\xij{i-1} + \tauj{i}\xij{i}}{2}
          - 2^{N \bmod 2} \tauj{n} \ln \frac{1}{10 \ell m \tauj{n}}
          + 2 \displaystyle{\sum_{j=1}^{i-1}}
            \tauj{j} \ln \frac{1}{10 \ell m \tauj{j}}
          - 2\displaystyle{\sum_{j=1}^{n-1}}
            \tauj{j} \ln \frac{1}{10 \ell m \tauj{j}},
    \end{aligned}
  \end{equation}
and the $N-1$ heights $\zzj{i}^K=\zzj{i}^K[N,k,\ell,m,\zeta,\xi]$ of the catenoids' centers
by
  \begin{equation}
  \label{catheights}
    \zzj{i}^K :=
      \tauj{i}\xij{i} + \tauj{i} \ln \frac{1}{10 \ell m \tauj{i}}
      + 2\displaystyle{\sum_{j=1}^{i-1}}
          \tauj{j} \ln \frac{1}{10 \ell m \tauj{j}}
      - 2\sum_{j=1}^{n-1}
          \tauj{j} \ln \frac{1}{10 \ell m \tauj{j}}
      - 2^{N \bmod 2} \tauj{n} \ln \frac{1}{10 \ell m \tauj{n}}
  \end{equation}
for $1 \leq i \leq N-1$.

Equivalently,
  \begin{equation}
    \begin{aligned}
      &\zzj{i}^K[N,k,\ell,m,\zeta,\xi]
      =
      \zzj{i}^K[N,k,\ell,m,\zeta,0]+\tauj{i}\xij{i}
        \mbox{ for $1 \leq i \leq N-1$}, \\
      &\zzj{i}[N,k,\ell,m,\zeta,\xi]
      =
      \zzj{i}[N,k,\ell,m,\zeta,0]
        +
      \begin{cases}
        \tauj{1}\xij{1} \mbox{ for $i=1$} \\
        \tauj{N-1}\xij{N-1} \mbox{ for $i=N$} \\
        \frac{1}{2}\left(\tauj{i-1}\xij{i-1}+\tauj{i}\xij{i}\right)
          \mbox{ for $1<i<N$,}
      \end{cases} \\
      &\zzj{i}^K[N,k,\ell,m,\zeta,0]
      =
      \zzj{i}[N,k,\ell,m,\zeta,0] + \tauj{i} \ln \frac{1}{10\ell m \tauj{i}}
        \mbox{ for $1 \leq i \leq N-1$}, \\
      &\zzj{i+1}[N,k,\ell,m,\zeta,0]
      =
      \zzj{i}[N,k,\ell,m,\zeta,0] + 2\tauj{i} \ln \frac{1}{10\ell m \tauj{i}}
        \mbox{ for $1 \leq i \leq N-1$, and} \\
      &\zzj{n}^K[2n,k,\ell,m,\zeta,0]
      =
      \zzj{n+1}[2n+1,k,\ell,m,\zeta,0]
      =
      0 \mbox{ (recalling \eqref{ndef})}.
    \end{aligned}
  \end{equation}
From definition \eqref{taudef}
  \begin{equation}
  \label{almosta}
    \ln \frac{1}{10\ell m \tauj{i}}
    =
    \frac{k \ell}{4\pi}\left(1-\frac{\cj{2}}{2}\right)m^2 - \zetaj{1}
    -\left(1-\delta_{i1}\right)\left(\ln \cj{i}+\frac{\zetaj{i}}{k \ell m^2}\right)
      \mbox{ for $1 \leq i \leq N-1$},
  \end{equation}
where
$\delta_{i1}=\left\{ \begin{smallmatrix}1 \mbox{ if $i=1$} \\0 \mbox{ otherwise,}\end{smallmatrix}\right.$
so by taking $m$ sufficiently large in terms of $\zeta$, $\xi$,
and each $\cj{i}$
(which according to \eqref{cisummary} satisfy bounds depending on just $N$, $k$, and $\ell$),
we can make every $\zzj{i}$ and $\zzj{i}^K$ as close to $0$ as desired
and we can also guarantee that
  \begin{equation}
    \zzj{i}[N,k,\ell,m,\zeta,\xi] < \zzj{i}^K[N,k,\ell,m,\zeta,0] < \zzj{i+1}[N,k,\ell,m,\zeta,\xi]
      \mbox{ for $1 \leq i \leq N-1$}.
  \end{equation}

The definitions \eqref{torheights} and \eqref{catheights}
can be understood as implementing the matching conditions
mentioned earlier as well as the vertical offsets introduced by the $\xi$ parameters.
Here the logarithm is used---merely because
it simplifies some expressions later---to capture
the dominant part of the inverse hyperbolic cosine function
for large values of its argument.
Thus each logarithmic term, ignoring any powers of $2$ appearing as prefactors,
represents the height achieved by a corresponding catenoid
above its waist plane a distance $\frac{1}{10 \ell m}$ from its axis,
where the catenoids are meant to transition to planes (tori under
$\Phi$). The factor of $10 \ell$ is chosen---$10$ somewhat arbitrarily
and $\ell$ because we assume $k \leq \ell$---to ensure the transition is completed
on the order---$m^{-1}$---of the lattice spacing but well away from neighboring
catenoids.
From the identity
  \begin{equation}
    \label{lnarcosh}
    \arcosh x = \ln x + \ln \left(1+\sqrt{1-x^{-2}}\right)
  \end{equation}
we have the estimate
  \begin{equation}
  \label{arcoshln}
    \abs{\tauj{i} \arcosh \frac{1}{10\ell m\tauj{i}}-\tauj{i} \ln \frac{1}{10\ell m\tauj{i}}}
    \leq
    \tauj{i} \ln 2,
  \end{equation}
so our very substitution of $\ln$ for $\arcosh$ introduces
additional height mismatch (beyond that contributed by the $\xi$ parameters)
where catenoidal annuli connect to tori,
but one whose ratio to $\tauj{1}$ is obviously bounded independently
of the $\zeta$ and $\xi$ parameters.

  \begin{sidewaysfigure}
    \centering
      \includegraphics[width=9in]{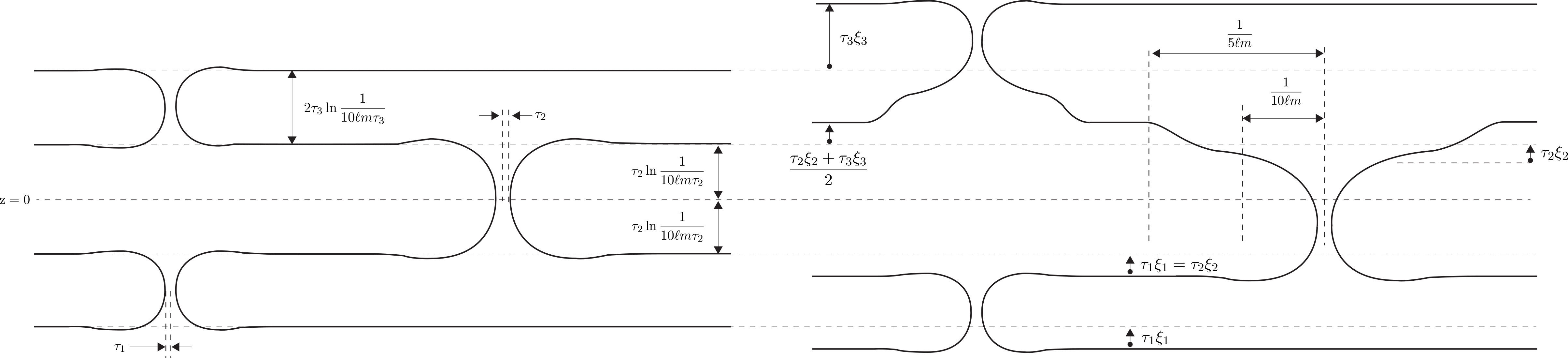}
      \caption{A schematic profile of two $N=4$ initial surfaces.
      The surface on the left has $\xij{1}=\xij{2}=\xij{3}=0$, while the one on the right
      has nonzero values for all three of these parameters,
      with $\tauj{1}\xij{1}=\tauj{2}\xij{2}\neq\tauj{3}\xij{3}$,
      resulting in dislocation at the third torus from the bottom
      and nowhere else.
      The figure is not drawn to scale.
      In fact (see especially \eqref{taudef} and \eqref{almosta})
      the vertical spacing $\tauj{i} \ln \frac{1}{10\ell m \tauj{i}}$
      is of order $m^2\tauj{1}$
      and $\lim_{m \to \infty} \frac{m^2\tauj{1}}{m^{-1}}=0$.
      Note also that the schematic does not attempt to convey
      the height mismatch, of order $\tauj{1}$, arising from
      the approximation of $\arcosh$ by $\ln$ (see \eqref{arcoshln}),
      because, unlike the dislocations controlled by $\xi$,
      it plays no role in the construction.}
  \end{sidewaysfigure}

\subsection*{Symmetries and horizontal specifications}
By the definition \eqref{Cliffordtorusdef} of $\T$ as the locus of points in $\Sph^3$
equidistant from $C_1$ and $C_2$,
the symmetry group $\Grp_{sym}[\T]$ of $\T$ in $\Sph^3$,
that is the subgroup of $O(4)$ preserving $\T$ as a set,
is precisely the symmetry group of $C_1 \cup C_2$.
Accordingly $\Grp_{sym}[\T]$ is generated by
the set of all reflections through great spheres containing either $C_1$ or $C_2$
together with the set of reflections through any great circle equidistant from $C_1$ and $C_2$
(so contained in a plane of the form
$\{z_1=e^{i\theta}z_2\}$ or $\{z_1=e^{i\theta}\overline{z}_2\}$)
and therefore lying on $\T$.
(Here reflection through a great sphere (or circle) refers to the element of $O(4)$
that identically preserves the $3$-space (or $2$-plane) containing that sphere (or circle)
and reflects the orthogonal complement through the origin.)
Consequently $\Grp_{sym}[\T]$ also includes
rotation in $C_1$ and $C_2$ by any angles
as well as reflection through any great circle orthogonally intersecting $\T$
(and $C_1$ and $C_2$).

We equip $\C^2$, $C_1$, and $C_2$ with their standard orientations and, given an oriented
circle $C$ in $\Sph^3$, we write $\rot_C^\theta$ for rotation by $\theta$ about $C$,
which by definition is the element of $O(4)$ fixing the plane containing $C$ pointwise
and rotating its orthogonal complement by angle $\theta$ in the direction consistent
with the induced orientation.
We write $\xbar$ for reflection
through the great sphere with equator $C_2$ and poles $(\pm 1,0)$,
$\ybar$ for reflection
through the great sphere with equator $C_1$ and poles $(0,\pm 1)$,
and $\zbar$ for reflection
through the great circle in the plane $\{z_1=z_2\}$
(so $\zbar=\rot_{\{z_1=z_2\} \cap \Sph^3}^\pi$).
Thus (writing $\overline{z}$ for the complex conjugate of $z \in \C$)
  \begin{equation}
    \begin{aligned}
      &\xbar(z_1,z_2) = (\overline{z_1},z_2),
      \qquad
      \ybar(z_1,z_2) = (z_1,\overline{z_2}),
      \qquad
      \zbar(z_1,z_2)=(z_2,z_1), \\
      &\rot_{C_2}^\theta (z_1,z_2) = (e^{i\theta}z_1,z_2),
      \qquad \mbox{and} \qquad
      \rot_{C_1}^\theta (z_1,z_2) = (z_1, e^{i\theta}z_2)
    \end{aligned}
  \end{equation}
and
  \begin{equation}
  \Grp_{sym}[\T]=\left\langle \xbar,\ybar,\zbar,\rot_{C_1}^{\theta_1}, \rot_{C_2}^{\theta_2}
     \; : \; \theta_1,\theta_2 \in \R \right\rangle,
  \end{equation}
where the right-hand side is the subgroup of $O(4)$
generated by the elements the angled brackets enclose.

Recalling \eqref{Phidef}, we also name some elements of $O(3)$ preserving
the domain of $\Phi$ and the pullback metric $\Phi^*\gsph$:
we write $\transe_{\xx\text{-axis}}^h$
and $\transe_{\yy\text{-axis}}^h$ for translation by (real number) $h$
in the (positive) $\xx$ and $\yy$ directions respectively,
$\xbarhat$ and $\ybarhat$ for reflection through the $\xx=0$ and $\yy=0$ planes respectively,
and $\zbarhat$ for reflection through the line $\xx=\yy$ in the plane $\zz=0$.
Thus
  \begin{equation}
    \begin{aligned}
      &\xbarhat(\xx,\yy,\zz)=(-\xx,\yy,\zz),
      \qquad
      \ybarhat(\xx,\yy,\zz)=(\xx,-\yy,\zz),
      \qquad
      \zbarhat(\xx,\yy,\zz)=(\yy,\xx,-\zz), \\
      &\transe_{\xx\text{-axis}}^h(\xx,\yy,\zz)=(\xx+h,\yy,\zz),
      \qquad \mbox{and} \qquad
      \transe_{\yy\text{-axis}}^h(\xx,\yy,\zz)=(\xx,\yy+h,\zz).
    \end{aligned}
  \end{equation}
It is then easy to verify the relations
  \begin{equation}
  \label{Phintertwine}
    \begin{aligned}
      \Phi \circ \xbarhat = \xbar \circ \Phi, \qquad
      \Phi \circ \ybarhat &= \ybar \circ \Phi, \qquad
      \Phi \circ \zbarhat = \zbar \circ \Phi, \\
      \Phi \circ \transe_{\xx\text{-axis}}^h = \rot_{C_2}^{\sqrt{2}h} \circ \Phi, &\qquad \mbox{and} \qquad
      \Phi \circ \transe_{\yy\text{-axis}}^h = \rot_{C_1}^{\sqrt{2}h} \circ \Phi
      \quad \mbox{ for every $h \in \R$},
    \end{aligned}
  \end{equation}
which play an essential role in identifying the symmetries of the initial surfaces,
defined via $\Phi$, and in exploiting these symmetries when solving the linearized problem.

Each initial surface will be built to be invariant under a certain subgroup of $\Grp_{sym}[T]$,
with the symmetries of $\T$ broken in several ways.
Recall that each initial surface is to be obtained by gluing together a collection
of constant-mean-curvature tori coaxial with $\T$ using a collection of catenoidal annuli.
Since among these tori only $\T$ itself is equidistant from $C_1$ and $C_2$,
the symmetry group of every other such torus is just the subgroup
$\left\langle \xbar,\ybar,\rot_{C_1}^{\theta_1}, \rot_{C_2}^{\theta_2}
  \; : \; \theta_1,\theta_2 \in \R \right\rangle$
of horizontal symmetries of $\T$.
If the toral components of an initial surface are
not arranged symmetrically about $\T$,
then of course the initial surface itself cannot possess any vertical symmetries.
The arrangement of the catenoidal tunnels may also be inconsistent with vertical symmetry
and inevitably breaks the continuous horizontal symmetries,
leaving only a discrete subgroup.
In this construction we impose the group
  \begin{equation}
  \label{Grpdef}
    \Grp
    =
    \Grp[k,\ell,m]
    :=
    \left\langle \rot_{C_2}^{\frac{2\pi}{km}}, \,
      \rot_{C_1}^{\frac{2\pi}{\ell m}}, \,
      \xbar, \, \ybar \right\rangle
    < \Grp_{sym}[\T] < O(4),
  \end{equation}
which is the subgroup of $O(4)$ preserving
the set of $km$\textsuperscript{th} roots of unity on $C_1$
as well as (separately) the set of $\ell m$\textsuperscript{th} roots of unity on $C_2$;
as such, $\Grp[k,\ell,m]$ is therefore isomorphic to $D_{km} \times D_{\ell m}$,
where $D_q$ is the dihedral group of order $2q$.
We will sometimes write $\Grp$ in place of $\Grp[k,\ell,m]$
when there is no danger of confusion.

The catenoidal tunnels connecting each pair of adjacent tori
will be placed, via $\Phi$, so as to take their centers on
the great circles orthogonally intersecting $\T$ at the sites of
a $\Grp[k,\ell,m]$-invariant $km \times \ell m$ rectangular lattice.
There are precisely four such lattices, namely the $\Grp[k,\ell,m]$ orbits
  \begin{equation}
  \label{lattices}
  \begin{aligned}
    &L_{\sigma_{\xx},\sigma_{\yy}}=L_{\sigma_{\xx},\sigma_{\yy}}[k,\ell,m]
    :=
    \Grp[k,\ell,m]
      \left(e^{i \sigma_{\xx}\pi/(km)},e^{i \sigma_{\yy}\pi/(\ell m)}\right)
    =
    \Phi \left(\lhat_{\sigma_{\xx},\sigma_{\yy}} \times \{0\} \right), \mbox{ where} \\
    &\lhat_{\sigma_{\xx},\sigma_{\yy}}=\lhat_{\sigma_{\xx},\sigma_{\yy}}[k,\ell,m]
    :=
      \left.
        \left\{
          \left(
            \frac{\sigma_{\xx}\pi}{\sqrt{2}km}+n_{\xx}\frac{\sqrt{2}\pi}{km},
            \frac{\sigma_{\yy}\pi}{\sqrt{2}\ell m}+n_{\yy}\frac{\sqrt{2}\pi}{\ell m}
          \right) \in \R^2
        \; \right| \; n_{\xx}, n_{\yy} \in \Z
        \right\},
  \end{aligned}
  \end{equation}
corresponding to the four choices of $\sigma_{\xx},\sigma_{\yy} \in \{0,1\}$.
To avoid further complicating the definition of the initial surfaces
in this article we make use of just the two lattices
$L_{0,0}$ and $L_{1,1}$,
but we emphasize that only obvious, minor modifications
to our procedure would be required
to take advantage of the other two lattices as well,
which freedom would allow us to produce minimal surfaces not congruent to the
ones we explicitly construct here.
(It would also be natural to attempt to adapt the construction
to permit more generally $\Grp[k,\ell,m]$-invariant
$n_i km \times n_i\ell m$ lattices refining the above four,
where the positive integer $n_i$ could be allowed to vary from layer to layer.)
See Remark \ref{latticeremark}.

Each initial surface will be invariant under the corresponding $\Grp[k,\ell,m]$,
and we will later admit only deformations respecting this group,
so in fact each minimal surface ultimately produced will also be invariant
under the corresponding $\Grp[k,\ell,m]$.
When $k \neq \ell$, $\Grp[k,\ell,m]$ is in fact the largest group
preserving each lattice
and the largest group preserving the set of centers of the catenoidal tunnels,
and it will consequently be the full symmetry group of the resulting minimal surface.
When $k=\ell$, however, there are choices of parameters $\zeta,\xi$
such that the corresponding initial surface is
also preserved by reflection through certain great circles on $\T$.
In this case it would be possible to modify the procedure that follows
by cutting in half the number of free parameters,
imposing constraints respecting these additional symmetries;
the construction would then produce minimal surfaces also invariant under the extra symmetries.
We will not pursue this modification here, instead enforcing just the smaller group
$\Grp[k,\ell,m]$ even in the square case.
From the analysis of this article alone
it is unclear whether or not the resulting minimal surfaces
nevertheless enjoy the additional symmetries.
See Remark \ref{verticalsymremark}.

\subsection*{Assembly and basic properties}
Suppose we are given the following data as above:
integers
$N \geq 2$, $\ell \geq k \geq 1$, and $m \geq 1$,
as well as two vectors $\zeta,\xi \in \R^{N-1}$.
We set $\cj{2}:=0$ if $N=2$ and $\cj{2}:=1$ if $N=3$,
and for $N \geq 4$ we accept a collection
$\{\cj{j}[N,k,\ell,m]\}_{j=2}^n \subset (1,\infty)$
as described in Lemma \ref{cilemma}
(namely the one determined in its proof).
The set of waist radii $\{\tauj{j}\}_{j=1}^{N-1}$ is then defined by \eqref{taudef},
the set of toral heights $\{\zzj{j}\}_{j=1}^N$ by \eqref{torheights},
and the set of catenoidal heights $\{\zzj{j}^K\}_{j=1}^{N-1}$ by \eqref{catheights}.
Setting
  \begin{equation}
    \label{RXYdef}
    R=R[\ell,m]:=\frac{1}{10 \ell m}, \qquad
    X=X[k,\ell,m]:=\frac{\pi}{\sqrt{2}km}, \quad \mbox {and} \quad
    Y=Y[k,\ell,m]:=\frac{\pi}{\sqrt{2}\ell m}
  \end{equation}
and recalling \eqref{Phidef}, \eqref{Text}, and \eqref{Tint},
for each integer $1 \leq i \leq N$ we define the parametrized surface-with-boundary $\Omega_i$
in $\Sph^3$ by
  \begin{equation}
  \label{Omega}
    \begin{aligned}
      &\Omega_N
      :=
      \Phi\left(T_{ext}[(N-2)(X,Y),
        \zzj{N-1}^K,\zzj{N},R,X,Y,\tauj{N-1}](\T_{X,Y,\tauj{N-1}})\right), \\
      &\Omega_1
      :=
      \Phi\left(
        T_{ext}[(0,0),\zzj{1}^K,\zzj{1},R,X,Y,\tauj{1}](\T_{X,Y,\tauj{1}})\right),
        \mbox{ and for } 1<i<N \\
      &\Omega_i
      :=
      \Phi\left(T_{int}
        \left[
          (2i-3)\left(\frac{X}{2},\frac{Y}{2}\right),
          \zzj{i-1}^K, \zzj{i}^K, \zzj{i}, R, X, Y, \tauj{i-1}, \tauj{i}
        \right]
        (\T_{X, Y, \tauj{i-1}, \tauj{i}})\right).
    \end{aligned}
  \end{equation}

Note that each $\Omega_i$ is a
$\frac{\sqrt{2}\pi}{km} \times \frac{\sqrt{2}\pi}{\ell m}$
rectangular patch of the constant-mean-curvature torus
at signed distance $\zzj{i}$ (increasing toward $C_1$) from $\T$,
within which patch
one (for $i=1$ or $i=N$) or two (for $1<i<N$) discs
have been replaced (via $\Phi$ and a cutoff function) by catenoidal annuli,
so that the boundary of $\Omega_i$ is the union of a rectangle
with one or two waist cicles,
the nearest-point projection onto $\T$
of the center of each deleted disc
is a site of either $L_{0,0}$ or $L_{1,1}$
(recalling \eqref{lattices}),
and $\bigcup_{i=1}^N \Omega_i$ is a smooth connected surface
whose boundary is the union of $N$ rectangles.

Finally we define the initial surface
  \begin{equation}
  \label{initsurfdef}
    \Sigma
    =
    \Sigma[N,k,\ell,m,\zeta,\xi]
    :=
    \Grp[k,\ell,m] \bigcup_{i=1}^N \Omega_i,
  \end{equation}
the orbit under $\Grp[k,\ell,m]$ of $\bigcup_{i=1}^N \Omega_i$.

\begin{prop}[Basic properties of the initial surfaces]
\label{initsurfprops}
Given a real number $c>0$ and integers $N \geq 2$ and $\ell \geq k \geq 1$,
there exists $m_0=m_0[N,k,\ell,c]>0$ such that for every integer
$m \geq m_0$ and every choice of parameters $\zeta,\xi \in [-c,c]^{N-1}$
the initial surface $\Sigma[N,k,\ell,m,\zeta,\xi]$
(defined by \eqref{initsurfdef}) is a smooth, closed surface, embedded in $\Sph^3$,
of genus $k \ell m^2(N-1) + 1$,
and invariant as a set under the action of
$\Grp[k,\ell,m]$ (defined in \eqref{Grpdef}).
\end{prop}

\begin{proof}
By \eqref{cisummary}, \eqref{taudef}, and \eqref{almosta}
for fixed $N$, $k$, $\ell$, and $c$,
we have
  \begin{equation}
    \lim_{m \to \infty}
      \left(
        \frac{\abs{c}}{m}
        +\sum_{i=1}^{N-1}10\ell m \tauj{i}
        +\sum_{i=1}^{N-1}\tauj{i}\ln \frac{1}{10\ell m \tauj{i}}
      \right)
    =
    0,
  \end{equation}
ensuring embeddedness of $\bigcup_{i=1}^N \Omega_i$.
All the claims are now clear from
\eqref{Phintertwine}, \eqref{Grpdef},
and the explicit construction of $\Sigma$.
\end{proof}

\begin{remark}[Smooth dependence on the parameters]
\label{initsurfcty}
Since for any fixed $N$, $k$, $\ell$, and $m$,
the quantities \eqref{taudef}, \eqref{torheights}, and \eqref{catheights}
all depend smoothly on the $\zeta,\xi$ parameters,
it follows from \eqref{initsurfdef}
and the supporting definitions, particularly \eqref{Text} and \eqref{Tint},
that the initial surface $\Sigma[N,k,\ell,m,\zeta,\xi]$
depends smoothly on $\zeta$ and $\xi$
in the sense that there exists a smooth map
$I=I[N,k,\ell,m]: \R^{N-1} \times \R^{N-1} \times \Sigma[N,k,\ell,m,0,0] \to \Sph^3$
such that
for any $\zeta,\xi \in \R^{N-1}$
the map $I[N,k,\ell,m](\zeta,\xi,\cdot)$ is an embedding
(provided $m$ is large enough in terms of $\abs{\zeta}$ and $\abs{\xi}$)
with image $\Sigma[N,k,\ell,m,\zeta,\xi]$.
In particular $I(\zeta,\xi,\cdot)$ is a diffeomorphism onto its image,
and, in casual abuse of notation,
we will routinely write $I(\zeta,\xi,\cdot)$
for this diffeomorphism,
so that
$I[N,k,\ell,m](\zeta,\xi,\cdot)^{-1}: \Sigma[N,k,\ell,m,\zeta,\xi] \to \Sigma[N,k,\ell,m,0,0]$
is the inverse diffeomorphism.
\end{remark}

Of course it remains to specify $\{\cj{j}\}_{j=2}^n \subset (1,\infty)$ when $N \geq 4$,
a gap filled in the next section.

\section{Forces and dislocations}

\subsection*{Forces}
As mentioned in Section 1, we will eventually discover (in Section \ref{linop})
that on each toral region (not yet defined),
as $m$ tends to infinity,
the Jacobi operator
converges,
after appropriate rescaling,
to a limit operator on a corresponding limit region
(likewise after rescaling),
and this limit operator has one-dimensional kernel.
The nontriviality of this kernel will compel us,
when attempting to prescribe mean curvature at the linear level by graphical deformation,
to modify the term to be prescribed by elements of the substitute kernel
(as outlined in Section 1 and formally defined in Section \ref{linop}).
The details appear in Section \ref{linop},
where the precise argument is couched in rather different terms,
but, in order to help motivate the study undertaken in the current section,
we offer in this paragraph the following rough sketch of the situation
from a purely geometric point of view.
In the next section (specifically \eqref{torrdef}) we will define
the $i$\textsuperscript{th} toral region $\torr[i]$
to be $\Omega_i$ \eqref{Omega}
less a certain portion of the half catenoid(s) attached.
Of course the area of $\Omega_i$ (and so of $\torr[i]$ too) shrinks to $0$ as $m$ goes to infinity,
but if we scale the ambient spherical metric $\gsph$ up to $m^2\gsph$,
then under the corresponding induced metric $\torr[i]$ tends 
to a flat $\frac{\sqrt{2}\pi}{k} \times \frac{\sqrt{2}\pi}{\ell}$ solid rectangle
(because the annuli deleted from $\Omega_i$ to define $\torr[i]$
are sized so as to vanish in this limit).
The corresponding limit Jacobi operator is just the flat Laplacian on this rectangle,
its Jacobi operator as a submanifold of Euclidean $\R^3$
(viewed in the domain of the map $\Phi$ \eqref{Phidef}).
By enforcing the symmetry group $\Grp$ \eqref{Grpdef} throughout the construction
we impose periodic boundary conditions on $\torr[i]$,
so that this limit Jacobi operator has kernel spanned by the constants.
The Jacobi field $1$ is induced by the Euclidean Killing field $\partial_{\zz}$
(viewed in the domain of $\Phi$) orthogonal to the rectangle.

Of course the vector field $\partial_{\zz}$ is not Killing relative to $\Phi^*\gsph$ \eqref{Phullback},
but it is approximately Killing relative to $m^2\Phi^*\gsph$,
and, conveniently, on a neighborhood of $\Omega_i$ in $\Sph^3$,
(viewed through $\Phi$) $\partial_{\zz}$ is itself approximated
by the exact Killing field $K$ on $(\Sph^3,\gsph)$ that generates rotation,
toward $C_1$ (defined just below \eqref{Cliffordtorusdef}),
along the great circle through the points $(1,0), (0,1) \in \C^2$.
The geodesic segment joining these points is simply
the closure of the image under $\Phi$ of the segment of the $\zz$-axis contained in the domain of $\Phi$;
thus $K \circ \Phi=\Phi_*\partial_{\zz}$ along this segment.
More generally (recalling \eqref{Phidef}),
\begin{equation}
  \label{Killing}
  \begin{aligned}
    K \circ \Phi = &-\frac{1}{\sqrt{2}} \cot \left(\zz + \frac{\pi}{4}\right)
          \sin \sqrt{2}\xx \, \cos \sqrt{2}\yy \; \Phi_*\partial_{\xx}
        + \frac{1}{\sqrt{2}} \tan \left(\zz + \frac{\pi}{4}\right)
          \cos \sqrt{2}\xx \,  \sin \sqrt{2}\yy \; \Phi_*\partial_{\yy} \\
    &+ \cos \sqrt{2}\xx \, \cos \sqrt{2}\yy \; \Phi_*\partial_{\zz}.
  \end{aligned}
\end{equation}

We will now calculate the $K$ force
(also called \emph{flux} in the literature) through various regions of the initial surface
and to study its dependence on the $\zeta,\xi$ parameters.
These forces will indirectly measure
the projection of the initial surface's mean curvature onto the substitute kernel,
so in the final section we will apply the results of this section
(specifically Lemma \ref{FD}) to help manage the substitute kernel
and complete the proof of the main theorem.
More immediately we will impose
\emph{balancing conditions} (\cite{KKS}, \cite{KapDelaunay}, \cite{KapClay})
on the initial surface,
such that the $K$ force on various regions vanishes,
at least within a margin on the order of the perturbations---by functions and
parameters---that we will be making.
This balancing will finally determine the waist radii, up to choice of $\zeta$,
thus completing the definition of the initial surfaces.

Let
  \begin{equation}
  \label{forcedef}
    \Fcal_i
    =\Fcal_i[N,k,\ell,m,\zeta,\xi]
    :=
    \int_{\partial \Omega_i} (\gsph \circ \iota)(K \circ \iota, \eta_i)
      \, \sqrt{\abs{\iota^*\gsph|_{\partial \Omega_i}}}
    =
    \int_{\Omega_i} (\gsph \circ \iota) (K \circ \iota,\mathbf{H}) \, \sqrt{\abs{\iota^*\gsph}},
  \end{equation}
the \emph{$K$ force} exerted by the region $\Omega_i$ (defined by \eqref{Omega})
on the rest of $\Sigma$,
where $\iota: \Sigma \to \Sph^3$ is the inclusion map of $\Sigma$ in $\Sph^3$,
$\eta_i$ is the outward conormal for $\Omega_i$,
$\mathbf{H}:=\tr_{\iota^*\gsph} Dd\iota$ is the vector-valued mean curvature of $\Sigma$
($D$ being the connection induced on $T^*\Sigma \otimes \iota^*T\Sph^3$ by $\gsph$ and $\iota$),
$\sqrt{\abs{\iota^*\gsph}}$
and $\sqrt{\abs{\iota^*\gsph|_{\partial \Omega_i}}}$
are respectively the area and length forms induced by $\iota$ and $\gsph$,
and the last equality follows trivially
from the formula for the first variation of area and the fact that $K$ is Killing.
Since the initial surface should be designed to be approximately minimal,
by virtue of the second equality
we can now impose the approximate (in a sense made precise below) balancing condition
  \begin{equation}
  \label{balance}
    \Fcal_i \approx 0 \mbox{ for each $1 \leq i \leq N$}
  \end{equation}
in order to estimate the necessary waist radii,
before varying the parameters or deforming the initial surface graphically.
We will see that this heuristic approach leads to \eqref{taubardef}.
We will also analyze how
by adjusting the parameters
we can fine-tune these radii as well as the heights of the catenoids and tori
in order to control the forces.

The computation of the forces is simple.
The boundary of each $\Omega_i$ consists of
a rectangle on a constant-mean-curvature torus coaxial with $\T$
and one (for $i \in \{1,N\}$)
or two (for $2 \leq i \leq N-1$) catenoidal waists.
For each such component,
by working in an $(\xx,\yy,\zz)$ coordinate system defined via $\Phi$ \eqref{Phidef},
we will estimate the corresponding integral
arising in \eqref{forcedef}.
Suppose $S$ is a catenoidal waist with center
$(\xx_{i}^K,\yy_{i}^K,\zzj{i}^K)$.
The height $\zzj{i}^K$ has been defined in \eqref{catheights},
and the horizontal coordinates $\xx_{i}^K,\yy{_i}^K$,
though not previously named, can be found in \eqref{Omega}.
The outward unit conormal $\eta_{_S}$ along $S$ is, recalling \eqref{Phullback},
simply $\pm \partial_{\zz}$
and, relative to $\geuc$, $S$ is a Euclidean circle of radius $\tauj{i}$,
so by \eqref{Phullback} and \eqref{Killing}
  \begin{equation}
  \label{waist}
    \int_{S} \gsph(K,\eta_{_S})
    =
    \int_0^{2\pi} \cos \left[\sqrt{2}\left(\xx_{i}^K+\tauj{i} \cos \theta\right)\right]
      \cos \left[\sqrt{2}\left(\yy_{i}^K+\tauj{i} \sin \theta\right)\right]
      \tauj{i} \sqrt{1-(\sin 2\zzj{i}^K)(\cos 2\theta)} \, d\theta.
  \end{equation}

Now suppose $T$ is a rectangular component of $\partial \Omega_i$ for some admissible $i$,
with outward conormal $\eta_{_T}$ and center $(\xx_{i},\yy_{i},\zzj{i})$.
The height $\zzj{i}$ has been defined in \eqref{torheights},
and the horizontal coordinates $(\xx_{i},\yy_{i})$
can be found in \eqref{Omega}.
Thus $T$ lies on the constant-mean-curvature torus $\{\zz=\zzj{i}\}$
and its complement in this torus has two connected components,
the smaller of which (in terms of area) we call $\overline{T}$,
a solid rectangle in $\{\zz=\zzj{i}\}$ having boundary
$\partial{\overline{T}}=T$.
Because of the way $\Omega_i$ is defined using cutoff functions,
there is a neighborhood $\mathcal{U}$ of $T$ in $\Sph^3$
such that $\mathcal{U} \cap \Omega_i=\mathcal{U} \cap \overline{T}$,
and therefore $\eta_{_T}$ is equally
the outward unit conormal for $\overline{T}$ along $T$.
By invoking the first-variation-of-area formula
(as in the second equality of \eqref{forcedef})
we have $\int_T \gsph(K,\eta_{_T})=\int_{\overline{T}} \gsph(K,\mathbf{H}_{_T})$,
where $\mathbf{H}_{_T}$ is the mean curvature of $\overline{T}$.
From \eqref{Phullback} we find that the area form on $\{\zz=\zzj{i}\}$
is $\cos 2\zzj{i} \, d\xx \, d\yy$,
whence $\mathbf{H}_{_T}=2 \tan 2\zzj{i} \partial_{\zz}$,
so using also \eqref{Omega} and \eqref{Killing}
we get
  \begin{equation}
  \label{rectangle}
    \int_T \gsph(K,\eta_{_T})
    =
    \int_{\xx_{i}-\frac{\sqrt{2}\pi}{km}}^{\xx_{i}+\frac{\sqrt{2}\pi}{km}}
      \int_{\yy_{i}-\frac{\sqrt{2}\pi}{\ell m}}^{\yy_{i}+\frac{\sqrt{2}\pi}{\ell m}}
        2 \cos (\sqrt{2}\xx) \cos (\sqrt{2}\yy) \tan (2\zzj{i}) \cos (2\zzj{i}) \, d\xx \, d\yy.
  \end{equation}

To estimate the integrals \eqref{waist} and \eqref{rectangle}
we will make the approximations
$\cos u \approx 1$, $\sin u \approx u$, and $\sqrt{1+u} \approx 1+u/2$.
On a heuristic level we could ignore the error in these approximations
and proceed with the calculation
to motivate \eqref{taubardef} as advertised.
On the other hand, the actual construction
will demand more detailed estimates,
so, to avoid repeating some calculations,
we will instead keep track of the error as we go,
predicating the estimates on \eqref{cisummary} and \eqref{taudef},
which themselves are suggested by the more cavalier approach.
From \eqref{cisummary}, \eqref{taudef}, \eqref{torheights}, \eqref{catheights},
and \eqref{Omega}
we see that whenever $m$ is sufficiently large in terms of $N$ for \eqref{cisummary} to hold,
there is some constant $C[N,k,\ell]>0$
(possibly larger than the one appearing in \eqref{cisummary}
but independent of $\zeta$, $\xi$, and $m$)
for which
  \begin{equation}
  \label{estimatesforforceestimate}
    \begin{aligned}
      &\abs{\xx_{i}^{(K)}}+\abs{\yy_{i}^{(K)}} \leq 10Nm^{-1}
        \mbox{ for $\xx_{i}$, $\yy_{i}$ and $\xx_{i}^K$, $\yy_{i}^K$
           appearing in \eqref{waist} and \eqref{rectangle}}, \\
      &\lim_{m \to \infty} m^q\tauj{1}=0
        \mbox{ uniformly in $\zeta \in [-c,c]^{N-1}$ for any fixed $q$, $c$, $N$, $k$, $\ell$}, \\
      &\max \{\tauj{i}/\tauj{1}+\tauj{1}/\tauj{i}\}_{i=2}^N
      \leq
      C[N,k,\ell]e^{\abs{\zetaj{i}}/m^2}, \mbox{ and} \\
      &\max \{\abs{\zzj{i}}\}_{i=1}^N + \max \{\abs{\zzj{i}^K}\}_{i=1}^{N-1}
      \leq
      C[N,k,\ell]N(m^2+\abs{\zeta}+\abs{\xi})\tauj{1}e^{\abs{\zeta}/m^2}
        \leq C[N,k,\ell]m^2\tauj{1},
    \end{aligned}
  \end{equation}
where $\abs{\zeta}$ and $\abs{\xi}$ are the Euclidean norms of $\zeta,\xi \in \R^{N-1}$
and for the final inequality we assume $m$ large in terms of $\abs{\zeta}$ and $\abs{\xi}$
(and allow a larger $C[N,k,\ell]$ on the right-hand side than previously needed).

Using \eqref{forcedef}, \eqref{waist}, \eqref{rectangle},
and \eqref{estimatesforforceestimate},
we conclude
  \begin{equation}
  \label{forcecomp}
    \begin{aligned}
      &\abs{2\pi\tauj{1} + \frac{8\pi^2}{k\ell m^2}\zzj{1} - \Fcal_1}
      \leq Cm^{-2}\tauj{1}, \\
      &\abs{-2\pi \tauj{N-1} + \frac{8\pi^2}{k\ell m^2}\zzj{N}-\Fcal_N}
      \leq
      Cm^{-2}\tauj{1}, \mbox{ and} \\
      &\abs{2\pi(\tauj{i} - \tauj{i-1}) + \frac{8\pi^2}{k\ell m^2}\zzj{i}-\Fcal_i}
      \leq
      Cm^{-2}\tauj{1} \mbox{ when } 2 \leq i \leq N-1
    \end{aligned}
  \end{equation}
for some constant $C=C[N,k,\ell]>0$ (independent of $m$ and $\zeta$ and $\xi$)
whenever $m$ is sufficiently large in terms of $N$, $\abs{\zeta}$, and $\abs{\xi}$.

It now follows from \eqref{forcecomp} and \eqref{torheights},
recalling \eqref{ndef},
that
  \begin{equation}
  \label{Fn}
    \begin{aligned}
      &\abs{2\pi \tauj{1} - \frac{8\pi^2}{k\ell m^2}\left(
        2^{N \bmod 2}\tauj{1} \ln \frac{1}{10 \ell m \tauj{1}}
        - \tauj{1}\xij{1}\right) - \Fcal_n}
      \leq
      Cm^{-2}\tauj{1} \mbox{ if } n=1, \\
      &\abs{2\pi(\tauj{n} - \tauj{n-1}) - \frac{8\pi^2}{k\ell m^2}\left(
        2^{N \bmod 2}\tauj{n} \ln \frac{1}{10 \ell m \tauj{n}}
        - \frac{\tauj{n-1}\xij{n-1} + \tauj{n}\xij{n}}{2}\right) - \Fcal_n} \\
      &\;\;\;\; \leq
      Cm^{-2}\tauj{1}
        \mbox{ if } n \geq 2,
    \end{aligned}
  \end{equation}
  \begin{equation}
  \label{Fn+1}
    \begin{aligned}
      &\abs{-2\pi\tauj{1}+\frac{8\pi^2}{k\ell m^2}\left(
        \tauj{1} \ln \frac{1}{10 \ell m \tauj{1}}
        + \tauj{1}\xij{1}\right)-\Fcal_{n+1}}
      \leq
      Cm^{-2}\tauj{1} \mbox{ if $N=2$}, \\
      &\abs{2\pi(\tauj{n+1} - \tauj{n}) + \frac{8\pi^2}{k\ell m^2}\left(
        [(N+1) \bmod 2] \tauj{n} \ln \frac{1}{10 \ell m \tauj{n}}
        + \frac{\tauj{n}\xij{n} + \tauj{n+1}\xij{n+1}}{2}\right) - \Fcal_{n+1}} \\
      &\;\;\;\; \leq
      Cm^{-2}\tauj{1} \mbox{ if } N \geq 3,
    \end{aligned}
  \end{equation}
and, when $N \geq 3$,
  \begin{equation}
  \label{forcediffs}
    \begin{aligned}
      &\abs{2\pi(\tauj{2} - 2\tauj{1}) + \frac{8\pi^2}{k\ell m^2}\left(
        2 \tauj{1} \ln \frac{1}{10 \ell m \tauj{1}}
        + \frac{\xij{2}\tauj{2} - \xij{1}\tauj{1}}{2}\right)
        -\left(\Fcal_{i+1}-\Fcal_i\right)}
      \leq
      Cm^{-2}\tauj{1}
        \mbox{ if } i=1, \\
      &\abs{2\pi(\tauj{i+1} - 2\tauj{i} + \tauj{i-1}) + \frac{8\pi^2}{k \ell m^2}\left(
        2 \tauj{i} \ln \frac{1}{10 \ell m \tauj{i}}
        + \frac{\xij{i+1}\tauj{i+1} - \xij{i-1}\tauj{i-1}}{2}\right)
        -\left(\Fcal_{i+1}-\Fcal_i\right)} \\
      &\;\;\;\; \leq
        Cm^{-2}\tauj{1}
        \mbox{ if } 2 \leq i \leq N-2, \mbox{ and } \\
      &\abs{2\pi(-2\tauj{N-1} + \tauj{N-2}) + \frac{8\pi^2}{k \ell m^2}\left(
        2 \tauj{N-1} \ln \frac{1}{10 \ell m \tauj{N-1}}
        + \frac{\xij{N-1}\tauj{N-1} - \xij{N-2}\tauj{N-2}}{2}\right)
        -\left(\Fcal_{i+1}-\Fcal_i\right)} \\
      &\;\;\;\; \leq
      Cm^{-2}\tauj{1} \mbox{ if } i=N-1
    \end{aligned}
  \end{equation}
(where $C$ can be taken to be twice the value of the $C$ appearing in \eqref{forcecomp}).

\subsection*{Balancing}
The force estimates \eqref{forcecomp}--\eqref{forcediffs}
will play an indispensable role, via Lemma \ref{FD}, in selecting viable parameter values
in the final steps of the construction,
but at the moment,
with the goal of completing the specification of the waist radii in the initial surfaces,
we set $\zeta=\xi=0$
and impose only the approximate balancing conditions \eqref{balance}
(since we will have to vary the parameters and graphically perturb
the initial surface anyway to achieve minimality):
temporarily ignoring the error bounded by the right-hand sides of
\eqref{Fn}--\eqref{forcediffs}, we demand
  \begin{equation}
  \label{balancing}
    \begin{aligned}
      &2\pi \taubarj{1} - 2^{N \bmod 2} \frac{8\pi^2}{k\ell m^2}
        \cdot \taubarj{1} \ln \frac{1}{10 \ell m \taubarj{1}}
      =0 \mbox{ if } n=1, \\
      &2\pi(\taubarj{n} - \taubarj{n-1}) - 2^{N \bmod 2} \frac{8\pi^2}{k\ell m^2}
        \cdot \taubarj{n} \ln \frac{1}{10 \ell m \taubarj{n}}
      =0 \mbox{ if } n \geq 2, \\
      &2\pi(\taubarj{n+1} - \taubarj{n}) + [(N+1) \bmod 2] \frac{8\pi^2}{k\ell m^2}
        \cdot \taubarj{n} \ln \frac{1}{10 \ell m \taubarj{n}}
      =0 \mbox{ if } N \geq 3, \\
      &2\pi(\taubarj{2} - 2\taubarj{1}) + \frac{8\pi^2}{k\ell m^2}
        \cdot 2 \taubarj{1} \ln \frac{1}{10 \ell m \taubarj{1}}
      =0 \mbox{ if $N \geq 3$}, \\
      &2\pi(\taubarj{i+1} - 2\taubarj{i} + \taubarj{i-1}) + \frac{8\pi^2}{k \ell m^2}
        \cdot 2 \taubarj{i} \ln \frac{1}{10 \ell m \taubarj{i}}
      =0 \mbox{ if $N \geq 3$ and $2 \leq i \leq N-2$, and} \\
      &2\pi(-2\taubarj{N-1} + \taubarj{N-2}) + \frac{8\pi^2}{k \ell m^2}
        \cdot 2 \taubarj{N-1} \ln \frac{1}{10 \ell m \taubarj{N-1}}
      =0 \mbox{ if $N \geq 3$}.
    \end{aligned}
  \end{equation}

From the third equation of \eqref{balancing}
we see that $\taubarj{n+1}=\taubarj{n}$ whenever $N \geq 3$ is odd,
while from this same equation together with the second equation
we see that $\taubarj{n+1}=\taubarj{n-1}$ whenever $N \geq 4$ is even
(and of course $\taubarj{n+1}$ is undefined for $N=2$);
thus $\taubarj{n+1}=\taubarj{N-(n+1)}$ whenever $N \geq 3$.
If $N$ is even, then $N=2n$, so obviously $\taubarj{n}=\taubarj{N-n}$,
while if $N$ is odd, then $N=2n+1$,
so $\taubarj{n}=\taubarj{N-(n+1)}$
and $\taubarj{N-n}=\taubarj{n+1}$
but we have just established that $\taubarj{N-(n+1)}=\taubarj{n+1}$;
thus we also have $\taubarj{n}=\taubarj{N-n}$ whenever $N \geq 2$.
It now follows by induction on $j$,
using the two equations
obtained by taking $i=j-1$ and $i=N-(j-1)$
in the penultimate line of \eqref{balancing},
having already dispensed in this paragraph with the cases $j=n$ and $j=n+1$,
that $\taubarj{j}=\taubarj{N-j}$ for each $j \in \Z \cap [n,N-1]$.
In fact it is clear that all the approximate balancing conditions
in \eqref{balancing} will be satisfied if and only if
we choose $\{\taubarj{i}\}_{i=1}^n$ satisfying
(for $n=1$) line 1 or
(for $n \geq 2$) lines 2, 4, and 5 (with $2 \leq i  \leq n-1$)
of \eqref{balancing}
and simultaneously set
  \begin{equation}
  \label{taubarrefl}
    \taubarj{i}:=\taubarj{N-i} \mbox{ for $n+1 \leq i \leq N-1$}.
  \end{equation}

For $n=1$ it therefore remains only to specify $\taubarj{1}$,
which is uniquely determined by imposing line 1 of \eqref{balancing}:
  \begin{equation}
  \label{taubar1lowN}
    \taubarj{1}
    :=
    \begin{cases}
      \frac{1}{10\ell m}e^{-\frac{k\ell m^2}{4\pi}} \mbox{ if } N=2 \\
      \frac{1}{10\ell m}e^{-\frac{k\ell m^2}{8\pi}} \mbox{ if } N=3.
    \end{cases}
  \end{equation}
For $n \geq 2$ (equivalently $N \geq 4$) we define
  \begin{equation}
  \label{cidef}
    \cj{i} := \taubarj{i}/\taubarj{1} \mbox{ for $1 \leq i \leq N-1$,}
  \end{equation}
so in particular $\cj{1}=1$;
dividing equations 4, 5, and 2 of \eqref{balancing} by $2\pi\taubarj{1}$,
we need now only solve the
$n$ equations
  \begin{equation}
  \label{intermediatecjsys}
    \begin{aligned}
      \cj{2} - 2
      &=
      \frac{8\pi}{k\ell m^2}
        \ln 10 \ell m \taubarj{1}, \\
      \cj{i+1} - 2\cj{i} + \cj{i-1}
      &=
      \frac{8\pi}{k \ell m^2}
        \cj{i} \ln 10 \ell m \cj{i}\taubarj{1}
        \mbox{ for } 2 \leq i \leq n-1, \mbox{ and} \\
      \cj{n-1} - \cj{n}
      &=
      \left(\frac{1}{2}\right)^{(N+1) \bmod 2} \frac{8\pi}{k\ell m^2}
        \cj{n} \ln 10 \ell m \cj{n}\taubarj{1}
    \end{aligned}
  \end{equation}
for the $n$ unknowns $\cj{2}, \ldots, \cj{n},$ and $\taubarj{1}$.

The first equation requires
  \begin{equation}
  \label{taubar1}
    \taubarj{1} := \frac{1}{10 \ell m}
    e^{-\frac{k \ell m^2}{4\pi}\left(1-\cj{2}/2\right)},
  \end{equation}
which we note even recovers \eqref{taubar1lowN} if we define
$\cj{2}=0$ for $N=2$ and $\cj{2}=1$ for $N=3$.
Assuming $N \geq 4$ (equivalently $n \geq 2$),
we derive a system
equivalent (presuming each $\cj{i} \neq 0$)
to the remaining $n-1$ equations of \eqref{intermediatecjsys}
by (i) for each $i \in \Z \cap [2,n-1]$ (a vacuous condition when $n=2$)
subtracting the middle equation of \eqref{intermediatecjsys}
from $\cj{i}$ times the top equation of \eqref{intermediatecjsys}
and (ii) subtracting $2^{(N+1) \bmod 2}$ times the bottom equation of \eqref{intermediatecjsys}
from $\cj{n}$ times the top equation of \eqref{intermediatecjsys}.
In this way (recalling $\cj{1}=1$) we obtain the system
  \begin{equation}
  \label{barbalance}
    \begin{aligned}
      &-\cj{i-1} + \cj{2}\cj{i} - \cj{i+1}
      =
      -\frac{8\pi}{k\ell m^2} \cj{i} \ln \cj{i}
        \mbox{ for } 2 \leq i \leq n-1 \mbox{ and} \\
      &-2^{(N+1) \bmod 2}\cj{n-1} + (\cj{2} - N \bmod 2)\cj{n}
      =
      -\frac{8\pi}{k\ell m^2}
        \cj{n} \ln \cj{n}.
    \end{aligned}
  \end{equation}

\begin{lemma}[Determination of the waist ratios by the approximate balancing conditions]
\label{cilemma}
Let $N \geq 4$ be a given integer and recall \eqref{ndef}.
There exist $n-1$ real numbers
$\djj{2}[N]<\djj{3}[N]<\cdots<\djj{n}[N]$,
with $\djj{2}[N] \in (1,2)$
and $\djj{2}[N]$ strictly increasing in $n$ (for a fixed parity of $N$),
and furthermore there exists $m_0=m_0[N]>0$
such that for each integer $m>m_0$
and for all integers $\ell \geq k \geq 1$
there are $n-1$ real numbers
$\cj{2}[N,k,\ell,m], \cj{3}[N,k,\ell,m], \ldots, \cj{n}[N,k,\ell,m]$
solving \eqref{barbalance}
and satisfying $\displaystyle{\lim_{m \to \infty}} \cj{i}[N,k,\ell,m] = \djj{i}[N]$
for any fixed $k$ and $\ell$.
\end{lemma}

\begin{proof}
Bear in mind that balancing has been accomplished by \eqref{taubar1lowN}
for $N=2$ and $N=3$.
Momentarily ignoring the logarithmic terms,
for $N=4$ the system \eqref{barbalance} reduces to $\djj{2}^2=2$,
so $\djj{2}[4]=\sqrt{2} \in (1,2)$,
while for $N=5$ we get $\djj{2}^2-\djj{2}-1 = 0$,
yielding $\djj{2}[5] = \frac{1+\sqrt{5}}{2} \in (1,2)$.
Now the functions $\cj{2} \mapsto \cj{2}^2$ and $\cj{2} \mapsto \cj{2}^2-\cj{2}-1$ have nonzero
derivatives at these respective values, so the lemma is established for $N=4$
and $N=5$ by applying the inverse function theorem and taking $m$ large.
Thus we may assume $n \geq 3$ and pursue an elaboration of the same strategy.

For real $\beta$ we define the $(n-1) \times (n-1)$ matrices
  \begin{equation}
  \label{Aevendef}
    A_{2n}(\beta)
    :=
    \begin{pmatrix}
      \beta & -1 & 0 & 0 & \cdots & 0 \\
      -1 & \beta & -1 & 0 & \cdots & 0 \\
      0 & -1 & \beta & -1 & \cdots & 0 \\
      \vdots & \vdots & \ddots & \ddots &\ddots &\vdots \\
      0 & 0 & \cdots & -1 & \beta & -1 \\
      0 & 0 & \cdots & 0 & -2 & \beta
    \end{pmatrix}
  \end{equation}
  and
  \begin{equation}
  \label{Aodddef}
    A_{2n+1}(\beta)
    :=
    \begin{pmatrix}
      \beta & -1 & 0 & 0 & \cdots & 0 \\
      -1 & \beta & -1 & 0 & \cdots & 0 \\
      0 & -1 & \beta & -1 & \cdots & 0 \\
      \vdots & \vdots & \ddots & \ddots &\ddots &\vdots \\
      0 & 0 & \cdots & -1 & \beta & -1 \\
      0 & 0 & \cdots & 0 & -1 & \beta-1
    \end{pmatrix},
  \end{equation}
so that the system obtained by temporarily replacing
all the logarithmic terms in \eqref{barbalance} by $0$
is equivalent (recalling $\cj{1}=1$) to the equation
  \begin{equation}
    A_N(\beta)
    \begin{pmatrix}
     \djj{2} \\
     \djj{3} \\
     \vdots \\
     \djj{n}
    \end{pmatrix}
    =
    \begin{pmatrix}
    1 \\
    0 \\
    \vdots \\
    0
    \end{pmatrix}
    \mbox{ subject to the constraints } \beta = \djj{2} \mbox{ and } \djj{i}>0
    \mbox{ for } 2 \leq i \leq n.
  \end{equation}

Using Cramer's rule and expansion by minors, we find
  \begin{equation}
  \label{di}
    \begin{aligned}
    \djj{i}[N]=&\frac{P_{n-i+1}[N \bmod 2](\beta)}{P_n[N \bmod 2](\beta)}
      \qquad \mbox{for $2 \leq i \leq n$, where } \\
    P_i[0](\lambda)=\det A_{2i}(\lambda)
    \qquad &\mbox{and} \qquad
    P_i[1](\lambda)=\det A_{2i+1}(\lambda) \qquad \mbox{for $i\geq 3$}, \\
    P_2[0](\lambda) = \lambda
    \qquad \qquad \quad &\mbox{and} \qquad
    P_2[1] = \lambda - 1, \mbox{ and} \\
    P_1[0](\lambda)=2
    \qquad \qquad \quad &\mbox{and} \qquad
    P_1[1](\lambda)=1.
    \end{aligned}
  \end{equation}
Further expansion by minors reveals the recursive relations
(independent of the parity of $N$)
  \begin{equation}
  \label{recur}
    P_{i+1}(\lambda) = \lambda P_i(\lambda) - P_{i-1}(\lambda) \mbox{ for }
    i \geq 2.
  \end{equation}
On the other hand,
by applying the constraint $\djj{2} = \beta$,
the expression for $\djj{2}$ given by \eqref{di}
can be rewritten as
  \begin{equation}
    \beta P_n(\beta) = P_{n-1}(\beta),
  \end{equation}
whence \eqref{recur} delivers
  \begin{equation}
    \label{betaeq}
    P_{n+1}(\beta) = 0.
  \end{equation}

We now claim that for each $i \geq 3$ (and either parity of $N$)
  \begin{enumerate}[(i)]
    \item $P_i=P_i[N \bmod 2]$ has a root strictly greater than $1$;
          if $\gamma_{i}=\gamma_{i}[N \bmod 2]$ is its largest such root, then
    \item $P_{i-1}(x) > 0$ whenever $x \geq \gamma_{i}$,
    \item $P_{i+1}(\gamma_{i}) < 0$, and
    \item $\gamma_{i}$ is strictly increasing in $i$.
  \end{enumerate}

These claims can be established by induction on $i$.
The case $i=3$ is easily verified:
$P_2[0](x)=x$, $P_3[0](x)=x^2-2$, and $P_4[0](x)=x^3-3x$,
so $\gamma_{3}[0] = \sqrt{2}$, $P_2[0](x \geq \gamma_{3}[0])>0$, and
$P_4[0](\gamma_{3}[0]) = 2\sqrt{2}-3\sqrt{2} < 0$, while
$P_2[1](x)=x-1$, $P_3[1](x)=x^2-x-1$, and $P_4[1](x)=x^3-x^2-2x+1$,
so $\gamma_{3}[1]=\frac{1+\sqrt{5}}{2}$,
$P_2[1](x \geq \gamma_{3}[1]) \geq \frac{\sqrt{5}-1}{2} > 0$,
and $P_4[1](\gamma_{3}[1]) = \frac{1-\sqrt{5}}{2} < 0$.
Now suppose that claims (i)-(iii) hold for $i=j$.
By claim (i) $\gamma_j$ exists and $\gamma_j>1$.
According to claim (iii) $P_{j+1}(\gamma_j)<0$,
but $P_{j+1}$ is clearly monic, so $P_{j+1}(x)>0$ for large $x$,
which implies that
$P_{j+1}$ has a root greater than $\gamma_j$,
so $\gamma_{j+1}$ exists
and $\gamma_{j+1}>\gamma_j>1$ (verifying claims (i) and (iv)).
Therefore
$P_j(x \geq \gamma_{j+1}) > 0$ by the maximality of $\gamma_j$ (verifying claim (ii)).
Finally,
using \eqref{recur},
$P_{j+2}(\gamma_{j+1}) = \gamma_{j+1}P_{j+1}(\gamma_{j+1}) - P_j(\gamma_{j+1})$,
which is negative, since the first term vanishes and the second has just
been established positive
(verifying claim (iii) and so completing the proof of claims (i)-(iv)).

Thus $\beta = \gamma_{n+1}$ solves \eqref{betaeq} and is strictly increasing in $n$
(for each fixed parity of $N$).
We have already checked that $\beta>1$;
now we claim that $\beta<2$.
In fact we assert that for each $i \geq 2$
(regardless of the parity of $N$)
  \begin{enumerate}[(v)]
    \item $P_i(x) - P_{i-1}(x) \geq 0$ and $P_i(x)>0$ whenever $x \geq 2$,
  \end{enumerate}
which is proven by induction on $i$. For $i=2$ and $x \geq 2$
we have $P_2(x)-P_1(x)=x-2 \geq 0$
(whatever the parity of $N$) and clearly both $P_2[0](x)=x>0$ and $P_2[1](x)=x-1>0$.
Assuming then that claim (v) holds for $i=j$, we get from \eqref{recur}, assuming
still $x \geq 2$,
  \begin{equation}
    \begin{aligned}
    P_{j+1}(x) - P_j(x) &= xP_j(x) - P_{j-1}(x) - P_j(x)
    = (x-1)P_j(x) - P_{j-1}(x) \\
    &\geq P_j(x) - P_{j-1}(x) \geq 0
    \end{aligned}
  \end{equation}
and therefore $P_{j+1}(x) > 0$ as well. We conclude that for every $i \geq 2$
we have $P_i(x)>0$ whenever $x \geq 2$, so all roots of $P_i$ lie to the left
of $2$, establishing the bound on $\djj{2}=\beta=\gamma_{n+1}$.

To show that for each fixed $N$ $\djj{j}[N]$ is a strictly increasing function of $j$
we claim that for any $n \geq 3$ and either parity of $N$
  \begin{enumerate}[(vi)]
    \item $P_{j+1}(\gamma_{n+1})-P_j(\gamma_{n+1})<0$ whenever $1 \leq j \leq n$,
  \end{enumerate}
which we prove by induction on $j$.
Since $P_2(x)-P_1(x)=x-2$ (whatever the parity of $N$),
the case $j=1$ follows immediately from the fact,
proved in the preceding paragraph,
that $\gamma_{n+1}<2$.
Assuming now that claim (vi) holds for a given $j \in \Z \cap [1,n-1]$
and applying \eqref{recur} with $i=j+1$
along with the same inequality $\gamma_{n+1}<2$, we get
  \begin{equation}
    P_{j+2}(\gamma_{n+1})-P_{j+1}(\gamma_{n+1})
    =
    (\gamma_{n+1}-1)P_{j+1}(\gamma_{n+1})-P_j(\gamma_{n+1})
    <
    P_{j+1}(\gamma_{n+1})-P_j(\gamma_{n+1})
    <
    0,
  \end{equation}
confirming (vi).
It then follows from \eqref{di} that
$1<\djj{2}<\djj{3}<\djj{4}<\cdots<\djj{n}$.

It remains to reintroduce the logarithmic terms.
Recalling the definitions \eqref{Aevendef} and \eqref{Aodddef},
for each integer $N \geq 2$ we define the function
$F_N: \R^{n-1} \to \R^{n-1}$ by
  \begin{equation}
    F_N
    \begin{pmatrix}
      x_2 \\
      x_3 \\
      \vdots \\
      x_n
    \end{pmatrix}
    :=
    A_N(x_2)
    \begin{pmatrix}
      x_2 \\
      x_3 \\
      \vdots \\
      x_n
    \end{pmatrix}
  \end{equation}
and calculate its derivative at $(\djj{i})_{i=2}^n$:
  \begin{equation}
    dF_N|_{(\djj{i})_{i=2}^n} = A_N(\djj{2}) +
    \begin{pmatrix}
      \djj{2} & 0 & 0 & \cdots & 0 \\
      \djj{3} & 0 & 0 & \cdots & 0 \\
      \vdots & \vdots & \vdots & \ddots & \vdots \\
      \djj{n} & 0 & 0 & \cdots & 0
    \end{pmatrix},
  \end{equation}
whose determinant is
  \begin{equation}
  \begin{aligned}
    &\det A_N(\djj{2}) + \det
    \begin{pmatrix}
      \djj{2} & -1 & 0 & 0 & \cdots & 0 \\
      \djj{3} & \djj{2} & -1 & 0 & \cdots & 0 \\
      \djj{4} & -1 & \djj{2} & -1 & \cdots & 0 \\
      \vdots & \vdots & \ddots & \ddots &\ddots &\vdots \\
      \djj{n-1} & 0 & \cdots & -1 & \djj{2} & -1 \\
      \djj{n} & 0 & \cdots & 0 & -2^{(N+1) \bmod 2} & \djj{2}-N \bmod 2
    \end{pmatrix} \\
    = &P_n(\djj{2}) + \sum_{i=2}^{n-1} \djj{i} P_{n-i+1}(\djj{2}) + \djj{n} > 0,
  \end{aligned}
  \end{equation}
using just expansion by minors along with the inequalities, proven above,
$\djj{i}>1>0$ for $2 \leq j \leq n$ and $P_j(\djj{2})>0$ for $1 \leq j \leq n$.
We conclude by invoking the inverse function theorem
and taking $m$ large in terms of $\djj{2},\ldots,\djj{n}$.
\end{proof}

Henceforth, given integers $N \geq 2$, $\ell \geq k \geq 1$, and $m \geq 1$,
along with $N-1$ real numbers $\zetaj{1}, \ldots, \zetaj{N-1}$,
we define
(i) $\cj{2}[N=2,k,\ell,m]:=0$, $\cj{2}[N=3,k,\ell,m]:=1$,
and $\{\cj{i}[N \geq 4, k, \ell, m]\}_{i=2}^n$
as in the proof of Lemma \ref{cilemma},
(ii) $\taubarj{1}[N,k,\ell,m]$ as defined by \eqref{taubar1},
(iii) $\taubarj{i}[N,k,\ell,m]:=\cj{i}[N,k,\ell,m] \taubarj{1}[N,k,\ell,m]$ for $2 \leq i \leq n$,
(iv) $\taubarj{i}[N,k,\ell,m]$ and $\cj{i}[N,k,\ell,m]$ for $n+1 \leq i \leq N-1$
(and $\cj{1}[N,k,\ell,m]=1$)
in accordance with \eqref{taubarrefl} and \eqref{cidef},
and (v) $\tauj{i}[N,k,\ell,m,\zeta]$ for $1 \leq i \leq N$ as defined by \eqref{taudef}.
These quantities
complete the specification of the initial surfaces defined in \eqref{initsurfdef}.

\subsection*{Parameter dependence}
The $\zeta$ and $\xi$ parameters influence the forces,
which will be analyzed later to manage the substitute kernel
(described in Section \ref{intro} and formally introduced in Section \ref{linop}),
as well as the \emph{dislocations}
  \begin{equation}
  \label{Ddef}
    \Dcal_i
    =
    \Dcal_i[N,k,\ell,m,\zeta,\xi]
    :=
    \begin{cases}
      \frac{1}{2}\tauj{i}\xij{i} - \frac{1}{2}\tauj{i-1}\xij{i-1}
        \mbox{ for } 2 \leq i \leq N-1 \\
      0 \mbox{ for $i=1$ and $i=N$},
    \end{cases}
  \end{equation}
which will be used to manage the extended part of the extended substitute kernel
(again see Sections \ref{intro} and \ref{linop})
and each of which measures the antisymmetric part of the vertical displacement of a pair of adjacent
inequivalent (under the action of $\Grp$ \eqref{Grpdef}) catenoidal regions
relative to the toral region they share.
Morally, the next lemma ensures
that, by adjusting the parameters,
we can freely prescribe any set of suitably bounded forces and dislocations.
Indeed this surjectivity assertion could be stated precisely
and proved as a corollary of the lemma
by applying the Brouwer fixed-point theorem.
Because, however, we will also need to allow graphical deformations of the initial surfaces,
we bypass this step and instead will more directly apply the lemma
to control the extended substitute kernel
in the proof of Theorem \ref{mainthm}.

\begin{lemma}[Parametric dependence of the forces and dislocations]
\label{FD}
Given $c>0$ and integers $N \geq 2$ and $\ell \geq k \geq 1$,
there exist real numbers $C=C[N,k,\ell], m_0=m_0[N,k,\ell,c]>0$
and an invertible linear map
$\Theta=\Theta[N,k,\ell,m]: \R^{2N-2} \to \R^{2N-2}$
such that whenever $m>m_0$
  \begin{equation}
  \tag{i}
    \norm{\Theta}+\norm{\Theta^{-1}} \leq C[N,k,\ell],
  \end{equation}
where $\norm{\cdot}$ denotes the operator norm
induced by the Euclidean norm on $\R^{2N-2}$, and
whenever $\zeta,\xi \in [-c,c]^{2N-2}$
  \begin{equation}
  \tag{ii}
    \abs{
      \tauj{1}^{-1}
      \begin{pmatrix}
        m^2\Fcal_1[N,k,\ell,m,\zeta,\xi] \\
        \vdots \\
        m^2\Fcal_N[N,k,\ell,m,\zeta,\xi] \\
        \Dcal_2[N,k,\ell,m,\zeta,\xi] \\
        \vdots \\
        \Dcal_{N-1}[N,k,\ell,m,\zeta,\xi]
      \end{pmatrix}
      -\Theta
      \begin{pmatrix}
        \zetaj{1} \\
        \vdots \\
        \zetaj{N-1} \\
        \xij{1} \\
        \vdots \\
        \xij{N-1}
      \end{pmatrix}
    }
    \leq
    C[N,k,\ell],
  \end{equation}
where $\abs{\cdot}$ denotes the Euclidean norm on $\R^{2N-2}$.
\end{lemma}

We emphasize that the estimates made in Lemma \ref{FD}
are independent of $m$ and the size of the parameters:
$C[N,k,\ell]$ does not depend on $m$ or $c$ (but $m_0[N,k,\ell,c]$ does depend on $c$).
This independence will be crucial 
in the proof of the main theorem
when defining
the nonlinear map whose fixed point will give us our final minimal surface.
In particular it will be needed to establish
that we can choose the parameter factor of the domain of this map to be compact,
an ingredient in the justification of the applicability of the Schauder fixed-point theorem
in the proof of Theorem \ref{mainthm}.
(Roughly, the independence from $c$ of these estimates ensures that attempts
to control the extended substitute kernel (via the forces and dislocations)
by varying the parameters
will not drive any of these parameters off to infinity.)

\begin{proof}
To begin, it is easy to see that the linear map $T=T[N,k,\ell]: \R^{N-1} \to \R^{N-1}$
defined by
  \begin{equation}
    T:
    \begin{pmatrix}
      m^2\Fcal_1[N,k,\ell,m,\zeta,\xi] \\
      m^2\Fcal_2[N,k,\ell,m,\zeta,\xi] \\
      \vdots \\
      m^2\Fcal_{N-1}[N,k,\ell,m,\zeta,\xi] \\
      m^2\Fcal_N[N,k,\ell,m,\zeta,\xi] \\
      \Dcal_2[N,k,\ell,m,\zeta,\xi] \\
      \vdots \\
      \Dcal_{N-1}[N,k,\ell,m,\zeta,\xi]
    \end{pmatrix}
    \mapsto
    \begin{pmatrix}
      \frac{k \ell m^2}{2\pi}(\Fcal_1-\Fcal_2) + 4\pi(\Dcal_1+\Dcal_2) \\
      \frac{k \ell m^2}{2\pi}(\Fcal_2-\Fcal_3) + 4\pi(\Dcal_2+\Dcal_3) \\
      \vdots \\
      \frac{k \ell m^2}{2\pi}(\Fcal_{N-1}-\Fcal_N) + 4\pi(\Dcal_{N-1}+\Dcal_N) \\
      \frac{k\ell m^2}{8\pi^2} \Fcal_1 \\
      2\Dcal_2 \\
      \vdots \\
      2\Dcal_{N-1}
    \end{pmatrix}
  \end{equation}
is invertible
(since, the lowest $N-2$ components on the right
determine all dislocations, which with the top $N$
components then determine all the forces too)
with inverse bounded independently of $m$ and $c$.
To prove the lemma it will therefore suffice to identify
an invertible linear map
$\widetilde{\Theta}=\widetilde{\Theta}[N,k,\ell,m]: \R^{2N-2} \to \R^{2N-2}$
such that
  \begin{equation}
  \label{Thestimates}
    \norm{\widetilde{\Theta}}+\norm{\widetilde{\Theta}^{-1}} \leq C[N,k,\ell]
    \qquad \mbox{and} \qquad
    \abs{
      \tauj{1}^{-1}T\begin{pmatrix} m^2\Fcal \\ \Dcal \end{pmatrix}
      -\widetilde{\Theta}\begin{pmatrix} \zeta \\ \xi \end{pmatrix}
    }
    \leq
    C[N,k,\ell]
  \end{equation}
for some constant $C[N,k,\ell]>0$ whenever
$\zeta,\xi \in [-c,c]^{N-1}$
and $m$ is sufficiently large in terms of $N$, $k$, $\ell$, and $c$.

In fact we will show that we can take
  \begin{equation}
  \label{Thetat}
    \widetilde{\Theta}\begin{pmatrix} \zeta \\ \xi \end{pmatrix}
    :=
    \begin{pmatrix} Z & 0 \\ B & \Xi \end{pmatrix}
      \begin{pmatrix} \zeta \\ \xi \end{pmatrix},
  \end{equation}
where $0$ is the $(N-1) \times (N-1)$ zero matrix,
$B=B[N,k,\ell,m]$ is an $(N-1) \times (N-1)$ matrix bounded independently of $c$ and $m$,
and $Z=Z[N,k,\ell,m]$ and $\Xi=\Xi[N,k,\ell,m]$ are the $(N-1) \times (N-1)$ matrices
  \begin{equation}
  \label{ZXidef}
    Z
    :=
    \begin{pmatrix}
      8\pi\cj{1} & -\cj{2} & 0 & \cdots & 0 & 0 \\
      8\pi\cj{2} & \cj{1}+\cj{3} & -\cj{3} & 0 & \cdots & 0 \\
      8\pi\cj{3} & -\cj{2} & \cj{2}+\cj{4} & -\cj{4} &  \ddots & \vdots \\
      8\pi\cj{4} & 0 & \ddots & \ddots & \ddots & 0 \\
      \vdots & \vdots & \ddots & -\cj{N-3} & \cj{N-3}+\cj{N-1} & -\cj{N-1} \\
      8\pi\cj{N-1} & 0 & \cdots & 0 & -\cj{N-2} & \cj{N-2}+\cj{N}
      \end{pmatrix}
  \end{equation}
and
  \begin{equation}
  \label{Xidef}
      \Xi
      :=
      \begin{pmatrix}
        \cj{1} & 0 & 0 & \cdots & 0 & 0 \\
        -\cj{1} & \cj{2} & 0 & \cdots & 0 & 0 \\
        0 & -\cj{2} & \cj{3} & \ddots & 0 & 0 \\
        \vdots & \ddots & \ddots & \ddots & \ddots & \vdots \\
        0 & 0 & \ddots & -\cj{N-3} & \cj{N-2} & 0 \\
        0 & 0 & \cdots & 0 & -\cj{N-2} & \cj{N-1}
      \end{pmatrix},
  \end{equation}
recalling \eqref{cisummary} and understanding $\cj{N}:=0$.
(We acknowledge that $\cj{1}=\cj{N-1}=1$
but refrain from making these substitutions above
so as to avoid obscuring the structure of the matrices.)
We will now verify that $\widetilde{\Theta}$ so defined
satisfies \eqref{Thestimates},
identifying the matrix $B$ along the way.

First we check that $Z$ and $\Xi$ are invertible.
Invertibility of $\Xi$ is obvious,
since it is lower-triangular
with all its diagonal entries strictly positive.
Next, we inductively alter the middle $N-3$ columns of $Z$,
starting with column $N-2$ and working our way to the left until column $2$,
by replacing each by its sum with the column immediately to its right;
we also divide the first column by $8\pi$.
The resulting matrix (for the computation of which we recall that $\cj{N}=0$)
  \begin{equation}
    \widetilde{Z}
    :=
    \begin{pmatrix}
      1 & -\cj{2} & 0 & \cdots & 0 & 0 \\
      \cj{2} & \cj{1} & -\cj{3} & 0 & \cdots & 0 \\
      \cj{3} & 0 & \cj{2} & -\cj{4} &  \ddots & \vdots \\
      \cj{4} & 0 & \ddots & \ddots & \ddots & 0 \\
      \vdots & \vdots & \ddots & 0 & \cj{N-3} & -\cj{N-1} \\
      \cj{N-1} & 0 & \cdots & 0 & 0 & \cj{N-2}
    \end{pmatrix}
  \end{equation}
is invertible if and only if $Z$ is.
The only nonzero entries of $\widetilde{Z}$ lie in (i) the first column,
whose entries are all strictly positive,
(ii) the diagonal, whose entries are also all strictly positive,
and (iii) the superdiagonal, whose entries are all strictly negative.
It is easy to see that any square matrix of this form
has strictly positive determinant.
For example we can compute the determinant by cofactor expansion along the bottom row.
Starting with the entry in the bottom row and first column
($\cj{N-1}$ in $\widetilde{Z}$),
we see that the submatrix obtained by deleting the bottom row and first column
is lower-triangular with all diagonal entries strictly negative.
Counting minus signs,
including possibly one contributed by the position
of the entry in question,
we find that the corresponding cofactor is strictly positive
and, the entry itself being strictly positive as well,
therefore the corresponding term in the expansion is also strictly positive.
The remaining term,
corresponding to the entry in the bottom row and last column
($\cj{N-2}$ in $\widetilde{Z}$)
is clearly the product of a strictly positive number, the entry itself,
with the determinant of a smaller matrix of the same form under consideration.
Since our claim is obvious in the $1 \times 1$ case, the general case now follows by induction.
Thus we see that $Z$ is invertible as well.
Consequently, for any choice of $B$ (to be identified shortly)
in \eqref{Thetat}, the map $\widetilde{\Theta}$ is indeed invertible,
and the first inequality in \eqref{Thestimates}
is now ensured by \eqref{cisummary}
(and the bound $\norm{B} \leq C[N,k,\ell]$ established below).

Now we estimate
$\frac{k\ell m^2}{2\pi \tauj{1}}(\Fcal_i-\Fcal_{i+1})+\frac{4\pi}{\tauj{1}}(\Dcal_i+\Dcal_{i+1})$
for $1 \leq i \leq N-1$.
When $N=2$, from \eqref{Fn} and \eqref{Fn+1} we find
  \begin{equation}
    \abs{\Fcal_1-\Fcal_2 - \left(4\pi\tauj{1}
      -\frac{16\pi^2}{k\ell m^2}\tauj{1}\ln \frac{1}{10\ell m\tauj{1}}\right)}
    \leq
    Cm^{-2}\tauj{1},
  \end{equation}
whence with \eqref{taudef}
  \begin{equation}
    \abs{\Fcal_1-\Fcal_2 - 2e^{\zetaj{1}}\left(2\pi\taubarj{1}
      -\frac{8\pi^2}{k\ell m^2}\taubarj{1}\left(-\zetaj{1}+\ln \frac{1}{10\ell m\taubarj{1}}\right)\right)}
    \leq
    Cm^{-2}\tauj{1},
  \end{equation}
so by the balancing condition \eqref{balancing}
  \begin{equation}
    \abs{\Fcal_1-\Fcal_2 - \frac{16\pi^2\tauj{1}}{k\ell m^2}\zetaj{1}}
    \leq
    Cm^{-2}\tauj{1},
  \end{equation}
while $\Dcal_1=\Dcal_2=0$ by \eqref{Ddef},
which proves 
  \begin{equation}
  \label{Zest1}
    \abs{\frac{k \ell m^2}{2\pi \tauj{1}}(\Fcal_1-\Fcal_2) + \frac{4\pi}{\tauj{1}}(\Dcal_1+\Dcal_2)
      -8\pi\zetaj{1}}
    \leq
    C \mbox{ when $N=1$},
  \end{equation}
taking the constant $C$ possibly larger (depending on $k$ and $\ell$) than the one immediately above
but still independent of $c$ and $m$.

When $N \geq 3$
and $1 \leq i \leq N-1$, from \eqref{forcediffs} and \eqref{Ddef} we find
  \begin{equation}
  \label{Fdiff}
    \begin{aligned}
      &\left|\frac{k \ell m^2}{2\pi\tauj{1}}\left(\Fcal_i-\Fcal_{i+1}\right)
        +\frac{4\pi}{\tauj{1}}\left(\Dcal_i+\Dcal_{i+1}\right)\right. \\
      &\;\;\;\;
      \left. -\frac{k\ell m^2}{2\pi\tauj{1}}\left(2\pi\left(-\tauj{i-1}+2\tauj{i}-\tauj{i+1}\right)
      -\frac{16\pi^2}{k\ell m^2} \tauj{i} \ln \frac{1}{10\ell m\tauj{i}}\right)\right|
      \leq
      C,
    \end{aligned}
  \end{equation}
understanding $\tauj{0}=\tauj{N}:=0$.
Using \eqref{taubardef}, \eqref{taudef}, and Taylor expansion,
for $2 \leq j \leq N-1$ we have
  \begin{equation}
  \label{Taylorest}
    \abs{\tauj{j}
      -e^{\zetaj{1}}\taubarj{j}-k^{-1}\ell^{-1}m^{-2}\cj{j}\zetaj{j}\tauj{1}}
    \leq
    m^{-4}c^2e^{2c/m^2}\cj{j}\tauj{1}
    \leq
    Cm^{-2}\tauj{1},
  \end{equation}
where to ensure the last inequality we take $m$ sufficiently large in terms of $c$.
Thus, when $3 \leq i \leq N-1$,
  \begin{equation}
  \label{Fdiff1st}
    \begin{aligned}
      &\left|\left(-\tauj{i-1}+2\tauj{i}-\tauj{i+1}\right)
        -e^{\zetaj{1}}\left(-\taubarj{i-1}+2\taubarj{i}-\taubarj{i+1}\right)\right. \\
      &\;\;\;\;
      \left. -\frac{\tauj{1}}{k \ell m^{2}}
          \left(-\cj{i-1}\zetaj{i-1}+2\cj{i}\zetaj{i}-\cj{i+1}\zetaj{i+1}\right)\right|
      \leq
      Cm^{-2}\tauj{1}
    \end{aligned}
  \end{equation}
and, further using \eqref{taubardef} and \eqref{taudef}, when $2 \leq i \leq N-1$
  \begin{equation}
  \label{Fdiff2nd}
    \abs{\tauj{i} \ln \frac{1}{10\ell m \tauj{i}}
      - \left(e^{\zetaj{1}}\taubarj{i} \ln \frac{1}{10\ell m \taubarj{i}}
      - \tauj{1}\cj{i}\zetaj{1}
      +\frac{\tauj{1}\cj{i}\zetaj{i}}{4\pi}
      -\frac{\tauj{1}\cj{2}\cj{i}\zetaj{i}}{8\pi}\right)}
    \leq
    C\tauj{1}
  \end{equation}
provided $m$ is sufficiently large in terms of $c$.
Note that by virtue of \eqref{taubardef} and \eqref{barbalance}
  \begin{equation}
  \label{Fdiffsupp}
    \abs{\cj{2}\cj{i}-(\cj{i-1}+\cj{i+1})} \leq Cm^{-2}
    \mbox{ for $1 \leq i \leq N-1$}
  \end{equation}
(understanding $\cj{0}=\cj{N}=0$),
so the last term on the left-hand side of \eqref{Fdiff2nd}
may be replaced by $-\frac{\tauj{1}(\cj{i-1}+\cj{i+1})\zetaj{i}}{8\pi}$.
By applying estimates \eqref{Fdiff1st}--\eqref{Fdiffsupp}
to \eqref{Fdiff}
and imposing the balancing condition \eqref{balancing}
we obtain
  \begin{equation}
  \label{Zest2}
    \abs{\frac{k\ell m^2}{2\pi \tauj{1}}(\Fcal_i-\Fcal_{i+1})
      +\frac{4\pi}{\tauj{1}}(\Dcal_i+\Dcal_{i+1})
      -\left(
        8\pi\cj{i}\zetaj{1}-\cj{i-1}\zetaj{i-1}+(\cj{i-1}+\cj{i+1})\zetaj{i}-\cj{i+1}\zetaj{i+1}
      \right)
    }
    \leq
    C
  \end{equation}
when $N \geq 3$ and $3 \leq i \leq N-1$
(understanding $\cj{N}:=0$).

In the remaining cases that $N \geq 3$ but $i \in \{1,2\}$
a similar computation
(using simply $\tauj{1}=e^{\zetaj{1}}\taubarj{1}$
in place of \eqref{Taylorest} when $j=1$)
reveals
  \begin{equation}
  \label{Zest3}
    \begin{aligned}
      &\abs{
      \frac{k\ell m^2}{2\pi \tauj{1}}(\Fcal_1-\Fcal_2)
        +\frac{4\pi}{\tauj{1}}(\Dcal_1+\Dcal_2)
      -\left(
          8\pi\cj{1}\zetaj{1} - \cj{2}\zetaj{2}
        \right)
      }
      \leq
      C \mbox{ and} \\
      &\abs{
      \frac{k\ell m^2}{2\pi \tauj{1}}(\Fcal_2-\Fcal_3)
          +\frac{4\pi}{\tauj{1}}(\Dcal_2+\Dcal_3)
      -\left(
        8\pi\cj{2}\zetaj{1}+(\cj{1}+\cj{3})\zetaj{2}-\cj{3}\zetaj{3}
      \right)
      }.
    \end{aligned}
  \end{equation}
Together, \eqref{Zest1}, \eqref{Zest2}, and \eqref{Zest3} show that
  \begin{equation}
  \label{Zest}
    \abs{
       \frac{k \ell m^2}{2\pi\tauj{1}}(\Fcal_i-\Fcal_{i+1})
        + \frac{4\pi}{\tauj{1}}(\Dcal_i+\Dcal_{i+1}) - (Z\zeta)_i
    }
    \leq
    C
  \end{equation}
for some constant $C=C[N,k,\ell]>0$
whenever $N \geq 2$,
$1 \leq i \leq N-1$,
$\zeta,\xi \in [-c,c]^{N-1}$,
and $m$ is sufficiently large in terms of $c$.

Next, it is obvious from \eqref{Ddef},
using \eqref{Taylorest} as necessary
and continuing to take $m$ large in terms of $c$,
that
  \begin{equation}
  \label{Xiest1}
    \abs{\frac{2}{\tauj{1}}\Dcal_i - \left(\cj{i} \xij{i} - \cj{i-1}\xij{i-1}\right)} \leq C
    \mbox{ for } 2 \leq i \leq N-1.
  \end{equation}
Furthermore, from \eqref{torheights}, \eqref{forcecomp}, and \eqref{balancing},
using \eqref{Fdiff2nd} and \eqref{Fdiffsupp} again
and still taking $m$ large in terms of $c$,
we have
  \begin{equation}
  \label{Xiest2}
    \begin{aligned}
      &\abs{
        \frac{k\ell m^2}{8\pi^2 \tauj{1}}\Fcal_1
        -\left(
           \cj{1}\xij{1}+2^{N \bmod 2} \zetaj{1}
         \right) 
      }
      \leq
      C
      \mbox{ if } 2 \leq N \leq 3 \mbox{ and} \\
      &\left|
      \frac{k\ell m^2}{8\pi^2 \tauj{1}}\Fcal_1
      -\left[
        \cj{1}\xij{1}
          +\left(2^{N \bmod 2}\cj{n}+2\sum_{j=1}^{n-1} \cj{j}\right)\zetaj{1} \right.\right. \\
      &\;\;\;\; \left. \left.
          +\frac{1}{4\pi}\sum_{j=2}^{n-1} \left(\cj{j-1}-2\cj{j}+\cj{j+1}\right)\zetaj{j}
          +\frac{2^{N \bmod 2}}{8\pi} \left(\cj{n-1}-2\cj{n}+\cj{n+1}\right)\zetaj{n}
      \right]
      \right|
      \leq C
      \mbox{ if $N \geq 3$}.
    \end{aligned}
  \end{equation}
Together, \eqref{Zest}, \eqref{Xiest1}, and \eqref{Xiest2}
establish the second inequality of \eqref{Thestimates}
with $\widetilde{\Theta}$ defined by \eqref{Thetat},
the entries of the first row of the matrix $B$
given by the coefficients of the components of $\zeta$
appearing in \eqref{Xiest2},
and the remaining entries of $B$ vanishing.
This completes the proof.
\end{proof}

\section{Estimates of the initial geometry}

\subsection*{Norms and spaces of sections}
To state the estimates for the geometry of the initial surfaces
and to carry out the rest of the construction
we must first identify certain norms and corresponding spaces of sections.
For the most part our notation is standard and speaks for itself.
Given a smooth vector bundle $E$ over a smooth manifold $M$ (possibly with boundary),
a nonnegative integer $j$, and a real number $\alpha \in (0,1)$,
we write $C^j_{loc}(E)$ and $C^{j,\alpha}_{loc}(E)$
for the space of sections of $E$ having
component functions of class $C^j_{loc}$ or $C^{j,\alpha}_{loc}$ respectively
relative to any (so every) smooth local chart and smooth trivialization;
we set $C^\infty(E):=\bigcap_{j=0}^\infty C^j_{loc}(E)$.
When $E$ is the trivial bundle $M \times \R$,
we write simply $M$ in place of $M \times \R$
in our notation for the spaces just defined and also for the spaces below,
as is standard for spaces of real-valued functions.

All of the vector bundles of interest to us are derived
from tangent bundles by a combination of
duality, tensor product, pullback, and projection;
a Riemannian metric on $M$
or on another manifold in which it is immersed
will determine canonical
metrics and connections on all these bundles.
When there is no danger of confusion,
we write simply $\abs{\cdot}$ for the corresponding pointwise norm
and $D$ for the connection.
Given a section $u$
of a bundle $E$ over $M$
thus equipped,
we define the standard global norms
  \begin{equation}
    \norm{u}_j
    =
    \norm{u}_{C^j(E,g)}
    =
    \norm{u: C^j(E,g)}
    :=
    \sum_{i=0}^j \sup_{p \in M} \abs{D^iu(p)}
  \end{equation}
as well as the H\"{o}lder seminorm
  \begin{equation}
   [u]_\alpha
    :=
    \sup_{\gamma: [0,1] \to M}
      \frac{\abs{u(\gamma(1))-P_0^1[\gamma]u(\gamma(0))}}{\abs{\gamma}^\alpha},
  \end{equation}
where the supremum is taken over all piecewise $C^1$ paths,
$\abs{\gamma}$ denotes the length of such a path,
and $P_0^1[\gamma]: E_{\gamma(0)} \to E_{\gamma(1)}$
is the parallel transport map along $\gamma$
from the fiber over $\gamma(0)$ to the fiber over $\gamma(1)$.

Then we can define also the H\"{o}lder norms
  \begin{equation}
    \norm{u}_{j,\alpha}
    =
    \norm{u: C^{j,\alpha}(E,g)}
    :=
    \norm{u}_j + \left[ D^j u \right]_\alpha.
  \end{equation}
Note that for functions on convex open subsets of Euclidean space
these H\"{o}lder norms agree with the conventional ones.
Generally, the spaces $C^{j,\alpha}(E,g)$ and $C^j(E,g)$ consisting
of sections with finite corresponding norm
enjoy many of the properties familiar from the Euclidean case;
in particular $C^{j,\beta}(E,g)$ embeds compactly in $C^{j,\alpha}(E,g)$
whenever $M$ is compact and $0<\alpha<\beta<1$.
In the construction we will routinely wish to compare norms of the above type
induced by different metrics on a single manifold. The definitions make it easy
to see that
$\norm{u: C^{j,\alpha}(E,h)} \leq C \norm{u: C^{j,\alpha}(E,g)}$,
where $C$ is controlled by the $g$-norms of $h$, its inverse,
and finitely many $g$-derivatives of $h$
(the maximum order needed depending in a transparent way on $j$, $\alpha$,
and on the bundle $E$);
of course if $M$ is not compact, then $C$ may blow up,
depending on $g$ and $h$.

If $M$ is a two-sided hypersurface immersed in a Riemannian manifold $N$
and $\Grp_M$ is a group of isometries of $N$ preserving $M$ as a set,
then $\Grp_M$ acts on a section $u$
of the normal bundle of $M$ by
$(\mathfrak{g}.u)(p):=\mathfrak{g}_*[u \left(\mathfrak{g}^{-1}(p) \right)]$
for each element $\mathfrak{g}: N \to N$ of $\Grp_M$.
Because this bundle is just the trivial $\R$
bundle over $M$, its sections can be identified with functions
(for us representing mean curvature or normal perturbations)
on which the corresponding action of $\Grp_M$ is given by
  \begin{equation}
  \label{Grpaction}
    (\mathfrak{g}.f)(p)
    :=
    \begin{cases}
      f\left(\mathfrak{g}^{-1}(p)\right) \mbox{ if $\mathfrak{g}$ preserves each side of $M$} \\
      -f\left(\mathfrak{g}^{-1}(p)\right) \mbox{ if $\mathfrak{g}$ reverses the sides of $M$}.
    \end{cases}
  \end{equation}
All the elements of the symmetry group
$\Grp=\Grp[k,\ell,m]$ (defined in \eqref{Grpdef})
of the construction
can be seen to fix each side of the initial surface
$\Sigma=\Sigma[N,k,\ell,m,\zeta,\xi]$ (defined in \eqref{initsurfdef}),
so its action is always given by the first line of \eqref{Grpaction}.
(As mentioned earlier (Remark \ref{verticalsymremark}),
when $k=\ell$, we have the option of enforcing a larger symmetry group
in the construction, one admitting reflections through great circles
on $\T$. When $N$ is odd, such reflections reverse the sides
of the initial surface $\Sigma$.)

\begin{notation}
\label{functionspaceappendages}
In general,
if $M$ is a two-sided hypersurface immersed in a Riemannian manifold $N$
and $\Grp_M$ is a group of isometries of $N$ preserving $M$ as a set,
we will append the subscript $\Grp_M$ to a space of functions
to designate the subspace
consisting of functions which are equivariant
under the $\Grp_M$ action \eqref{Grpaction}.
\end{notation}
Finally we will often wish to work with weighted versions of the above norms.
For this construction the following definition suffices:
  \begin{equation}
  \label{weightedhoelderdef}
    \norm{u: C^{j,\alpha}(E,g,f)}
    :=
    \sup_{p \in M}
      \frac{\norm{u: C^{j,\alpha}\left(E|_{B[p,1,g]},g\right)}}{f(p)},
  \end{equation}
where $f: M \to (0,\infty)$ is a given weight function
and $B[p,1,g] \subseteq M$ is
the $g$ metric ball of radius $1$ centered at $p \in M$.
We will also make use
of weighted $C^j$ norms, with the obvious definition.

\subsection*{The $\chi$ metric}
It is the primary task of this section to estimate the intrinsic
and extrinsic geometry of the initial surfaces.
We continue to write
$\iota=\iota[N,k,\ell,m,\zeta,\xi]: \Sigma[N,k,\ell,m,\zeta,\xi] \to \Sph^3$
for the inclusion map of $\Sigma$ in $\Sph^3$,
and we define
  \begin{equation}
  \label{gdef}
    g=g[N,k,\ell,m,\zeta,\xi]:=\iota^*\gsph,
  \end{equation}
the metric induced on $\Sigma$ by $\gsph$ and $\iota$.
To fix the extrinsic quantities
we pick on each initial surface $\Sigma$
the global unit normal $\nu=\nu[N,k,\ell,m,\zeta,\xi]$ which is
directed toward $C_1$ at the points of $\Sigma$ closest to $C_1$.
We then define
  \begin{equation}
  \label{AHdef}
    A=A[N,k,\ell,m,\zeta,\xi]:=(\gsph \circ \iota)(Dd\iota,\nu)
    \qquad \mbox{and} \qquad
    H=H[N,k,\ell,m,\zeta,\xi]:=\tr_{\iota^*g}A,
  \end{equation}
respectively
the scalar-valued second fundamental form and mean curvature
of $\Sigma$ relative to $\nu$ and $\gsph$,
$D$ being the connection induced on $T^*\Sigma \otimes \iota^*T\Sph^3$ by $\gsph$ and $\iota$.
In particular $H=(\gsph \circ \iota)(\mathbf{H},\nu)$,
recalling the vector-valued mean curvature $\mathbf{H}$
defined below \eqref{forcedef}.

Every initial surface admits by virtue of its construction
a natural decomposition into overlapping regions,
each of which resembles either a portion of a torus coaxial with $\T$
or (via $\Phi$ \eqref{Phidef}) a truncated catenoid.
Modulo the horizontal symmetries, there are $N$ such toral regions, one
for each torus incorporated in the construction, and there are $N-1$
catenoidal regions, one for each pair of adjacent tori.
Definitions are made in the subsections below. The estimates
will then be obtained by treating the catenoidal regions as perturbations of
Euclidean catenoids and the toral regions as graphs over the Clifford torus.

Because all these regions shrink with increasing $m$ and because even on a fixed
initial surface the characteristic scale $m^{-1}$ of the toral regions
dwarfs the characteristic scale $\tauj{1}$ near the waists,
it will be advantageous to uniformize the problem (and flatten $\Sigma$)
by working
with a metric $\chi$ on each initial surface
conformal to the natural one $g=\iota^*\gsph$.
We will set
  \begin{equation}
  \label{chidef}
    \chi=\chi[N,k,\ell,m,\zeta,\xi]:=\rho^2 g,
  \end{equation}
where the conformal factor 
$\rho=\rho[N,k,\ell,m,\zeta,\xi]: \Sigma \to \R$
is defined so that its reciprocal
$\rho^{-1}$ is
a $\Grp$-equivariant function (i) measuring
on each catenoidal region the ${\Phi^{-1}}^*\geuc$ distance to the axis
and (ii) transitioning smoothly to the constant $m^{-1}$ by the edge of
the toral regions.

To be precise, we first define
$\widetilde{\rho}[\zzj{1},\zzj{2}]: \R^3 \to \R$,
for given $\zzj{1}<\zzj{2} \in \R$, by
  \begin{equation}
    \widetilde{\rho}[\zzj{1},\zzj{2}](\xx,\yy,\zz)
    :=
    \begin{cases}
      \left(\left(\xx^2+\yy^2\right)^{-1/2}-m\right)
        \cutoff{(5 \ell m)^{-1}}{(10 \ell m)^{-1}}
        \left(\sqrt{\xx^2+\yy^2}\right)
        \mbox{ if } \zzj{1} \leq \zz \leq \zzj{2} \\
      0 \mbox{ otherwise }
    \end{cases},
  \end{equation}
recalling \eqref{cutoff},
and $\widehat{\rho}: \R^3 \to \R$ by
  \begin{equation}
    \widehat{\rho}(\xx,\yy,\zz)
    :=
    m + \sum_{i=0}^{N-2} \sum_{(\xx_0,\yy_0) \in \lhat_{i \bmod 2,i \bmod 2}}
      \widetilde{\rho}\left[\zzj{i}^K,\zzj{i+2}^K\right]
      \left(\xx - \xx_0, \yy - \yy_0 ,\zz\right),
  \end{equation}
recalling \eqref{catheights} and \eqref{lattices} and
taking $\zzj{0}^K:=-\frac{\pi}{5}$ and $\zzj{N}^K:=\frac{\pi}{5}$.
Then, recalling \eqref{Phidef},
$\rho \in C_{\Grp}^\infty(\Sigma)$
is uniquely defined by requiring
  \begin{equation}
  \label{rhodef}
    \rho \circ \Phi
    =
    \widehat{\rho} \quad \mbox{on $\Phi^{-1}(\Sigma)$}.
  \end{equation}

Evidently from \eqref{cisummary}, \eqref{taudef}, and \eqref{initsurfdef}
  \begin{equation}
  \label{rhobounds}
    \frac{m}{C[N,k,\ell]} \leq \rho \leq \frac{C[N,k,\ell]}{\tauj{1}}
  \end{equation}
for some constant $C[N,k,\ell]>0$ whenever $\zeta,\xi \in [-c,c]^{N-1}$
and $m$ is sufficiently large in terms of $N$, $k$, $\ell$, and $c$.
Equipped with the $\chi$ metric, each catenoidal region
tends with large $m$ to the flat cylinder of radius $1$,
while each toral region tends, away from the catenoids adjoining it,
to a flat $\frac{\sqrt{2}\pi}{k} \times \frac{\sqrt{2}\pi}{\ell}$
rectangle. Details are provided in the next two subsections.

Before proceeding, 
we briefly mention
a couple differences of our approach from \cite{KY}.
First, our catenoidal and toral regions (defined below)
correspond to their \emph{extended standard regions},
but since we never view their \emph{standard regions}
or \emph{transition regions} in isolation,
we omit the modifier \emph{extended}.
Second, our use of the $\chi$ metric follows theirs
to study the mean curvature equation on the initial surfaces globally,
but whereas Kapouleas and Yang introduce another metric
(the $h$ metric) conformal to $g$
in order to analyze the approximate kernel,
we will apply the $\chi$ metric to this problem as well, in the next section.

\subsection*{Catenoidal regions}
We define the standard cylinder
  \begin{equation}
    \label{cyldef}
    \cyl:=\R \times \Sph^1,
  \end{equation}
where $\Sph^1:=\{z \in \C \; : \; \abs{z}=1\}$.
We write $t$ for the standard coordinate on the $\R$ factor
and $\theta$ for the standard coordinate
on the universal cover $\R$ of the $\Sph^1$ factor given by
$\theta \mapsto e^{i\theta}$.
We will routinely and implicitly define functions on $\cyl$
by defining functions on this universal cover that are
invariant under the deck transformations.
We equip $\cyl$ with its usual flat metric
  \begin{equation}
  \label{chihatK}
    \chihatK:=dt^2+d\theta^2
  \end{equation}
and we define the embedding
  \begin{equation}
  \label{stdcat}
     \cate: \cyl \to \R^3
     \qquad \mbox{by} \qquad
     \cate(t,\theta)=(\cosh t \, \cos \theta, \cosh t \, \sin \theta, t),
  \end{equation}
whose image is the catenoid with waist radius $1$, axis of symmetry
the $\zz$-axis, and horizontal plane of the symmetry the $\zz=0$ plane.
It is easy to see that $\cate$ is conformal, with
  \begin{equation}
  \label{catconf}
    \cate^*\geuc = \cosh^2 t \, \chihatK,
  \end{equation}
$\geuc$ being the Euclidean metric on $\R^3$,
and that its unit normal pointing outward at the waist is
  \begin{equation}
  \label{catenormal}
    \widehat{\nu}
    =
    \sech t \, \cos \theta \, \partial_{\xx}
    +\sech t \, \sin \theta \, \partial_{\yy}
    -\tanh t \, \partial_{\zz}.
  \end{equation}
Of course the catenoid is famously minimal,
and it is elementary to check that,
more specifically,
it has second fundamental form (relative to $\widehat{\nu}$)
  \begin{equation}
  \label{cateA}
    \widehat{A}
    =
    dt^2 - d\theta^2.
  \end{equation}

Given $a>0$, we define the truncated cylinder
  \begin{equation}
    \label{cyladef}
    \cyl_a := [-a,a] \times \Sph^1,
  \end{equation}
and, given also $(\xx_{0},\yy_{0},\zzj{0}) \in \R^3$ and $\tau>0$,
we define the embedding
  \begin{equation}
  \label{catedef}
    \begin{aligned}
      &\cate[(\xx_{0},\yy_{0},\zzj{0}),\tau,a]: \cyl_a \to \R^3 \mbox{ by} \\
      &\cate[(\xx{_0},\yy{_0},\zzj{0}),\tau,a](t,\theta)
      :=
      (\xx_{0},\yy_{0},\zzj{0}) + \tau (\cosh t \, \cos \theta, \cosh t \, \sin \theta,t),
    \end{aligned}
  \end{equation}
whose image is a truncated, translated, and scaled catenoid with vertical axis of symmetry.
For $1 \leq i \leq N-1$ we set
  \begin{equation}
  \label{adef}
    \begin{aligned}
      a_i
      &=
      a_i[N,k,\ell,m,\zeta,\xi]
      :=
      \arcosh \frac{1}{10 \ell m \tauj{i}} \\
      &=
      \frac{k \ell}{4\pi}\left(1-\frac{\cj{2}}{2}\right)m^2 - \zetaj{1}
        -\left(1-\delta_{i1}\right)\left(\ln \cj{i}+\frac{\zetaj{i}}{k \ell m^2}\right)
        +\ln \left(1+\sqrt{1-100\ell^2m^2\tauj{i}^2}\right)
    \end{aligned}
  \end{equation}
(writing $\delta_{ij}$ for the Kronecker delta
and recalling \eqref{cisummary}, \eqref{taudef}, \eqref{almosta}, and \eqref{lnarcosh})
and, recalling \eqref{Phidef}, we define the map
  \begin{equation}
  \label{catdef}
    \begin{aligned}
      &\cat_i=\cat_i[N,k,\ell,m,\zeta,\xi]: \cyl_{a_i} \to \Sph^3 \mbox{ by} \\
      &\cat_i
      :=
      \Phi \circ \cate
      \left[
        \left(
          \frac{(i-1)\pi}{\sqrt{2}km}, \frac{(i-1)\pi}{\sqrt{2} \ell m}, \zzj{i}^K
      \right), \tauj{i}, a_i
    \right].
    \end{aligned}
  \end{equation}
Then
(referring to \eqref{initsurfdef} and particularly \eqref{Omega})
the image of $\cat_i$ is entirely contained in the initial surface $\Sigma$
and defines the catenoidal region
  \begin{equation}
  \label{catrdef}
    \catr[i]=\catr[i;N,k,\ell,m,\zeta,\xi]:= \cat_i \left(\cyl_{a_i}\right),
  \end{equation}
so that $\cat_i$ is a diffeomorphism onto its image;
in innocuous abuse of notation
we will routinely write $\cat_i^{-1}: \catr[i] \to \cyl_{a_i}$
for the inverse of this diffeomorphism.

Note that, recalling \eqref{rhodef},
  \begin{equation}
  \label{rhocat}
    {\cat_i}^*\rho=\tauj{i}^{-1} \sech t
  \end{equation}
and, by applying \eqref{cisummary} and \eqref{taudef} to \eqref{adef},
  \begin{equation}
  \label{aest}
    C[N,k,\ell]^{-1}m^2 \leq a_i \leq C[N,k,\ell]m^2
  \end{equation}
for some constant $C[N,k,\ell]>0$
whenever $\zeta,\xi \in [-c,c]^{N-1}$
and $m$ is sufficiently large in terms of $N$, $k$, $\ell$, and $c$.
From \eqref{taudef}, \eqref{catheights}, and \eqref{catrdef}
it is clear that
  \begin{equation}
    \catr[i] \cap \catr[i'] \neq \emptyset
      \mbox{ if and only if $i=i'$}.
  \end{equation}

Since on small scales
the covering map $\Phi$ is an approximate isometry,
we expect that in a suitably rescaled sense
each catenoidal region will
converge in the large-$m$ limit to an exact catenoid in Euclidean space.
The next proposition quantifies this convergence.

\begin{prop}[Estimates of the geometry of the catenoidal regions]
\label{catest}
Given a real number $c>0$ and integers $N \geq 2$, $\ell \geq k \geq 1$, and $j \geq 0$,
there exist real numbers $m_0=m_0[N,k,\ell,c]>0$ and $C=C[N,k,\ell,j]>0$
such that whenever
$\zeta,\xi \in [-c,c]^{N-1}$ and $m>m_0$,
for $1 \leq i \leq N-1$
  \begin{enumerate}[(i)]
    \item $\norm{{\cat_i}^*\chi - \chihatK:
      C^j\left(T^*\cyl_{a_i}^{\otimes 2},\chihatK \right)} \leq Cm^2\tauj{1}$,

    \item $\norm{\rho:C^j(\catr[i],\chi,\rho)}+\norm{\rho^{-1}:C^j(\catr[i],\chi,\rho^{-1})}
           \leq C$

    \item $\norm{\zz: C^j(\catr[i],\chi)} \leq Cm^2\tauj{1}$,

    \item $\norm{A-(-1)^{N-i}{\cat_i^{-1}}^*\tauj{i}\left(dt^2-d\theta^2\right):
           C^j\left(T^*\catr[i]^{\otimes 2},\chi, \tauj{1}\abs{\zz}+\rho^{-2} \right)}
           \leq
           C$,

    \item $\norm{\rho^{-2}\abs{A}_g^2-2\tau^2 \rho^2:
           C^j(\catr[i],\chi,\tauj{1}+\abs{\zz}+\rho^{-2})}
           \leq C$, and

    \item $\norm{\rho^{-2}H : C^j\left(\catr[i],\chi,
      \tauj{1}\abs{\zz}+\rho^{-2}\abs{\zz}+\tauj{1}^2\right)} \leq C$,
  \end{enumerate}
where $\zz: \Sph^3 \to \R$ is defined via \eqref{Phidef}
and we also recall \eqref{taudef}, \eqref{weightedhoelderdef}, \eqref{gdef}, \eqref{AHdef},
\eqref{chidef}, \eqref{rhodef}, \eqref{chihatK}, \eqref{cyladef},
\eqref{catdef}, and \eqref{catrdef}.
\end{prop}

\begin{proof}
From \eqref{Phullback}, \eqref{catedef}, \eqref{catdef}, and \eqref{rhocat}
we calculate
  \begin{equation}
  \label{chidiff}
    {\cat_i}^*\chi - \chihatK
    =
    \left(\sin 2\left(\zzj{i}^K+\tauj{i} t\right)\right)
      \left(
        \tanh^2 t \, \cos 2 \theta \, dt^2
        - 2 \tanh t \, \sin 2\theta \, dt \, d\theta
        - \cos 2\theta \, d\theta^2
      \right).
  \end{equation}
It follows from \eqref{taudef},
\eqref{estimatesforforceestimate} (the bottom line),
and \eqref{aest}
that $\abs{\zzj{i}^K+\tauj{i}t} \leq C[N,k,\ell]m^2\tauj{1}$
everywhere on $\cyl_{a_i}$
for some $C[N,k,\ell]>0$ 
whenever $m$ is sufficiently large in terms of $N$, $k$, $\ell$, and $c$.
Since $\chihatK=dt^2+d\theta^2$ is flat,
it is trivial to differentiate \eqref{chidiff},
so item (i) now follows immediately,
as do items (ii) and (iii) in turn, using also \eqref{rhocat} and
  \begin{equation}
    {\cat_i}^*\zz=\zzj{i}^K+\tauj{i}t.
  \end{equation}

Next, using \eqref{catdef}, for any vector fields $V$ and $W$ on $\cyl_{a_i}$ we have
  \begin{equation}
  \label{Acomp}
    \begin{aligned}
      {\cat_i}^*A(V,W)
      &=
      A[\cate_i,\Phi^*\gsph](V,W)
      =
      \nu[\cate_i,\Phi^*\gsph]_c\left(D[\Phi^*\gsph]_{d\cate_iV}d\cate_iW\right)^c \\
      &=
      \frac{\nu[\cate_i,\geuc]_c}{\abs{\nu[\cate_i,\geuc]}_{\Phi^*\gsph \circ \cate_i}}
        \left(D[\geuc]_{d\cate_i V} d\cate_i W
          +\left(D[\Phi^*\gsph]-D[\geuc]\right)(d\cate_i V,d\cate_i W)\right)^c \\
      &=
      \abs{\nu[\cate_i,\geuc]}_{\Phi^*\gsph \circ \cate_i}^{-1}
        \left(A[\cate_i,\geuc](V,W)+B(V,W)\right),
    \end{aligned}
  \end{equation}
where
$\cate_i:=\cate\left[
  \left(\frac{(i-1)\pi}{\sqrt{2}km},\frac{(i-1)\pi}{\sqrt{2} \ell m}, \zzj{i}^K\right),
  \tauj{i}, a_i \right]$
(recalling \eqref{catedef}),
$\nu[\cate_i,\Phi^*\gsph]$ is the unit normal for $\cate_i$ relative to $\gsph$
directed so that $d\Phi \nu[\cate_i,\Phi^*\gsph]=\nu \circ \cat_i$
and $\nu[\cate_i,\Phi^*\gsph]_c$ is its $\Phi^*\gsph$ metric dual,
$\nu[\cate_i,\geuc]$ is the unit normal for $\cate_i$ relative to $\geuc$
directed so that its $\geuc$ metric dual $\nu[\cate_i,\geuc]_c$
is a positive multiple of $\nu[\cate_i,\Phi^*\gsph]_c$,
$A[\cate_i,\cdot]$ is the second fundamental of $\cate_i$
with respect to
the ambient metric $\cdot$
and the unit normal $\nu[\cate_i,\cdot]$,
$D[\cdot]$ is the Levi-Civita connection induced by the metric $\cdot$,
and $B$ is the symmetric tensor
  \begin{equation}
  \label{Bdef}
    \begin{aligned}
      &B_{\alpha \beta}
      :=
      \nu[\cate_i,\geuc]_c (\Gamma^c_{\;\; ab} \circ \cate_i)
       (d\cate_i)^a_{\;\; \alpha}(d\cate_i)^b_{\;\; \beta}
        \mbox{ with} \\
      &\Gamma^c_{\;\;ab}
      :=
      \frac{1}{2}(\Phi^*\gsph)^{cd}
      \left(
        D[\geuc]_b(\Phi^*\gsph)_{ad}+D[\geuc]_a(\Phi^*\gsph)_{bd}-D[\geuc]_d(\Phi^*\gsph)_{ab}
      \right).
    \end{aligned}
  \end{equation}

Recalling that we have chosen the unit normal $\nu$ for $\Sigma$ pointing toward $C_1$
at the points closest to $C_1$
(so $\nu$ has positive inner product with $\partial_{\zz}$
at the top of $\Sigma$,
in an $(\xx,\yy,\zz)$ coordinate system defined via $\Phi$),
we see from \eqref{catedef} that
  \begin{equation}
    \nu[\cate_i,\geuc]=(-1)^{N-i}\widehat{\nu},
  \end{equation}
recalling \eqref{catenormal},
so in particular, using \eqref{Phullback},
  \begin{equation}
  \label{normalization}
    \abs{\nu[\cate_i,\geuc]}_{\Phi^*\gsph \circ \cate_i}(t,\theta)
    =
    \frac{\sqrt{1-\sech^2 t \, \cos 2\theta \, \sin 2{\cat_i}^*\zz
       -\tanh^2 t \,  \sin^2 2{\cat_i}^*\zz}}{\cos 2{\cat_i}^*\zz}.
  \end{equation}
It is also easy to see from \eqref{Phullback}
that the only Christoffel symbols (in $(\xx,\yy,\zz)$ coordinates)
for $\Phi^*\gsph$ not vanishing identically are
  \begin{equation}
    \Gamma^{\zz}_{\;\; \xx \xx}=-\Gamma^{\zz}_{\;\; \yy \yy}=-\cos 2\zz,
    \qquad
    \Gamma^{\xx}_{\;\; \xx \zz}=\Gamma^{\xx}_{\;\; \zz \xx}=\frac{\cos 2\zz}{1+\sin 2\zz},
    \qquad
    \Gamma^{\yy}_{\;\; \yy \zz}=\Gamma^{\yy}_{\;\; \zz \yy}=\frac{-\cos 2\zz}{1-\sin 2\zz}.
  \end{equation}
Returning to \eqref{Bdef} (and again using \eqref{catenormal} and \eqref{catedef}) we now find
  \begin{equation}
  \label{Bcomp}
    \begin{aligned}
      (-1)^{N-i}B_{tt}
      &=
      {\cat_i}^*\rho^{-2} \tanh^3 t \, \cos 2\theta \, \cos 2{\cat_i}^*\zz
        +2\tauj{i}^2 \tanh t \, \cos 2\theta \, \sec \, 2{\cat_i}^*\zz 
        -2\tauj{i}^2 \tanh t \, \tan 2{\cat_i}^*\zz, \\
      (-1)^{N-i}B_{t \theta}
      &=
      -{\cat_i}^*\rho^{-2} \tanh^2 t \, \sin 2\theta \, \cos 2{\cat_i}^*\zz
        -\tauj{i}^2 \sin 2\theta \, \sec 2{\cat_i}^*\zz, \mbox{ and} \\
      (-1)^{N-i}B_{\theta \theta}
      &=
      -{\cat_i}^*\rho^{-2}\tanh t \, \cos 2\theta \, \cos 2{\cat_i}^*\zz.
    \end{aligned}
  \end{equation}

Thus, applying \eqref{normalization} and \eqref{Bcomp} in \eqref{Acomp},
we have computed $A$ on $\catr[i]$,
proving (iv).
From \eqref{chidiff} we also compute
  \begin{equation}
  \label{catchiinverse}
    \left(\cate_i^*\Phi^*\gsph\right)^{-1}
    =
    {\cat_i}^*\rho^2
    \frac{
      \partial_t^2+\partial_\theta^2
      +(\sin 2{\cat_i}^*\zz)
      \left(
        -\cos 2\theta \, \partial_t^2
        +2\tanh t \, \sin 2\theta \, \partial_t \, \partial_\theta
        +\tanh^2 t \, \cos 2\theta \, \partial_\theta^2
      \right)
      }
      {1-\sech^2 t \, \cos 2\theta \, \sin 2{\cat_i}^*\zz-\tanh^2 t \,  \sin^2 2{\cat_i}^*\zz},
  \end{equation}
which in conjunction with (iv) immediately yields the estimate (v).
The estimate (vi) requires slightly more care,
but from \eqref{cateA}, \eqref{Acomp}, \eqref{normalization}, and \eqref{catchiinverse}
we compute the exact mean curvature
  \begin{equation}
    \begin{aligned}
      {\cat_i}^*H
      =
      &(-1)^{N-i}
      \left(1-\sech^2 t\,\cos 2\theta\,\sin 2{\cat_i}^*\zz-\tanh^2 t\,\sin^2 2{\cat_i}^*\zz\right)^{-3/2} \left({\cat_i}^*\rho^2\right) \\
      &\times
      \left[
        -\frac{1}{2}\tauj{i}(\sin 4 {\cat_i}^*\zz)(1+\tanh^2 t) \cos 2\theta
        -{\cat_i}^*\rho^{-2}(\sin 4{\cat_i}^*\zz)(\cos 2{\cat_i}^*\zz)\tanh^3 t \right. \\
      &\;\;\;\;\;\; \left. +\tauj{i}^2 \tanh t \, \cos 2\theta
        -4\tauj{i}^2(\sin 2{\cat_i}^*\zz)\tanh t
        +3\tauj{i}^2(\sin^2 2{\cat_i}^*\zz) \tanh t \, \cos 2\theta \vphantom{\frac{1}{2}}
      \right],
    \end{aligned}
  \end{equation}
delivering (vi) and completing the proof.
\end{proof}

\subsection*{Graphs over immersions}
The estimates away from the catenoidal regions will be obtained
by treating the initial surface there as a graph over the torus,
as an application of the following lemma,
which will
be used again in the final section
to estimate the contributions
to the mean curvature of the perturbed surface which are nonlinear
in the perturbing function
and to estimate the perturbation to the corresponding forces \eqref{forcedef}.
We first clarify some notation used in the statement of the lemma.
Suppose $(M,g)$ is a Riemannian manifold.
If $u \in C^2_{loc}(M)$,
we write ${D[g]^2}_{ab}u$ and $\Delta_g u:=g^{ab}{D[g]^2}_{ab}u$
for the Hessian and Laplacian of $u$ under $g$.
Given $p \in M$ and $r>0$,
we write $\overline{B}[p,r,(M,g)]$ for the closed metric ball in $(M,g)$
with center $p$ and radius $r$.
We adopt the sign and indexing conventions that
the Riemann curvature tensor $R_{abcd}$ of $(M,g)$
is defined by
  \begin{equation}
  \label{curvatureconvention}
    R_{abcd}V^aW^bX^cY^d
    =
    g\left(D[g]_V D[g]_W X - D[g]_W D[g]_V X - D[g]_{[V,W]} X, Y\right),
  \end{equation}
for any vector fields $V,W,X,Y$ on $M$;
then $R_{ab}:=R_{ac\;\;\;b}^{\;\;\;\; c}={R_{cab}}^c$
is the Ricci curvature of $(M,g)$.
Suppose also that $\phi: S \to M$ is a $C^2_{loc}$ codimension-one
immersion of a manifold $S$ into $M$
and that $\nu$ is a global unit normal for $\phi$.
We write $A[\phi,\nu]:=(g \circ \phi)(\nu,D[g]d\phi)$
and $H[\phi,\nu]:=\tr_{\phi^*g}A[\phi,\nu]$
for the corresponding second fundamental form and mean curvature
(here $D[g]$ being the canonical connection
on the bundle $\phi^*TM \otimes T^*S$ defined by the Levi-Civita
connections induced by $g$ and $\phi^*g$).
Finally we point out that we reserve the right
to denote evaluation of a section $X$ at a point $p$
by either of the standard options $X|_p=X(p)$.

\begin{lemma}[Graphs over immersions]
\label{perimm}
Let $\phi: S \to M$ be a smooth two-sided (codimension-one) immersion
of a smooth manifold $S$ into a smooth complete Riemannian manifold $M$
with smooth metric $g$.
Let $\nu \in \phi^*(TM)$ be a global unit normal for $\phi$
and write $A:=A[\phi,\nu]$ and $H:=H[\phi,\nu]$
for the corresponding second fundamental form
and mean curvature.
For each $t \in \R$ and $u \in C^2_{loc}(S)$ we define the maps
$\phi_t, \phi[u]: S \to M$ by
  \begin{equation}
    \phi_t(p):=\exp_{\phi(p)} t\nu(p)
    \qquad \mbox{and} \qquad
    \phi[u](p):=\exp_{\phi(p)} u(p)\nu(p)=\phi_{u(p)}(p),
  \end{equation}
where $\exp: TM \to M$ is the exponential map of $(M,g)$.
Suppose that for a given $u \in C^2_{loc}(S)$ and $p \in S$
  \begin{equation}
  \label{immineq}
    \abs{u(p)}\abs{A(p)}_{\phi^*g}
      +\abs{u(p)}^2\norm{\abs{R}_g: C^0(\overline{B}[p,\abs{u(p)},(M,g)])}
    <
    \frac{1}{3}.
  \end{equation}
Then $\phi[u]$ is a $C^2_{loc}$ immersion
on a neighborhood of $p$,
as is $\phi_t$ for every $t$ between $0$ and $u(p)$.
On this neighborhood $\phi_t$ and $\phi[u]$
admit respective $C^1_{loc}$ unit normals $\nu_t$ and $\nu[\phi[u]]$
satisfying $\nu_t(p)=\frac{d}{dt} \exp_{\phi(p)} t\nu(p)$
and $g\left(\nu[\phi[u]]|_p, \nu_t|_p\right)>0$.
If (working near $p$)
we set $g^t:=\phi_t^*g$,
$g_t:=\left(\phi_t^*g\right)^{-1}$,
and $A^t:=A[\phi_t,\nu_t]$,
then
  \begin{equation}
  \tag{i}
    \begin{aligned}
      &\partial_t {g^t}_{\alpha \beta}
      =
      -2{A^t}_{\alpha \beta}
      \qquad \mbox{and} \qquad
      \partial_t {A^t}_{\alpha \beta}
      =
      {\nu_t}^a{\nu_t}^d{{\phi_t}^b}_{,\alpha}{{\phi_t}^c}_{,\beta}R_{abcd} \circ \phi_t
        -{A^t}_{\alpha \gamma}{A^t}_{\beta \delta}{g_t}^{\gamma \delta}.
    \end{aligned}
  \end{equation}
If we also set
$A[u]:=A[\phi[u],\nu[\phi[u]]]$ and $H[u]:=H[\phi[u],\nu[\phi[u]]]$
and we define on $S$ the symmetric $2$-tensors
$g^u:=g^{u(\cdot)}$, $g_u:=g_{u(\cdot)}$, and $A^u:=A^{u(\cdot)}$,
as well as the function $H^u:={g_u}^{\alpha \beta}{A^u}_{\alpha \beta}$
and the section
${{\phi_u}^c}_{,\alpha}:=\left.{{\phi_{t}}^c}_{,\alpha}\right|_{t=u(\cdot)}$
of $\phi[u]^*TM$,
then
  \begin{equation*}
  \tag{ii}
    \phi[u]^*g = g^u +  du \otimes du,
  \end{equation*}

  \begin{equation*}
  \tag{iii}
    A[u]_{\alpha \beta}
    = \frac{
      {A^u}_{\alpha \beta}
      +{D[g^u]^2}_{\alpha \beta}u
      -2{g_u}^{\gamma \delta}{A^u}_{\gamma (\alpha}u_{,\beta)}u_{,\delta}
    }
    {\sqrt{1+\abs{du}_{g^u}^2}}, \mbox{ and}
  \end{equation*}

  \begin{equation*}
  \tag{iv}
    \begin{aligned}
      H[u]
      =
      &H + \left(\Delta_{\phi^*g} + \abs{A}_g^2 + \nu^a\nu^bR_{ab} \circ \phi\right)u \\
      &+u_{,\gamma}{g_u}^{\gamma \delta}{g_u}^{\alpha \beta}
          \int_0^1 \left(2u{A^{tu}}_{\alpha \delta;\beta}
            +2u_{,\beta}{A^{tu}}_{\alpha \delta}
            -u{A^{tu}}_{\alpha \beta; \delta}
            -u_{,\delta}{A^{tu}}_{\alpha \beta}\right) \, dt \\
      &-\frac{
           u_{,\gamma}u_{,\delta}
           {g_u}^{\alpha \gamma}
           {g_u}^{\beta \delta}
        }
        {\sqrt{1+\abs{du}_{g^u}^2}}
        \left(
          \frac{H^u+\Delta_{g^u}u}
            {1+\sqrt{1+\abs{du}_{g^u}^2}}
            {g^u}_{\alpha \beta}
          +\frac{{D[g^u]^2}_{\alpha \beta}u + 3{A^u}_{\alpha \beta}}
            {1+\abs{du}_{g^u}^2}
        \right) \\
      &+u^2\int_0^1 (t-1)
        \left[
        \vphantom{\frac{1}{2}}
          2{g_{tu}}^{\alpha \beta}{g_{tu}}^{\gamma \delta}{A^{tu}}_{\alpha \beta}
            {\nu_{tu}}^a{\nu_{tu}}^d{{\phi_{tu}}^b}_{,\gamma} {{\phi_{tu}}^c}_{,\delta}
            R_{abcd} \circ \phi[u]
          \right. \\
      &\left.
          +2{A^{tu}}_{\alpha \beta}{A^{tu}}_{\gamma \delta}{A^{tu}}_{\epsilon \zeta}
            {g_{tu}}^{\beta \gamma}{g_{tu}}^{\delta \epsilon}{g_{tu}}^{\alpha \zeta}
          +{\nu_{tu}}^a{\nu_{tu}}^b{\nu_{tu}}^c R_{ab|c} \circ \phi[u]
        \vphantom{\frac{1}{2}}
        \right] \, dt,
    \end{aligned}
  \end{equation*}
where the vertical bar $|$ and semilcolon $;$ before an index
indicate differentiation under the Levi-Civita connection
induced by $g$ and $\phi^*g$ respectively.
\end{lemma}

\begin{remark}
\label{spherical}
Note that in the special case of Lemma \ref{perimm}
that $(M,g)=(\Sph^3,\lambda^2 \gsph)$ is the round $3$-sphere of radius $\lambda>0$
we have
  \begin{equation}
    {\nu_t}^a{\nu_t}^d{{\phi_t}^b}_{,\alpha}{{\phi_t}^c}_{,\beta}R_{abcd} \circ \phi_t
    =\lambda^{-2}{g^t}_{\alpha \beta},
    \qquad
    \nu^a \nu^b R_{ab} \circ \phi = 2\lambda^{-2},
    \qquad \mbox{and} \qquad
    R_{ab|c}=0.
  \end{equation}
\end{remark}

\begin{proof}
We begin with a few basic generalities
concerning connections on pullbacks of vector bundles.
Suppose $\varphi: P \to M$ is a smooth map (not necessarily an immersion)
between smooth manifolds
and the target $M$ is equipped with a smooth Riemannian metric $g$.
We will write $D[TM]$ for the Levi-Civita connection on $TM$ induced by $g$.
We omit the elementary verification of the following observations.
There is a unique connection $D[\varphi^*TM]$
on the pullback bundle $\varphi^*TM$
satisfying the chain rule
  \begin{equation}
    D[\varphi^*TM]_V (X \circ \varphi)=\left(D[TM]_{\varphi_*V}X\right) \circ \varphi
  \end{equation}
for all $V \in C^0_{loc}(TP)$ and $X \in C^1_{loc}(TM)$.
Moreover, $D[\varphi^*TM]$ is torsion-free
in the sense that
  \begin{equation}
  \label{torsionfree}
    D[\varphi^*TM]_V \varphi_*W - D[\varphi^*TM]_W \varphi_*V = \varphi_*[V,W]
  \end{equation}
for all $V,W \in C^1_{loc}(TP)$;
$D[\varphi^*TM]$ is compatible with $g$ in the sense that
  \begin{equation}
  \label{metriccompatible}
    V(g \circ \varphi)(X,Y)=(g \circ \varphi)(D[\varphi^*TM]_V X, Y)+(g \circ \phi)(X,D[\varphi^*TM]_V Y)
  \end{equation}
for all $V \in C^0(TP)_{loc}$ and $X,Y \in C^1_{loc}(\varphi^*TM)$;
and $D[\varphi^*TM]$ inherits the curvature of $M$:
for all $V,W \in C^1_{loc}(TP)$ and $X \in C^2_{loc}(\varphi^*TM)$
  \begin{equation}
  \label{inherittheR}
    D[\varphi^*TM]_V D[\varphi^*TM]_W X - D[\varphi^*TM]_W D[\varphi^*TM]_V X - D[\varphi^*TM]_{[V,W]}X
    =
    (R \circ \varphi)(\varphi_*V,\varphi_*W)X.
  \end{equation}

Now let $(M,g)$, $S$, $\phi$, and $\nu$
be as in the statement of the lemma.
We define the map
  \begin{equation}
    \Phi: S \times \R \to M
    \quad \mbox{by} \quad
    \Phi(p,t):=\exp_{\phi(p)} t\nu(p),
  \end{equation}
so that $\phi_t=\Phi(\cdot,t)$.
Suppose $V\in C^\infty(TS)$ and write
$\mathbb{V}$ for the unique vector field on $\Sigma \times \R$
such that
$(\mathbb{V}f)(p,t)=(Vf(\cdot,t))(p)$
for all $f \in C^\infty(S \times \R)$, $p \in S$, and $t \in \R$.
Then $\Phi_*\mathbb{V}|_{(p,t)}=d\phi_t V|_p$
and $\Phi_*\partial_t|_{(p,t)}=\frac{d}{dt} \exp_{\phi(p)} t\nu(p)$.
Given $s,t \in \R$, we write $P_s^t: \Phi^*TM|_{(\cdot,s)} \to \Phi^*TM|_{(\cdot,t)}$
for the map of parallel translation (relative to $D[\Phi^*TM]$)
along the $\R$ cross-sections of $S \times \R$.
Using \eqref{torsionfree} and \eqref{inherittheR} 
as well as the fact that $D[\Phi^*TM]_{\partial_t}\Phi_*\partial_t=0$,
we have
  \begin{equation}
    \begin{aligned}
      \left.\frac{d}{dt}\right|_{t=0}P_t^0 \left(d\phi_t V\right)
      &=
      \left.\frac{d}{dt}\right|_{t=0}P_t^0 \left(\Phi_* \mathbb{V}\right)
      =
      D[\Phi^*TM]_{\partial_t}\Phi_* \mathbb{V}|_{(\cdot,0)}
      =
      D[\Phi^*TM]_{\mathbb{V}} \Phi_*\partial_t|_{(\cdot,0)} \\
      &=
      D[\phi^*TM]_V \nu \mbox{ and} \\
      \frac{d^2}{dt^2}P_t^0 \left(d\phi_t V\right)
      &=
      \frac{d^2}{dt^2}P_t^0 \left(\Phi_*\mathbb{V}\right)
      =
      \frac{d}{dt}P_t^0 \left(D[\Phi^*TM]_{\partial_t} \Phi_*\mathbb{V}\right)
      =
      \frac{d}{dt}P_t^0 \left(D[\Phi^*TM] D_{\mathbb{V}} \Phi_*\partial_t\right) \\
      &=
      P_t^0 \left(D[\Phi^*TM]_{\partial_t} D_{\mathbb{V}} \Phi_*\partial_t\right)
      =
      P_t^0 \left((R \circ \phi_t)(\Phi_*\partial_t, d\phi_t V)\Phi_*\partial_t \right).
    \end{aligned}
  \end{equation}
Thus
  \begin{equation}
    P_t^0 d\phi_t V
    =
    d\phi V + tD[\phi^*TM]_V \nu
      + \int_0^t (t-s)
      P_s^0 \left((R \circ \phi_s)(\Phi_*\partial_t, d\phi_s V)\Phi_*\partial_t \right) \, ds,
  \end{equation}
so, noting that $\Phi_*\partial_t$ is unit,
for all $p \in S$ and $t \in \R$ (replacing $[0,t]$ below by $[t,0]$ if $t<0$)
  \begin{equation}
    \sup_{s \in [0,t]} \abs{d\phi_sV_p}_g
    \leq
    \frac{1+\abs{t}\abs{A(p)}_{\phi^*g}}{1-t^2\norm{\abs{R}_g : C^0(\overline{B}[p,\abs{t},(M,g)])}}
      \abs{d\phi V_p}_g
  \end{equation}
and also consequently
  \begin{equation}
  \label{phitmetlowerbound}
    \frac{\abs{d\phi_t V_p}_g}{\abs{d\phi V_p}_g}
      \geq
      1-\abs{t}\abs{A(p)}_{\phi^*g}-t^2\norm{\abs{R}_g : C^0(\overline{B}[p,\abs{t},(M,g)])}
        \frac{1+\abs{t}\abs{A(p)}_{\phi^*g}}{1-t^2\norm{\abs{R}_g : C^0(\overline{B}[p,\abs{t},(M,g)])}},
  \end{equation}
which confirms $\phi_t$ is an immersion near $p$
whenever
$\abs{t}\abs{A(p)}_{\phi^*g}+t^2\norm{\abs{R}_g : C^0(\overline{B}[p,\abs{t},(M,g)])}<1/3$
(which condition is obviously not sharp).

Using \eqref{torsionfree} and \eqref{metriccompatible}
as well as the fact that $D[\Phi^*TM]_{\partial_t}\Phi_*\partial_t=0$,
we also compute
  \begin{equation}
    \begin{aligned}
      &(g \circ \Phi)(D[\Phi^*TM]_{\mathbb{V}} \Phi_*\partial_t,\partial_t)
      =
      \frac{1}{2}\mathbb{V}(g \circ \Phi)(\Phi_*\partial_t,\Phi_*\partial_t)
      =
      0 \mbox{ and} \\
      &\frac{d}{dt}(g \circ \Phi)(\Phi_*\mathbb{V},\Phi_*\partial_t)
      =
      (g \circ \Phi)(D[\Phi^*TM]_{\partial_t}\Phi_*\mathbb{V},\Phi_*\partial_t)
      =
      \mathbb{V}(g \circ \Phi)(\Phi_*\partial_t,\Phi_*\partial_t)
      =
      0.
    \end{aligned}
  \end{equation}
Since $\Phi_*\partial_t|_{(\cdot,0)}=\nu$,
in fact $\Phi_*\partial_t$ is everywhere and always orthogonal
to $\Phi_*\mathbb{V}$.
Thus wherever $\abs{t}$ is small enough that $\phi_t$ is locally an immersion,
$\Phi_*\partial_t$ is a smooth local unit normal for $\phi_t$,
designated $\nu_t$ in the statement of the lemma.
Then, letting $W$ be another vector field on $S$
with $\mathbb{W}$ the canonically corresponding vector field on $S \times \R$,
\eqref{torsionfree}--\eqref{inherittheR},
  \begin{equation}
    \begin{aligned}
      &\frac{d}{dt}\phi_t^*g(V,W)
      =
      (g \circ \Phi)(D[\Phi^*TM]_{\mathbb{V}} \Phi_*\partial_t,\Phi_*\mathbb{W})
        +(g \circ \Phi)(\Phi_*\mathbb{V},D[\Phi^*TM]_{\mathbb{W}} \Phi_*\partial_t) \mbox{ and} \\
      &\frac{d}{dt}A[\phi_t,\nu_t](V,W)
      =
      (g \circ \Phi)
        \left(
          (R \circ \Phi)(\Phi_*\partial_t,\Phi_*\mathbb{V})\mathbb{W}
            +D[\Phi^*TM]_{\mathbb{V}}D[\Phi_*TM]_{\mathbb{W}}\Phi_*\partial_t,
          \Phi_*\partial_t
        \right),
    \end{aligned}
  \end{equation}
proving item (i) of the lemma.

Now suppose $u \in C^2_{loc}(S)$
and take $\phi[u]$ as defined in the statement of the lemma.
Let $\pi: S \times \R \to S$ be the canonical projection onto $S$
and let $\nu_u$ be the section of $\phi[u]^*TM$
defined by
  \begin{equation}
    \nu_u:=\nu_{u(\cdot)}
  \end{equation}
(so $\nu_u(p)=\nu_{u(p)}(p)$, $\nu_t$ having been defined in the statement of the lemma).
Since $\phi[u](p)=\Phi(p,u(p))$,
  \begin{equation}
    d\phi[u]V
    =
    d\Phi \left(\mathbb{V}+(\mathbb{V}\pi^*u)\partial_t\right)
    =
    \left.\left(d\phi_{t}\right)\right|_{t=u(\cdot)}V + (Vu)\nu_u,
  \end{equation}
which implies item (ii) of the lemma.
In particular, because $du \otimes du$ is nonnegative,
$\phi[u]$ is an immersion on a neighborhood of $p$
whenever $\phi_{u(p)}$ is,
so in particular, in view of \eqref{phitmetlowerbound},
provided \eqref{immineq} holds.
It also follows that the corresponding metric on the cotangent space satisfies
  \begin{equation}
  \label{phiudualmetric}
    (\phi[u]^*g)^{-1}
    =
    g_u
      -\frac{1}{1+\abs{du}_{g^u}^2}
        \nabla_{g^u}u \otimes \nabla_{g^u}u,
  \end{equation}
where $g_u$ and $g^u$ are as defined in the statement of the lemma
and $\nabla_{g^u}$ is the gradient operator on $S$
induced by the metric $g^u$.
(Equation \eqref{phiudualmetric} is a trivial consequence of item (ii) of the lemma
at points where $du$ vanishes;
at any other point $p$ 
it is easily derived by working relative to a $g^u$ orthogonal basis
one of whose elements is
$\nabla_{g^u}u|_p$.)

Clearly the $1$-form $dt-d\pi^*u$ on $S \times \R$
annihilates all tangent vectors to the graph of $u$ in $S \times \R$,
so, relative to the metric $\Phi^*g=\phi_t^*g+dt^2$,
the upward unit normal to this graph is
  \begin{equation}
    \frac{\partial_t - \nabla_{\Phi^*g} \pi^*u}{\sqrt{1+\abs{d\pi^*u}_{\Phi^*g}^2}},
  \end{equation}
Noting that
$d\Phi \left(\nabla_{\Phi^*g} \pi^*u\right)|_{(p,t)} = d\phi_t \left(\nabla_{\phi_t^*g}u\right)|_p$,
we see that the unit normal $\nu[\phi[u]]$ for $\phi[u]$ identified
in the statement of the lemma satisfies
  \begin{equation}
  \label{phiunormal}
    \nu[\phi[u]]
    =
    \frac{\nu_u - \left.\left(d\phi_t\right)\right|_{t=u(\cdot)}\nabla_{g^u}u}
      {\sqrt{1+\abs{du}_{g^u}^2}}.
  \end{equation}
The corresponding second fundamental form is
$A[u]=(g \circ \phi[u])(\nu[\phi[u]],D[\phi[u]^*TM]_V d\phi[u]W)$,
but
  \begin{equation}
    \begin{aligned}
      D[\phi[u]^*TM]_V d\phi[u]W |_p
      &=
      D[\Phi^*TM]_{\mathbb{V}+(\mathbb{V}\pi^*u)\partial_t}
        \left(\Phi_*\mathbb{W}+(\mathbb{W}\pi^*u)\Phi_*\partial_t\right)|_{(p,u(p)} \\
      &=
      D[\phi_{u(p)}^*TM]_V \left.\left(\left.d\phi_t\right|_p\right)\right|_{t=u(p)}W|_p
        +(VWu)\nu[\phi_{u(p)}]|_p \\
        &\;\;\;\; +(Wu)D[\phi_{u(p)}^*TM]_V \nu[\phi_{u(p)}]|_p
        +(Vu)D[\phi_{u(p)}^*TM]_W \nu[\phi_{u(p)}]|_p,
    \end{aligned}
  \end{equation}
whose inner product with \eqref{phiunormal} yields item (iii) of the lemma.
By contracting item (ii) of the lemma with item (iii)
we obtain
  \begin{equation}
    \begin{aligned}
      H[u]
      &=
      \frac{H^u+\Delta_{g^u}u}{\sqrt{1+\abs{du}_{g^u}^2}}
        -u_{,\gamma} u_{,\delta}
         \frac{{D[g^u]^2}_{\alpha \beta}u+3{A^u}_{\alpha\beta}}
         {\left(1+\abs{du}_{g^u}^2\right)^{3/2}}
         {g_u}^{\alpha \gamma} {g_u}^{\beta \delta} \\
      &=
      H^u+\Delta_{g^u}u
        -\frac{
           u_{,\gamma}u_{,\delta}
           {g_u}^{\alpha \gamma}
           {g_u}^{\beta \delta}
        }
        {\sqrt{1+\abs{du}_{g^u}^2}}
        \left(
          \frac{H^u+\Delta_{g^u}u}
            {1+\sqrt{1+\abs{du}_{g^u}^2}}
            {g^u}_{\alpha \beta}
          +\frac{{D[g^u]^2}_{\alpha \beta}u + 3{A^u}_{\alpha \beta}}
            {1+\abs{du}_{g^u}^2}
        \right),
    \end{aligned}
  \end{equation}
but, using item (i) of the lemma,
  \begin{equation}
    \begin{aligned}
      \Delta_{g^u}u
      &=
      \Delta_{\phi^*g}u
        +u_{,\gamma}{g_u}^{\gamma \delta}{g_u}^{\alpha \beta}
          \left(\frac{1}{2}D[\phi^*g]_\delta{g^u}_{\alpha \beta}
          -D[\phi^*g]_{\beta}{g^u}_{\alpha \delta}\right) \\
      &=
      \Delta_{\phi^*g}u
        +u_{,\gamma}{g_u}^{\gamma \delta}{g_u}^{\alpha \beta}
          \int_0^1 \left(2u{A^{tu}}_{\alpha \delta;\beta}
            +2u_{,\beta}{A^{tu}}_{\alpha \delta}
            -u{A^{tu}}_{\alpha \beta; \delta}
            -u_{,\delta}{A^{tu}}_{\alpha \beta}\right) \, dt
    \end{aligned}
  \end{equation}
and
  \begin{equation}
    \begin{aligned}
      H^u
      =
      &H+u\left(\abs{A}_g^2+(R_{ab} \circ \phi)\nu^a \nu^b\right) \\
      &+u^2 \int_0^1 (t-1)
        \left[
        \vphantom{\frac{1}{2}}
          2{g_{tu}}^{\alpha \beta}{g_{tu}}^{\gamma \delta}{A^{tu}}_{\alpha \beta}
            (R_{abcd} \circ \phi[u]){\nu_{tu}}^a{\nu_{tu}}^d
              {{\phi_{tu}}^b}_{,\gamma} {{\phi_{tu}}^c}_{,\delta}
          \right. \\
      &\left.
          +2{A^{tu}}_{\alpha \beta}{A^{tu}}_{\gamma \delta}{A^{tu}}_{\epsilon \zeta}
            {g_{tu}}^{\beta \gamma}{g_{tu}}^{\delta \epsilon}{g_{tu}}^{\alpha \zeta}
          +\left(R_{ab|c} \circ \phi[u]\right){\nu_{tu}}^a{\nu_{tu}}^b{\nu_{tu}}^c
        \vphantom{\frac{1}{2}}
        \right] \, dt,
    \end{aligned}
  \end{equation}
establishing item (iv) and completing the proof.
\end{proof}

\subsection*{Toral regions}
Recalling \eqref{Tone}, \eqref{Ttwo}, \eqref{taudef}, and \eqref{RXYdef}
we define, for $1 \leq i \leq N$, the closed domains $\T_i \subset \R^2$ by
  \begin{equation}
  \label{Tipre}
    \T_i=\T_i[N,k,\ell,m,\zeta,\xi]
    :=
    \begin{cases}
      m\T_{X,Y,\sqrt{\tauj{1}}} \mbox{ if $i=1$} \\
      m\T_{X,Y,\sqrt{\tauj{N-1}}} \mbox{ if $i=N$} \\
      m\T_{X,Y,\sqrt{\tauj{i-1}}, \sqrt{\tauj{i}}} \mbox{ if $2 \leq i \leq N-1$},
    \end{cases}
  \end{equation}
so that $\T_i$ is a $\sqrt{2}\pi/k \times \sqrt{2}\pi/\ell$ rectangle
with one or two discs removed, each having radius of order $m\sqrt{\tauj{1}}$.
By virtue of the second line of \eqref{estimatesforforceestimate}
we see that $\T_i$ tends with large $m$ to
  \begin{equation}
  \label{That}
    \That_i
    =
    \That_i[N,k,\ell]
    :=
    \left.\left[-\frac{\pi}{\sqrt{2}k},\frac{\pi}{\sqrt{2}k}\right]
      \times \left[-\frac{\pi}{\sqrt{2}\ell},\frac{\pi}{\sqrt{2}\ell}\right]
    \right\backslash
      \begin{cases}
        \{(0,0)\} \mbox{ if } i \in \{1,N\} \\
        \left\{\pm \left(\frac{\pi}{2\sqrt{2}k},\frac{\pi}{2\sqrt{2}\ell}\right)\right\}
          \mbox{ if } 1<i<N.
      \end{cases}
  \end{equation}
Recalling \eqref{Phidef}, \eqref{Text}, \eqref{Tint}, \eqref{torheights}, and \eqref{catheights},
we also define the maps $T_i=T_i[N,k,\ell,m,\zeta,\xi]: \T_i \to \Sph^3$ by
  \begin{equation}
  \label{Tdef}
    \begin{aligned}
      &T_1(\xx,\yy)
      :=
      \Phi
        \left(
          T_{ext}\left[(0,0), \zzj{1}^K, \zzj{1}, R, X,Y, \tauj{1}\right]
            \left(\frac{\xx}{m},\frac{\yy}{m}\right)
        \right), \\
      &T_N(\xx,\yy)
      :=
      \Phi
        \left(
          T_{ext}\left[(N-2)(X,Y),\zzj{N-1}^K,\zzj{N},R,X,Y,\tauj{N-1}\right]
            \left(\frac{\xx}{m},\frac{\yy}{m}\right)
        \right), \mbox{ and otherwise} \\
      &T_i(\xx,\yy)
      :=
      \Phi
        \left(
          T_{int}
            \left[
              (2i-3)\left(\frac{X}{2},\frac{Y}{2}\right),
              \zzj{i-1}^K, \zzj{i}^K, \zzj{i}, R, X, Y,\tauj{i-1},\tauj{i}
            \right]
              \left(\frac{\xx}{m},\frac{\yy}{m}\right)
        \right).
    \end{aligned}
  \end{equation}
Then, referring to \eqref{initsurfdef} and particularly \eqref{Omega},
the image of each $T_i$ is entirely contained in
the initial surface $\Sigma[N,k,\ell,m,\zeta,\xi]$
and defines the corresponding toral region
  \begin{equation}
  \label{torrdef}
    \torr[i]
    =
    \torr[i;N,k,\ell,m,\zeta,\xi]
    :=
    T_i\left(\T_i\right),
  \end{equation}
so that $T_i$ a diffeomorphism onto its image.
Abusing notation slightly
we denote the inverse of this diffeomorphism
by $T_i^{-1}$.
From \eqref{taudef}, \eqref{torheights}, \eqref{catheights}, \eqref{Grpdef},
\eqref{Omega}, 
\eqref{initsurfdef}, \eqref{catrdef}, and \eqref{torrdef}
it is clear that
  \begin{equation}
    \begin{aligned}
      \torr[i] \cap \torr[i'] \neq \emptyset
        \mbox{ if and only if $i=i'$},
      &\qquad
      \torr[i] \cap \catr[i'] \neq \emptyset
        \mbox{ if and only if $i-i' \in \{0,1\}$}, \\
      \bigcup_{i=1}^N \Omega_i
      =
      \bigcup_{i=1}^{N-1} \catr[i] \cup \bigcup_{i=1}^N \torr[i],
      \qquad &\mbox{and} \qquad
      \Sigma
      =
      \bigcup_{i=1}^{N-1} \Grp \catr[i] \cup \bigcup_{i=1}^N \Grp \torr[i].
    \end{aligned}
  \end{equation}

Each limit region $\That_i \supset \T_i$
naturally carries the flat metric
$\geuc=d\xx^2+d\yy^2$,
but we also equip it with the conformal metric
  \begin{equation}
  \label{chihatidef}
    \begin{aligned}
      &\chihat_i
      =
      \chihat_i[N,k,\ell]
      :=
      \rhohat_i^2\geuc, \mbox{ having conformal factor} \\
      &\rhohat_i
      =
      \rhohat_i[N,k,\ell]: \That_i \to (0,\infty)
        \mbox{ defined by} \\
      &\rhohat_i(\xx,\yy)
      :=
      \cutoff{\frac{1}{10\ell}}{\frac{1}{5\ell}}(d_i(\xx,\yy))
        +\frac{1}{d_i(\xx,\yy)}
          \cdot \cutoff{\frac{1}{5\ell}}{\frac{1}{10\ell}}(d_i(\xx,\yy)),
          \mbox{ where} \\
      &d_i(\xx,\yy)
        \mbox{ is the Euclidean distance in $\R^2$ from the set}
      \begin{cases}
        \{(0,0)\} \mbox{ if $i \in \{1,N\}$} \\
        \left\{\pm \left(\frac{mX}{2},\frac{mY}{2}\right) \right\} \mbox{ if $1<i<N$}, \\
      \end{cases}
    \end{aligned}
  \end{equation}
recalling \eqref{cutoff}.
Under the $\chihat_i$ metric
$\That_i$ looks like
a flat $\sqrt{2}\pi/k \times \sqrt{2}\pi/\ell$ rectangle
with one or two discs of radius $1/5\ell$ replaced
by one or two infinite half-cylinders of radius $1$,
each attached smoothly along an annulus. 
We emphasize that $\rhohat_i$ is independent of $m$
as well as the parameters $\zeta$, $\xi$,
and, in view of \eqref{rhodef}, we observe that on each domain $\T_i$
  \begin{equation}
  \label{rhorhohat}
    T_i^*\rho=m\rhohat_i.
  \end{equation}

In the next section we will define the \emph{extended substitute kernel}
needed to complete the construction, as outlined in Section 1.
Then, in the final section, the role of the dislocations will become clear:
the dislocation $\Dcal_i$ (recalling \eqref{Ddef})
on the toral region $\torr[i]$ will be varied
to cancel the ``extended'' portion of the extended substitute kernel
supported there. For this reason it is necessary to isolate the dominant
contribution of each dislocation to the mean curvature, and to that end
we define
$v_i \in C^\infty(\That_i)$ by
  \begin{equation}
  \label{vdef}
    \begin{aligned}
      v_i(\xx,\yy)
      :=
      &\cutoff{\frac{1}{5\ell}}{\frac{1}{10\ell}}
        \left(\sqrt{\left(\xx-\frac{mX}{2}\right)^2
          + \left(\yy-\frac{mY}{2}\right)^2}\right) \\
      &- \cutoff{\frac{1}{5\ell}}{\frac{1}{10\ell}}
        \left(\sqrt{\left(\xx+\frac{mX}{2}\right)^2
          + \left(\yy+\frac{mY}{2}\right)^2}\right)
        \mbox{ for $2 \leq i \leq N-1$},
    \end{aligned}
  \end{equation}
and we define $\wbar_i \in C^\infty_{\Grp}(\Sigma)$
to be the unique $\Grp$-invariant (recalling \eqref{Grpdef} and \eqref{Grpaction})
function satisfying
  \begin{equation}
  \label{wbar}
    T_i^*\wbar_i
    :=
    \begin{cases}
      (-1)^{N-i}\left.\left(\Delta_{\chihat_i}v_i\right)\right|_{\T_i} \mbox{ if $1<i<N$} \\
      0 \mbox{ if $i \in \{1,N\}$},
    \end{cases}
    \quad \mbox{and} \quad
    \wbar_i|_{\Sigma \backslash \Grp \torr[i]}:=0,
  \end{equation}
the alternating sign included to account for the alternating direction
of the unit normal on the toral regions
and
the exceptional cases $i=1$ and $i=N$ included merely for convenience of notation.
(We could have alternatively built the alternating sign
into the definition of the dislocations.)
The function $v_i|_{\T_i} \circ T_i^{-1}$
should be regarded
as the section of the normal bundle
graphically generating dislocations on $\torr[i]$,
and in the following proposition
we will see that the function $\wbar_i$ then captures the principal
effect of dislocation on the mean curvature.
Later the collection $\{\wbar_i\}_{i=2}^{N-1}$ will reappear as the defining
basis for the extended part of the extended substitute kernel.
Right now we estimate the geometry of the toral regions.

\begin{prop}[Estimates of the geometry of the toral regions]
\label{torest}
Given a real number $c>0$
and integers $N \geq 2$, $\ell \geq k \geq 1$, and $j \geq 0$,
there exist real numbers
$m_0=m_0[N,k,\ell,c]>0$
and $C=C[N,k,\ell,j]>0$
such that whenever $\zeta,\xi \in [-c,c]^{N-1}$ and $m>m_0$,
for $1 \leq i \leq N$
  \begin{enumerate}[(i)]
    \item $\norm{\chi - {T_i^{-1}}^*\chihat_i:
            C^j\left(T^*\torr[i]^{\otimes 2},\chi \right)}
            \leq Cm^2\tauj{1}$;

    \item $\norm{\rho: C^j(\torr[i],\chi,\rho)}
           +\norm{\rho^{-1}: C^j(\torr[i],\chi,\rho)}
           \leq
           C$,

    \item $\norm{A-(-1)^{N-i}m^{-2}{T_i^{-1}}^*(d\yy^2-d\xx^2)
           :C^j\left(\torr[i] \backslash \left(\catr[i-1] \cup \catr[i]\right),\chi\right)}
           \leq Cm\tauj{1}$,

    \item $\norm{\rho^{-2}\abs{A}_g^2 : C^j\left(\torr[i],\chi \right)}
      \leq Cm^{-2}$, and

    \item $\norm{\rho^{-2}H - \Dcal_i\wbar_i:
      C^j\left(\torr[i],\chi,m^2\rho^{-2}\tauj{1}+m^2\tauj{1}^2\right)} \leq C$,
  \end{enumerate}
recalling \eqref{taudef}, \eqref{Ddef},
\eqref{weightedhoelderdef}, \eqref{gdef}, \eqref{AHdef},
\eqref{chidef}, \eqref{rhodef}, \eqref{Tdef}, \eqref{torrdef},
\eqref{chihatidef}, and \eqref{wbar}.
\end{prop}

\begin{proof}

We first observe,
using \eqref{chihatK}, \eqref{catdef}, \eqref{rhocat}, 
\eqref{Tdef}, \eqref{chihatidef}, and \eqref{rhorhohat},
that if $\torr[i] \cap \catr[i'] \neq \emptyset$,
then
  \begin{equation}
    \left.\chihatK \right|_{\cat_{i'}^{-1}(\torr[i])}
      -\left. \cat_{i'}^* {T_i^{-1}}^* \chihat_i \right|_{\torr[i] \cap \catr[i']}
    =
    \sech^2 t \, dt^2
    =
    \tauj{i'}^2 \left(\rho \circ \cat_{i'}\right)^2 \, dt^2, \\
  \end{equation}
but
$\rho|_{\cat_{i'}^{-1}(\torr[i])} \leq \tauj{i'}^{-1/2}$
by \eqref{torrdef},
so by applying items (i) and (ii) of Proposition \ref{catest}
we have proven item (i) of the present proposition
on the overlap of the toral and catenoidal regions.
On this same intersection the remaining items (with (iii) obviously excluded)
also follow from the corresponding ones in Proposition \ref{catest}.

We will finish the proof by verifying the estimates
on
$\torr[i] \cap \{\rho \leq 10 \ell m\}
 =\torr[i] \backslash \left(\catr[i-1] \cup \catr[i]\right)$
(understanding $\catr[0]=\catr[N]=\emptyset$).
For this we set
  \begin{equation}
    D_i:=T_i^{-1}\left(\{\rho \leq 10 \ell m\}\right) \subset \T_i
  \end{equation}
and apply Lemma \ref{perimm},
viewing $T_i|_{D_i}$ as a perturbation $\phi_i[\uj{i}]$
of the embedding
$\phi_i:=\varpi \circ T_i: D_i \to (\Sph^3,m^2\gsph)$
of $D_i$ into the Clifford torus $\T$
with $m^2\gsph$ unit normal $\nu[\phi_i]$ directed toward $C_1$;
here $\varpi: \Sph^3 \backslash (C_1 \cup C_2) \to \T$
is nearest-point projection
in $(\Sph^3,m^2\gsph)$ onto $\T$
and the function $\uj{i}$ generating the perturbation
is identified below.
Thus, recalling \eqref{Phidef} and \eqref{Phullback},
  \begin{equation}
    \begin{aligned}
      &\varpi \left(\Phi(\xx,\yy,\zz)\right)=\Phi(\xx,\yy,0),
      \qquad
      \phi_i(\xx,\yy)=\Phi\left(\frac{\xx}{m}+\xx_i,\frac{\yy}{m}+\yy_i,0\right),
      \qquad
      \phi_i^*m^2\gsph=\geuc, \\
      &\nu[\phi_i]|_{(\xx,\yy)}
        =
        \left.m^{-1}\Phi_*\partial_{\zz}\right|_{\left(\frac{\xx}{m}+\xx_i,\frac{\yy}{m}+\yy_i,0\right)}
      \quad \mbox{and} \quad
      \nu[\phi_i[\uj{i}]]|_{(\xx,\yy)}
        =
        (-1)^{N-i}m^{-1}\nu|_{T_i(\xx,\yy)},
    \end{aligned}
  \end{equation}
where $\xx_i,\yy_i$ give the appropriate lattice site appearing in \eqref{Tdef},
$\nu[\phi_i[\uj{i}]]$ is the $m^2\gsph$ unit normal for $\phi_i[\uj{i}]$
specified in Lemma \ref{perimm},
and $\nu$ is the $\gsph$ unit normal we chose for $\Sigma$ just above \eqref{AHdef}.
Writing $A[\uj{i}]$ and $H[\uj{i}]$
for the second fundamental form and mean curvature of $\phi_i[\uj{i}]$
relative to $m^2\gsph$ and $\nu[\phi_i[\uj{i}]]$,
as in Lemma \ref{perimm},
and recalling \eqref{AHdef} and the definition of $\iota: \Sigma \to \Sph^3$
as the inclusion map of the initial surface in $\Sph^3$,
it follows that
  \begin{equation}
  \label{torscale}
    \begin{aligned}
      &\iota^*\gsph|_{T_i(D_i)}=m^{-2}\phi_i[\uj{i}]^*m^2\gsph, \\
      &A|_{T_i(D_i)}=(-1)^{N-i}m^{-1}A[\uj{i}], \mbox{ and} \\
      &H|_{T_i(D_i)}=(-1)^{N-i}mH[\uj{i}].
    \end{aligned}
  \end{equation}

Referring to \eqref{Tdef} and the supporting definitions
(including in particular \eqref{torheights} and \eqref{catheights})
and recalling \eqref{Ddef},
we have as the function generating the perturbation
(that is playing the role of $u$ in the statement of Lemma \ref{perimm})
$\uj{i}: D_i \to \R$ given by
  \begin{equation}
      m^{-1}\uj{1}(\xx,\yy)
      =
      \zzj{1} +
        \tauj{1} \left(
          \ln \frac{1}{10 \ell m \tauj{1}} - \arcosh \frac{\rrj{0}(\xx,\yy)}{m\tauj{1}}
        \right)
        \cutoff{\frac{1}{5\ell}}{\frac{1}{10\ell}}(\rrj{0}(\xx,\yy)),
  \end{equation}
  \begin{equation}
      m^{-1}\uj{N}(\xx,\yy)
      =
      \zzj{N} +
        \tauj{N-1}\left(
          \arcosh \frac{\rrj{0}(\xx,\yy)}{m\tauj{N-1}} - \ln \frac{1}{10\ell m \tauj{N-1}}
        \right)
        \cutoff{\frac{1}{5\ell}}{\frac{1}{10\ell}}(\rrj{0}(\xx,\yy)),
  \end{equation}
and for $1<i<N$
  \begin{equation}
    \begin{aligned}
      m^{-1}\uj{i}(\xx,\yy)
        =
        \zzj{i} &+
        \left(
          -\Dcal_i
          +\tauj{i-1}\arcosh \frac{\rrj{-1}(\xx,\yy)}{m\tauj{i-1}}
          -\tauj{i-1} \ln \frac{1}{10\ell m \tauj{i-1}}
        \right)
        \cutoff{\frac{1}{5\ell}}{\frac{1}{10\ell}}(\rrj{-1}(\xx,\yy)) \\
        &
        +
        \left(
          \Dcal_i
          +\tauj{i}\ln \frac{1}{10\ell m \tauj{i}}
          -\tauj{i} \arcosh \frac{\rrj{1}(\xx,\yy)}{m\tauj{i}}
        \right)
        \cutoff{\frac{1}{5\ell}}{\frac{1}{10\ell}}(\rrj{1}(\xx,\yy)),
    \end{aligned}
  \end{equation}
where
  \begin{equation}
    \rrj{0}(\xx,\yy):=\sqrt{\xx^2+\yy^2}
    \qquad \mbox{and} \qquad
    \rrj{\pm 1}(\xx,\yy):=\rrj{0}\left(\xx \mp \frac{mX}{2},\yy \mp \frac{mY}{2}\right).
  \end{equation}

Using the estimates
  \begin{equation}
    \begin{aligned}
    \label{estsforuests}
      &\norm{\cutoff{\frac{1}{5\ell}}{\frac{1}{10\ell}} \circ \rrj{0}:
        C^j(T^*{D_1}^{\otimes j},\geuc)}
      \leq
      C[\ell,j] \\
      &\norm{ D[\geuc]^j \arcosh \frac{\rrj{0}}{m\tauj{i}}:
        C^0({T^*D_1}^{\otimes j},\geuc)}
      \leq
      C[\ell,j] \mbox{ whenever $j>0$}, \\
      &\sum_{i=1}^{N-1} \norm{\tauj{i} \arcosh \frac{\rrj{0}}{m\tauj{i}}:
        C^0(D_1)}
      \leq C[N,k,\ell]m^2\tauj{1}
        \mbox{ (by \eqref{cisummary} and \eqref{taudef})}, \\
      &\eqref{arcoshln}, \mbox{ the last item of \eqref{estimatesforforceestimate}, and }
        \abs{\Dcal_i} \leq C[N,k,\ell]m\tauj{1}
        \mbox{ (by \eqref{cisummary} and \eqref{taudef} whenever $m \geq c$)},
    \end{aligned}
  \end{equation}
we obtain
  \begin{equation}
  \label{ujest}
    \begin{aligned}
      &\norm{\uj{i}: C^0(D_i)} \leq C[N,k,\ell]m^3\tauj{1}
        \mbox{ and} \\
      &\norm{D[\geuc]^j\uj{i}: C^0({T^*D_i}^{\otimes j},\geuc)} \leq C[N,k,\ell,j]m^2\tauj{1}
        \mbox{ for each $j \geq 1$}.
    \end{aligned}
  \end{equation}

  Using \eqref{Phullback} and recalling notation of Lemma \ref{perimm}
  (remembering in particular that we are taking $g$ in its statement
  to be $m^2\gsph$) we also have
  \begin{equation}
  \label{guest}
    \begin{aligned}
      &g^{\uj{i}}=(1+\sin 2m^{-1}\uj{i}) \, d\xx^2 + (1-\sin 2m^{-1}\uj{i}) \, d\yy^2,
        \mbox{ so by \eqref{ujest}} \\
      &\norm{g^{\uj{i}}-\geuc: C^j\left(T^*D_i^{\otimes 2},\geuc \right)}
      \leq
      C[N,k,\ell,j]m^2\tauj{1}.
    \end{aligned}
  \end{equation}
Item (ii) of Lemma \ref{perimm} now yields
  \begin{equation}
    \norm{m^2T_i^*\gsph - \geuc: C^j\left(T^*D_i^{\otimes 2},\geuc \right)}
    \leq
    C[N,k,\ell,j]m^2\tauj{1}.
  \end{equation}
The proofs of (i) and (ii) are now completed by
\eqref{chihatidef}, \eqref{rhorhohat},
and the observation that
  \begin{equation}
    \norm{\rhohat_i: C^j(D_i,\geuc)}
      +\norm{\rhohat_i^{-1}: C^j(D_i,\geuc)}
    \leq
    C[\ell,j].
  \end{equation}
Furthermore,
again using \eqref{Phullback} and the notation of Lemma \ref{perimm},
  \begin{equation}
  \label{Auest}
    \begin{aligned}
      &A^{\uj{i}}=-m^{-1}\left(\cos 2m^{-1}\uj{i}\right) \, \left(d\xx^2-d\yy^2\right),
        \mbox{ so by \eqref{ujest}}, \\
      &\norm{A^{\uj{i}}-m^{-1}(d\yy^2-d\xx^2): C^j\left(T^*D_i^{\otimes 2},\geuc \right)}
      \leq
      C[N,k,\ell,j]m\tauj{1},
    \end{aligned}
  \end{equation}
while from \eqref{ujest} and \eqref{guest}
  \begin{equation}
  \label{Duguest}
    \norm{D[g^{\uj{i}}]^2\uj{i}: C^j\left(T^*D_i^{\otimes 2},\geuc \right)}
      +\norm{\abs{d\uj{i}}_{g^{\uj{i}}} : C^j\left(D_i,\geuc\right)}
    \leq
    C[N,k,\ell,j]m^2\tauj{1}.
  \end{equation}
Item (iii) of the present proposition
is now proved by applying item (i), \eqref{Auest}, and \eqref{Duguest}
(and \eqref{ujest} again) in item (iii) of Lemma \ref{perimm},
keeping in mind \eqref{torscale}.
Item (iv) follows in turn.

Finally,
from the identity
  \begin{equation}
    \Delta_{\geuc} \arcosh \frac{\rrj{0}}{m\tauj{i}}
    =
    -\frac{m^2\tauj{i}^2}{\rrj{0}\left(\rrj{0}^2-m^2\tauj{i}^2\right)^{3/2}}
  \end{equation}
along with \eqref{arcoshln},
the first two estimates of \eqref{estsforuests},
definition \eqref{wbar},
and the fact that $\Delta_{\chihat_i}=\rhohat_i^{-2}\Delta_{\geuc}$
(by \eqref{chihatidef} and the two-dimensionality of $\That_i$)
we find
  \begin{equation}
  \label{Lapest}
    \norm{\Delta_{\geuc}\uj{i} - (-1)^{N-i}m\rhohat_i^2\Dcal_i T_i^*\wbar_i:
      C^j\left(D_i, \geuc\right)}
    \leq
    C[N,k,\ell,j]m\tauj{1}.
  \end{equation}
(Note that without subtracting the dislocation term on the left
it would be necessary
to allow $C$ on the right-hand side of \eqref{Lapest}
to depend on $c$,
or, if we were to apply the assumption we have used repeatedly
that $c \leq m$,
to allow the exponent on $m$ to increase to $3$.)
We now apply item (iv) of Lemma \ref{perimm}.
In doing so we make use of
\eqref{Lapest},
\eqref{ujest}, \eqref{guest}, \eqref{Auest}, and \eqref{Duguest};
we also take note of Remark \ref{spherical}
and of course the facts that $\T$ itself is minimal
and $m^2\tauj{1}<1$.
We thereby obtain
  \begin{equation}
    \norm{H[\uj{i}]-(-1)^{N-i} m\rhohat_i^2 \Dcal_i T_i^*\wbar_i: C^j(D_i,\geuc)}
    \leq
    C[N,k,\ell,j]m\tauj{1},
  \end{equation}
and the proof is completed by \eqref{rhorhohat} and \eqref{torscale}.
\end{proof}

\subsection*{Decay norms and a global estimate of the mean curvature}
As mentioned in Section 1, because the characteristic size $\tauj{1}$
of the catenoidal waists is so much smaller than the characteristic size $m^{-1}$
of the toral regions, we must allow perturbing functions to be much larger
on the toral regions than on the core of the catenoidal regions.
For this reason we will weight our norms by powers of the factor $m \rho^{-1}$,
which takes the value $1$ a maximal distance from the catenoidal regions
and is of order $m\tauj{1}$ at the waists. Specifically, for each $\alpha \in (0,1)$,
$\gamma \in [0,\infty)$, and nonnegative integer $j$, we define the norm
  \begin{equation}
  \label{globalweightednorm}
    \norm{\cdot}_{j,\alpha,\gamma}
    =
    \norm{\cdot}_{C^{j,\alpha,\gamma}(\Sigma)}
    :=
    \norm{\cdot: C^{j,\alpha}\left(\Sigma,\chi, \frac{m^\gamma}{\rho^\gamma}\right)},
  \end{equation}
(recalling \eqref{weightedhoelderdef} and \eqref{rhodef})
and the corresponding Banach space
along with its (closed) $\Grp$-invariant subspace
(recalling \eqref{Grpdef} and \eqref{Grpaction})
  \begin{equation}
  \label{globalweightedspace}
    \begin{aligned}
      &C^{j,\alpha,\gamma}(\Sigma)
      :=
      \left. \left\{
        u \in C^{j,\alpha}(\Sigma,\chi) \; \right| \;
        \norm{u}_{j,\alpha,\gamma}<\infty
      \right\} \mbox{ and} \\
      &C^{j,\alpha,\gamma}_{\Grp}(\Sigma)
      :=
      \left. \left\{
        u \in C^{j,\alpha,\gamma}(\Sigma) \; \right| \;
      \mathfrak{g}.u=u \circ \mathfrak{g}^{-1}
      \mbox{ for all } \mathfrak{g} \in \Grp
      \right\},
    \end{aligned}
  \end{equation}
in accordance with Notation \ref{functionspaceappendages}.

\begin{remark}
\label{topspaces}
Of course, since each initial surface $\Sigma$ is compact,
$C^{j,\alpha,\gamma}(\Sigma)$,
$C^{j,\alpha}(\Sigma,\chi)$,
and $C^{j,\alpha}_{loc}(\Sigma)$
all refer to the same topological vector space,
which we more simply call $C^{j,\alpha}(\Sigma)$
(forgetting the norm structures of the first two spaces
and dropping the superfluous subscript of the third).
\end{remark}

\begin{definition}[Continuity in the parameters]
\label{functionalcty}
If $f=f[N,k,\ell,m,\zeta,\xi] \in C^{j,\alpha}(\Sigma[N,k,\ell,m,\zeta,\xi])$
is a family of functions on the initial surfaces
and we make the usual assumption that
$\zeta,\xi \in [-c,c]^{N-1}$ for some $c>0$,
we say that $f$ \emph{depends continuously on $(\zeta,\xi)$}
if the map
$(\zeta,\xi) \to f[N,k,\ell,m,\zeta,\xi] \circ I[N,k,\ell,m](\zeta,\xi,\cdot)$
is continuous from
$[-c,c]^{N-1} \times [-c,c]^{N-1}$
to $C^{j,\alpha}(\Sigma[N,k,\ell,m,0,0])$,
where $I=I[N,k,\ell,m]$ is as described in Remark \ref{initsurfcty}.
Note that this definition does not depend on the particular choice of $I$.

Similarly, if
$\Acal[\zeta,\xi]=\Acal[N,k,\ell,m,\zeta,\xi]:
 C^{J,\alpha}(\Sigma[N,k,\ell,m,\zeta,\xi])
 \to
 C^{j,\alpha}(\Sigma[N,k,\ell,m,\zeta,\xi])$
 (or $\R$)
is a family of (not necessarily linear) continuous maps,
we call the associated map 
$(u,\zeta,\xi) \mapsto \Acal[\zeta,\xi][u]$
\emph{continuous} for fixed $N$, $k$, $\ell$, and $m$
if the map
  \begin{equation}
    \begin{aligned}
      C^{J,\alpha}(\Sigma[N,k,\ell,m,0,0]) \times [-c,c]^{N-1} \times [-c,c]^{N-1}
      \to
      C^{j,\alpha}(\Sigma[N,k,\ell,m,0,0]) \mbox{ (or $\R$)} \\
      (u,\zeta,\xi)
      \mapsto 
      I(\zeta,\xi,\cdot)^*\left(\Acal[\zeta,\xi]\left({I(\zeta,\xi,\cdot)^{-1}}^*u\right)\right)
        \mbox{ (or $\Acal[\zeta,\xi]\left({I(\zeta,\xi,\cdot)^{-1}}^*u\right)$)}
    \end{aligned}
  \end{equation}
is continuous.
\end{definition}

In order to secure acceptable decay estimates
for solutions to the linearized problem
we will need the following estimate for the initial mean curvature.

\begin{cor}[Global weighted estimate of the initial mean curvature]
\label{Hest}
Given real numbers $c>0$ and $\gamma \in (0,1)$
as well as integers $N \geq 2$ and $\ell \geq k \geq 1$,
there exist $C=C[N,k,\ell]>0$
and $m_0=m_0[N,k,\ell,c,\gamma]$
such that for each integer $m>m_0$, each $\zeta, \xi \in [-c,c]^{N-1}$,
and each $\alpha \in (0,1)$
  \begin{equation}
    \norm{\rho^{-2}H-\sum_{i=2}^{N-1} \Dcal_i \wbar_i}_{0,\alpha,\gamma}
    \leq
    C\tauj{1},
  \end{equation}
using the norm \eqref{globalweightednorm}
and recalling \eqref{Ddef}, \eqref{AHdef}, \eqref{rhodef}, and \eqref{wbar}.
Moreover $\rho^{-2}H$ is $\Grp$-invariant
(recalling \eqref{Grpdef} and \eqref{Grpaction})
and depends continuously,
as an element of $C^{0,\alpha,\gamma}(\Sigma)$,
on $(\zeta,\xi)$
in the sense of Definition \ref{functionalcty}.
\end{cor}

\begin{proof}
The continuity claim is obvious,
as in fact both $\rho \circ I$ and $H \circ I$ are manifestly smooth.
The $\Grp$-invariance of $\rho$ also follows directly from its definition,
while that of $H$ follows from the $\Grp$-invariance of $\Sigma$ itself,
establishing that $\rho^{-2}H$ is $\Grp$-invariant as well.
(Note that in this construction all elements of our symmetry group $\Grp$
act on all functions under consideration according to the first line of \eqref{Grpaction}.
Of course, were we to enforce also symmetries reversing the sides of $\Sigma$,
it would be natural to consider a different action of $\Grp$ on $\rho$
from that defined by \eqref{Grpaction} on $H$,
since the former represents a true scalar field,
while $H$ represents a section of the normal bundle.
Specifically, the appropriate action on $\rho$ would be to follow the first line
of \eqref{Grpaction} even for elements reversing the sides of $\Sigma$.
In this case too we would conclude the appropriate $\Grp$-equivariance of $\rho^{-2}H$.)
As for the estimate, from item (v) of Proposition \ref{torest} 
and items (iii) and (vi) of Proposition \ref{catest}
we get
  \begin{equation}
  \label{C1Hest}
    \norm{\rho^{-2}H-\sum_{i=2}^{N-1} \Dcal_i \wbar_i:
      C^1\left(\Sigma, \chi,m^2\rho^{-2}\tauj{1}+m^2\tauj{1}^2\right)}
      \leq
      C[N,k,\ell]
  \end{equation}
for some constant $C[N,k,\ell]>0$ whenever $m$ is sufficiently large
in terms of $N$, $k$, $\ell$, and $c$,
but
  \begin{equation}
    \frac{m^2\rho^{-2}\tauj{1}+m^2\tauj{1}^2}{m^\gamma \rho^{-\gamma} \tauj{1}}
    =
    \left(\frac{m}{\rho}\right)^{2-\gamma} + m^{2-\gamma}\tauj{1}\rho^\gamma
    \leq
    \left(\frac{m}{m}\right)^{2-\gamma} + m^{2-\gamma}\tauj{1}^{1-\gamma},
  \end{equation}
using \eqref{rhobounds} for the last inequality,
and the estimate now follows
from the second item of \eqref{estimatesforforceestimate}.
\end{proof}

\section{The linearized operator}
\label{linop}

We continue to write
$\iota=\iota[N,k,\ell,m,\zeta,\xi]: \Sigma[N,k,\ell,m,\zeta,\xi] \to \Sph^3$
for the inclusion map of the initial surface
$\Sigma=\Sigma[N,k,\ell,m,\zeta,\xi]$ in $\Sph^3$
and $\nu=\nu[N,k,\ell,m,\zeta,\xi]: \Sigma \to \iota^*T\Sph^3$
for the unit normal which points toward $C_1$ at the points of $\Sigma$
nearest to $C_1$
(or equivalently which points upward at the top of $\Sigma$
as viewed via coordinates obtained through the map $\Phi$ defined in \eqref{Phidef}).
Fixing the data $N \geq 2$, $\ell \geq k \geq 1$,
and $m$ sufficiently large,
we consider deformations of $\iota[N,k,\ell,m,0,0]$
obtained by varying the parameters $\zeta$ and $\xi$
and by additionally perturbing the resulting initial surface $\Sigma[N,k,\ell,m,\zeta,\xi]$
in the normal direction by a prescribed function.
Specifically we define
$\iota[N,k,\ell,m,\zeta,\xi,u]:
 \Sigma[N,k,\ell,m,\zeta,\xi] \to \Sph^3$
by
  \begin{equation}
  \label{deformediota}
    \iota[u]
    =
    \iota[N,k,\ell,m,\zeta,\xi,u](p)
    :=
    \exp_{\iota[N,k,\ell,m,\zeta,\xi](p)} u(p)\nu[N,k,\ell,m,\zeta,\xi](p),
  \end{equation}
where $\exp: T\Sph^3 \to \Sph^3$ is the exponential map for $(\Sph^3,\gsph)$.
As asserted in Lemma \ref{perimm},
$\iota[N,k,\ell,m,\zeta,\xi,u]$ is an immersion for sufficiently small $u$.
In this case
we write $\nu[N,k,\ell,m,\zeta,\xi,u]: \Sigma \to \iota[\zeta,\xi,u]^*T\Sph^3$
for the unit normal of $\iota[N,k,\ell,m,\zeta,\xi,u]$
whose value at each $p \in \Sigma$
has positive inner product with the vector
$\frac{d}{dt} \exp_{\iota[N,k,\ell,m,\zeta,\xi](p)} t\nu[N,k,\ell,m,\zeta,\xi](p)$
and
we write
$\Hcal[u]=\Hcal[u,\zeta,\xi]=\Hcal[N,k,\ell,m,\zeta,\xi,u]:
 \Sigma[N,k,\ell,m,\zeta,\xi] \to \R$ for the corresponding
mean curvature
  \begin{equation}
  \label{Hcaldef}
    \Hcal[u]
    =
    \Hcal[u,\zeta,\xi]
    =
    \Hcal[N,k,\ell,m,\zeta,\xi,u]
    :=
    H[\iota[N,k,\ell,m,\zeta,\xi,u],\nu[N,k,\ell,m,\zeta,\xi,u]],
  \end{equation}
with the notation and conventions introduced just below \eqref{curvatureconvention}.

Our goal is to find $\zeta,\xi \in \R^{N-1}$
and $u \in C^\infty_{\Grp}(\Sigma)$
(recalling Notation \ref{functionspaceappendages}) solving
  \begin{equation}
  \label{vanishingH}
    \Hcal[N,k,\ell,m,\zeta,\xi,u]=0
  \end{equation}
for each given $N \geq 2$, $\ell \geq k \geq 1$,
and $m$ sufficiently large,
with $u$ small enough that the resulting minimal surface
(the image of $\iota[N,k,\ell,m,\zeta,\xi,u]$)
is a small perturbation of the initial surface $\Sigma[N,k,\ell,m,\zeta,\xi]$,
so in particular embedded.
A major step toward the solution of \eqref{vanishingH} consists in the study of
the initial surface's Jacobi operator
$\Lcal=\Lcal[N,k,\ell,m,\zeta,\xi]$ defined by
  \begin{equation}
  \label{Lcaldef}
    \Lcal u = \left. \frac{d}{dt} \right|_{t=0} \Hcal[N,k,\ell,m,\zeta,\xi,tu]
    = \left( \Delta_g + \abs{A}_g^2 + 2 \right) u,
  \end{equation}
recalling that $g=\iota[N,k,\ell,m,\zeta,\xi]^*\gsph$.
Actually, because of the uniformity afforded by the $\chi$ metric \eqref{chidef},
it is much more convenient to study instead
the linear operator
  \begin{equation}
  \label{Lchidef}
    \Lchi
    =
    \Lchi[N,k,\ell,m,\zeta,\xi]
    :=
    \rho^{-2}\Lcal
    =
    \Delta_\chi + \rho^{-2}\abs{A}_g^2 + 2\rho^{-2},
  \end{equation}
which clearly takes $\Grp$-equivariant functions
(as defined by \eqref{Grpdef} and \eqref{Grpaction}) to $\Grp$-equivariant functions
and which, by virtue of the estimates of $\rho^{-2}\abs{A}_g^2$
in Propositions \ref{catest} and \ref{torest},
defines (for any $\alpha,\gamma \in (0,1)$) a linear map 
$\Lchi: C^{2,\alpha,\gamma}_{\Grp}(\Sigma) \to C^{0,\alpha,\gamma}_{\Grp}(\Sigma)$
bounded independently of $m$ and $c$.

In this section we construct a likewise bounded right inverse $\Rcal$ to $\Lchi$,
modulo the extended substitute kernel described in Section 1
and formally defined below.
We do this by first analyzing
$\Lchi$ ``semilocally'', meaning on the toral and catenoidal regions
individually,
and by observing that on each of these regions
$\Lchi$ has a simple limit as $m \to \infty$.
Significantly, because adjacent toral and catenoidal regions overlap,
when attempting to solve the equation $\Lchi u = f$
on a toral region $\torr[i]$,
we may assume that $f$ is supported
outside the intersection of $\torr[i]$ with the adjoining catenoidal region(s).
We will find we can invert these regional limits of $\Lchi$
(modulo extended substitute kernel in the toral cases)
and so produce approximate semilocal inverses to $\Lchi$,
which will be applied iteratively,
using decay properties of the solutions they yield,
to construct $\Rcal$.

\subsection*{Approximate solutions on the catenoidal regions}
Recalling \eqref{cyldef} and \eqref{chihatK},
we define the operator
  \begin{equation}
  \label{LhatKdef}
    \LhatK
    :=
    \Delta_{\chihatK} + 2 \sech^2 t
  \end{equation}
on functions on $\cyl$.
Note that $\LhatK$ is simply $\cosh^2 t$ times the Jacobi operator
of the standard catenoid \eqref{stdcat}.
From items (i) and (v) of Proposition \ref{catest} we see that
(recalling \eqref{Lchidef})
  \begin{equation}
  \label{catoplimit}
    \lim_{m \to \infty} {\cat_i}^* \Lchi {{\cat_i}^*}^{-1}
    =
    \LhatK,
  \end{equation}
where the convergence is to be interpreted in the following sense.
For any given bounded subset $\Omega$ of the cylinder $\cyl$
the operator on the left-hand side of \eqref{catoplimit}
is defined as a map $C^2_{loc}(\Omega) \to C^0_{loc}(\Omega)$
whenever $m$ is taken sufficiently large
in terms of the diameter of $\Omega$ and $\abs{\zeta}$
(noting $\lim_{m \to \infty} a_i=\infty$ by \eqref{adef})
and its difference from the operator on the right-hand side
is a first-order operator $X[m]+f[m]$ on $\Omega$
satisfying
$\lim_{m \to \infty}
 \left(\norm{X[m]: C^j(T\Omega,\chihatK)}+\norm{f[m]: C^j(\Omega,\chihatK)}\right)
 =0$
for each nonnegative integer $j$.

Recall that each catenoidal region $\catr[i]$
is defined in \eqref{catrdef} via \eqref{catdef}
as the image under $\Phi$ \eqref{Phidef}
of a certain catenoid in $\R^3$.
Of course this last catenoid has an axis of symmetry,
a line in $\R^3$
whose intersection with the domain of $\Phi$ has image under $\Phi$
a quarter great circle in $\Sph^3$,
which circle (at least in this paragraph)
we will call \emph{the axis of $\catr[i]$}.
It follows from
\eqref{Phintertwine}, \eqref{Grpdef}, and \eqref{catrdef}
that the subgroup of $\Grp$
preserving a given catenoidal region $\catr[i]$ as a set
is isomorphic to the dihedral group $D_2$ of order $4$
(also isomorphic to $\Z_2 \times \Z_2$ of course,
but we favor the more concise and geometric nomenclature),
consisting of
(i) the identity element $I$ of $O(4)$,
(ii) reflection $\xbar_i$ through the great sphere
     containing $C_2$ and the axis of $\catr[i]$,
(iii) reflection $\ybar_i$ through the great sphere
      containing $C_1$ and the axis of $\catr[i]$, and
(iv) rotation $\xbar_i\ybar_i=\ybar_i\xbar_i$
     through angle $\pi$ (also called reflection)
     through the axis of $\catr[i]$.

Using \eqref{Phintertwine} again,
we see that $\cat_i$ intertwines
the above $D_2$ action on $\catr[i]$
with the natural action of the $D_2$ subgroup 
of symmetries of $(\cyl,\chihatK)$
  \begin{equation}
  \label{GrpKdef}
      \GrpK
      :=
      \left\{
        \widehat{I}_K, \xbarhat_K, \ybarhat_K,
        \xbarhat_K\ybarhat_K
      \right\},
      \mbox{ where }
      \widehat{I}_K(t,\theta):=(t,\theta), \mbox{ }
      \xbarhat_K(t,\theta):=(t,\pi-\theta), \mbox{ and }
      \ybarhat_K(t,\theta):=(t,-\theta),
  \end{equation}
in the sense that
  \begin{equation}
    \cat_i \circ \xbarhat_K = \xbar_i \circ \cat_i
    \qquad \mbox{and} \qquad
    \cat_i \circ \ybarhat_K = \ybar_i \circ \cat_i
  \end{equation}
(and these elements generate the two groups).
Since $I$, $\xbarhat$, and $\ybarhat$ all preserve each side
(choice of unit normal) of $\Sigma$,
the natural action (recalling \eqref{Grpaction})
of any element $\mathfrak{g} \in \GrpK$ on a function $f$ on $\cyl$
is simply $\mathfrak{g}.f=f \circ \mathfrak{g}$.

Next, having in mind \eqref{globalweightednorm}
and using \eqref{cisummary}, \eqref{taudef}, and \eqref{rhocat},
we also note that on each $\cyl_{a_i}$
  \begin{equation}
    C[N,k,\ell]^{-1}m\tauj{1}
    \leq
    \frac{m{\cat_i}^*\rho^{-1}}{e^{\abs{t}}}
    \leq
    C[N,k,\ell]m\tauj{1},
  \end{equation}
so that the pullback to $\cyl_{a_i}$ by $\cat_i$
of the weight $m^\gamma \rho^{-\gamma}$
appearing in our global norm \eqref{globalweightednorm} on $\Sigma$
is comparable to the weight $m^\gamma \tauj{1}^\gamma e^{\gamma \abs{t}}$ on $\cyl$,
where obviously the first two factors are constant on $\cyl$.
All the above considerations motivate us to introduce,
for each nonnegative integer $j$ and $\alpha \in [0,1)$ and $\gamma \in (0,1)$
the norms
  \begin{equation}
  \label{catnorm}
    \norm{\cdot}_{C^{j,\alpha,\gamma}(\cyl)}
    :=
    \norm{\cdot, C^{j,\alpha}\left(\cyl,\chihatK,e^{\gamma \abs{t}}\right)}
  \end{equation}
(recalling \eqref{weightedhoelderdef})
and corresponding Banach spaces of $\GrpK$-even functions
  \begin{equation}
  \label{catspace}
    C_{\GrpK}^{j,\alpha,\gamma}(\cyl)
    :=
    \left. \left\{
      u \in C^{j,\alpha}(\cyl,\chihatK) \; \right| \;
      \norm{u}_{C^{j,\alpha,\gamma}(\cyl)}<\infty
      \mbox{ and }
      \mathfrak{g}.u=u \circ \mathfrak{g}
      \mbox{ for all } \mathfrak{g} \in \GrpK
    \right\}.
  \end{equation}
Clearly
$\LhatK: C^{2,\alpha,\gamma}_{\GrpK}(\cyl) \to C_{\GrpK}^{0,\alpha,\gamma}(\cyl)$
is bounded (independently of $\alpha, \gamma \in (0,1)$).
The following proposition presents a suitable inverse.

\begin{prop}[Solutions to the model problem on the catenoid]
\label{catsol}
There exists a linear map
  \begin{equation}
    \RhatK: C^{0,\alpha,\gamma}_{\GrpK}(\cyl) \to C^{2,\alpha,\gamma}_{\GrpK}(\cyl)
  \end{equation}
defined for all $\alpha,\gamma \in (0,1)$,
and, given any $\alpha, \gamma \in (0,1)$,
there exists a constant $C=C[\alpha,\gamma]>0$
such that whenever $f \in C^{0,\alpha,\gamma}_{\GrpK}(\cyl)$,
  \begin{equation}
    \LhatK \RhatK f = f
    \qquad \mbox{and} \qquad
    \norm{\RhatK f}_{C^{2,\alpha,\gamma}(\cyl)}  \leq C \norm{f}_{C^{0,\alpha,\gamma}(\cyl)},
  \end{equation}
recalling \eqref{cyldef}, \eqref{LhatKdef}, \eqref{GrpKdef},
\eqref{catnorm}, and \eqref{catspace}.
\end{prop}

\begin{proof}
Let $f \in C^{0,\alpha,\gamma}_{\GrpK}(\cyl)$ for some $\alpha,\gamma \in (0,1)$.
For each nonnegative integer $n$ we define the functions $f_n^{\pm}: \R \to \R$ by
  \begin{equation}
    f_n^+(t):=\int_0^{2\pi} f(t,\theta) \cos n\theta \, d\theta
     \qquad \mbox{ and } \qquad
    f_n^-(t):=\int_0^{2\pi} f(t,\theta) \sin n\theta \, d\theta,
  \end{equation}
but $f$ is $\GrpK$-even,
so by \eqref{GrpKdef} and \eqref{catspace}
$f_n^-(t) \equiv 0$ for every $n$
and $f_n^+(t) \equiv 0$ for every odd $n$,
so that
$f(t,\theta)=\frac{1}{2\pi}f_0^+(t)+\frac{1}{\pi}\sum_{n=1}^\infty f_{2n}^+(t) \cos 2n\theta$,
at least distributionally.
From the identities
  \begin{equation}
    (\partial_t - \tanh t)(\partial_t + \tanh t) + 1 = \partial_t^2+2\sech^2 t
    \qquad \mbox{and} \qquad
    (\partial_t + \tanh t)(\partial_t - \tanh t) + 1 = \partial_t^2
  \end{equation}
we find that for $n \geq 2$ the kernel (without any restriction on the rate of growth)
of $\partial_t^2+2\sech^2 t-n^2$ is spanned by the functions
$(\partial_t - \tanh t)e^{\pm nt}=(\pm n - \tanh t)e^{\pm nt}$,
while for $n=0$ it is spanned by the functions
$-(\partial_t-\tanh t)1=\tanh t$
and
$(\partial_t-\tanh t)t=1- t \, \tanh t$
(the Jacobi fields on the catenoid \eqref{stdcat}
induced respectively by
vertical translation and dilations about the origin),
and for $n=1$ (not needed for this construction)
the kernel is spanned by the functions
$(\partial_t-\tanh t)\sinh t = \sech t$
and
$(\partial_t-\tanh t)t \sinh t=\sinh t + t \sech t$
(which, multiplied by linear combinations of $\cos \theta$ and $\sin \theta$,
respectively generate horizontal translations and rotations about horizontal axes
through the center of \eqref{stdcat}).

It follows (and is straightforward to check directly)
that if for each nonnegative integer $n \neq 1$ we define
the function $u_n: \R \to \R$ by
  \begin{equation}
  \label{solmodes}
    \begin{aligned}
      u_0(t)
      :=
      &\int_0^t \left[(t-s) \tanh s \tanh t + (\tanh t - \tanh s)\right]f_n(s) \, ds
        \mbox{ and for $n>1$} \\
      u_n(t)
      :=
      &\frac{n+\tanh t}{2n(1-n^2)}e^{-nt} \int_{-\infty}^t (n-\tanh s) e^{ns} f_n(s) \, ds
        +\frac{n-\tanh t}{2n(1-n^2)}e^{nt} \int_t^\infty (n+\tanh s)e^{-ns} f_n(s) \, ds,
    \end{aligned}
  \end{equation}
then $u_n$ solves
$\left(\partial_t^2+2\sech^2 t - n^2\right)u_n=f_n$
with $u_0(0)=\dot{u}_0(0)=0$ and $u_n$ bounded whenever $f_n$
is compactly supported and $n>1$.
Therefore the distribution
  \begin{equation}
    u:=\frac{1}{2\pi}u_0+\frac{1}{\pi}\sum_{n=1}^\infty u_{2n}
  \end{equation}
solves $\LhatK u = f$,
at least in the distributional sense,
and is even (also as a distribution)
under the reflections
$\xbarhat_K$ and $\ybarhat_K$ (defined in \eqref{GrpKdef}).
It is elementary to verify from \eqref{solmodes} that
  \begin{equation}
    \begin{aligned}
      &\abs{u_n(t)}
      \leq
      \frac{C[\gamma]}{n^2+1}\norm{f}_{C^{0,\alpha,\gamma}(\cyl)}e^{\gamma \abs{t}}, \mbox{ so}
      &\norm{u}_{C^{0,0,\gamma}(\cyl)}
      \leq
      C[\gamma]\norm{f}_{C^{0,\alpha,\gamma}(\cyl)}
    \end{aligned}
  \end{equation}
for some constant $C[\gamma]$ independent of the data $f$.
Standard elliptic theory, using in particular the Schauder estimates,
then implies that in fact $u$ is a classical solution satisfying
  \begin{equation}
    \norm{u}_{C^{2,\alpha,\gamma}(\cyl)}
    \leq
    C[\alpha,\gamma]\norm{f}_{C^{0,\alpha,\gamma}(\cyl)}
  \end{equation}
for some constant $C[\alpha,\gamma]>0$ independent of the data $f$,
and we have already observed that $u$ is $\GrpK$-even.
Taking $\RhatK f:=u$ thus concludes the proof.
\end{proof}

\subsection*{Approximate solutions on the toral regions}
Recalling \eqref{That} and \eqref{chihatidef},
note that both $\That_i$ and $\chihat_i$
are, for all $i$, independent of $m$ and the $\zeta, \xi$ parameters.
Recalling also \eqref{Tipre} and \eqref{Lchidef},
from items (i) and (ii) of Proposition \ref{torest}
we see that on $\T_i$
  \begin{equation}
  \label{toroplimit}
    \lim_{m \to \infty} T_i^* \Lchi {T_i^*}^{-1}
    =
    \Delta_{\chihat_i},
  \end{equation}
where the convergence is to be interpreted
along the lines described for \eqref{catoplimit},
using in this case, in addition to Proposition \ref{torest},
the fact that $\lim_{m \to \infty} m\sqrt{\tauj{i}}=0$,
as follows from the second line of \eqref{estimatesforforceestimate}.
Note additionally that,
by \eqref{Phintertwine}, \eqref{Grpdef}, \eqref{initsurfdef}, and \eqref{Tdef},
the pullback by $T_i$
of any $\Grp$-invariant function on $\Sigma$
to $\T_i$ must satisfy periodic boundary conditions
on the rectangular part of the boundary
and must moreover respect a $D_2$ group of symmetries.

Specifically, for $1 \leq i \leq N$, we define the quotients
  \begin{equation}
  \label{quotients}
    \left(\T_i / \mathord \sim\right)
    \subset
    \That_i / \mathord \sim,
    \quad
    \mbox{ where $(\xx,\yy) \sim (\xx',\yy') \Leftrightarrow 
      (\xx-\xx',\yy-\yy') \in \frac{\sqrt{2}\pi}{k}\Z \times \frac{\sqrt{2}\pi}{\ell} \Z$},
  \end{equation}
so that $\T_i / \mathord \sim$ (or $\That_i / \mathord \sim$) is a $\sqrt{2}\pi/k \times \sqrt{2}\pi/\ell$ torus
with one disc (or point) deleted if $i \in \{1,N\}$ and two otherwise. 
We also define the $D_2$ subgroup of symmetries of $\That_i / \mathord \sim$
  \begin{equation}
  \label{GrpTdef}
    \begin{aligned}
      &\GrpT{i}
      :=
      \left\{\widehat{I}_T,\xbarhat_{T_i},\ybarhat_{T_i},\xbarhat_{T_i}\ybarhat_{T_i}\right\},
        \qquad \mbox{where }
        \widehat{I}_T(\xx,\yy):=(\xx,\yy), \\
      &\xbarhat_{T_i}(\xx,\yy)
        :=
        \begin{cases}
          (-\xx,\yy) \mbox{ if $i \in \{1,N\}$} \\
          \left(\frac{\pi}{\sqrt{2}k}-\xx,\yy \right) \mbox{ if $1<i<N$},
        \end{cases}
      \quad \mbox{and} \qquad
      \ybarhat_{T_i}(\xx,\yy)
        :=
        \begin{cases}
          (\xx,-\yy) \mbox{ if $i \in \{1,N\}$} \\
          \left(\xx,\frac{\pi}{\sqrt{2}\ell}-\yy\right) \mbox{ if $1<i<N$}
        \end{cases}
    \end{aligned}
  \end{equation}
(using coordinates on the universal cover of $\That_i / \mathord \sim$
to define the symmetries).
Obviously $\GrpT{i}$ preserves both $\That_i / \mathord \sim$ and $\T_i / \mathord \sim$.

\begin{remark}
\label{symmetricextension}
Recalling Notation \ref{functionspaceappendages},
it now follows from
\eqref{Phintertwine}, \eqref{Grpdef}, \eqref{initsurfdef},
\eqref{catdef}, \eqref{catrdef}, \eqref{Tdef}, \eqref{torrdef},
\eqref{GrpKdef}, and \eqref{GrpTdef}
that, for any $\alpha \in [0,1)$ and nonnegative integer $j$,
a function
$f \in C^{j,\alpha}_{loc}
 \left(\bigcup_{i=1}^{N-1} \catr[i] \cup \bigcup_{i=1}^N \torr[i]\right)$
 extends (uniquely) to a function in $C^{j,\alpha}_{loc,\Grp}(\Sigma)$
 if and only if
 ${\cat_i}^*f \in C^{j,\alpha}_{loc,\GrpK}\left(\cyl_{a_i}\right)$
for each $1 \leq i \leq N-1$
and $T_i^*f$ descends to a function in
$C^{j,\alpha}_{loc,\GrpT{i}}\left(\T_i / \mathord \sim \right)$
for each  $1 \leq i \leq N$.
\end{remark}

Motivated also by \eqref{rhorhohat} and definition \eqref{globalweightednorm},
we are led to define,
for any nonnegative integer $j$, $\alpha \in [0,1)$, and $\gamma \in (0,\infty)$,
the norms
  \begin{equation}
  \label{tornorm}
    \norm{\cdot}_{C^{j,\alpha,\gamma}\left(\That_i / \mathord \sim\right)}
    :=
    \norm{\cdot : C^{j,\alpha}\left(\That_i / \mathord \sim,\chihat_i,{\rhohat_i}^{-\gamma}\right)}
  \end{equation}
(recalling \eqref{weightedhoelderdef})
and the corresponding Banach spaces of $\GrpT{i}$-even functions
  \begin{equation}
  \label{torspace}
    C^{j,\alpha,\gamma}_{\GrpT{i}}\left(\That_i / \mathord \sim\right)
    :=
    \left. \left\{
      u \in C^{j,\alpha}\left(\That_i / \mathord \sim,\chihat_i\right) \; \right| \;
      \norm{u}_{C^{j,\alpha,\gamma}\left(\That_i / \mathord \sim\right)}<\infty
      \mbox{ and }
      \mathfrak{g}.u=u \circ \mathfrak{g}
      \mbox{ for all } \mathfrak{g} \in \GrpT{i}
    \right\}.
  \end{equation}

Clearly
$\Delta_{\chihat_i}: C^{2,\alpha,\gamma}_{\GrpT{i}}\left(\That_i / \mathord \sim\right)
 \to C^{0,\alpha,\gamma}_{\GrpT{i}}\left(\That_i / \mathord \sim\right)$
is bounded (independently of $\alpha,\gamma \in (0,1)$).
Proposition \ref{torsol} below presents a suitable inverse,
modulo the extended substitute kernel
and with a support condition on the source function.
The support assumption
(expressed below
using the function $d_i$ from \eqref{chihatidef}
and writing $\spt f$ for the support of a function $f$)
we can afford to make
because in practice we will apply Proposition \ref{catsol}
before applying Proposition \ref{torsol}.
The extended substitute kernel
we formally define right now.
We first recall the definition \eqref{wbar} of $\wbar_i$,
remembering in particular that $\wbar_1$ and $\wbar_N$ vanish identically,
and for $1 \leq i \leq N$ we introduce
$w_i \in C^\infty_{\Grp}(\Sigma)$
defined by
  \begin{equation}
  \label{w}
    T_i^*w_i
    :=
    \cutoff{\frac{1}{10\ell}}{\frac{1}{5\ell}} \circ d_i
    \qquad \mbox{and} \qquad
    T_i|_{\Sigma \backslash \Grp \torr[i]}
    :=
    0
    \qquad \mbox{for} \qquad
    1 \leq i \leq N,
  \end{equation}
recalling the function $d_i$ from \eqref{chihatidef}.
Finally we define the \emph{extended substitute kernel}
to be the linear span in $C^\infty(\Sigma)$
of $\{w_i,\wbar_i\}_{i=1}^N$.

\begin{prop}[Solutions to the model problems on the torus]
\label{torsol}
Let $\ell \geq k \in \Z \cap [1,\infty)$ and $i \in \Z \cap [1,N]$.
There exists a linear map
  \begin{equation}
    \RhatT{i}=\RhatT{i}[k,\ell]:
    \left.\left\{f \in C^{0,\alpha}_{\GrpT{i}}\left(\That_i / \mathord \sim,\chihat_i \right)
      \; \right| \;
      \spt f \subset \left\{d_i > \frac{1}{20\ell}\right\} \right\}
    \to
    C^{2,\alpha,2}_{\GrpT{i}}\left(\That_i / \mathord \sim \right) \times \R \times \R
  \end{equation}
defined for all $\alpha \in (0,1)$,
and, given any $\alpha \in (0,1)$,
there exists a constant $C=C[k,\ell,\alpha]>0$
such that whenever $f$ belongs to the domain of $\RhatT{i}$ above
and $(u,\mu,\mubar)=\RhatT{i}f$, then
  \begin{equation}
    \Delta_{\chihat_i}u = f + \mu T_i^*w_i + \mubar T_i^*\wbar_i
    \qquad \mbox{and} \qquad
    \norm{u}_{C^{2,\alpha,2}\left(\That_i / \mathord \sim\right)}
      +\abs{\mu}+\abs{\mubar}
    \leq
    C \norm{f: C^{0,\alpha}\left(\That_i / \mathord \sim,\chihat_i\right)},
  \end{equation}
recalling \eqref{That}, \eqref{chihatidef}, \eqref{wbar}, \eqref{quotients},
\eqref{GrpTdef}, \eqref{tornorm}, \eqref{torspace}, and \eqref{w}.
\end{prop}

\begin{proof}
Suppose $f \in C^{0,\alpha}_{\GrpT{i}}(\That_i / \mathord \sim, \chihat_i)$
has support contained in the set $U:=\{1/d_i < 20 \ell\}$.
We intend to apply a conformal change of metric
and attack the corresponding problem on the flat torus
$(\Tbar_i,\geuc)$,
where $\Tbar_i$ is simply $\That_i / \mathord \sim$
with the missing point(s) filled in
and $\geuc=d\xx^2+d\yy^2$ is the standard flat metric.
By definition \eqref{chihatidef} (and because $\That_i$ is two-dimensional)
the equation $\Delta_{\chihat_i}u=f$ on $\That_i / \mathord \sim$ is equivalent to
$\Delta_{\geuc}u=\rhohat_i^2 f$.
Note that $\GrpT{i}$ acts by isometries on $(\Tbar_i,\geuc)$
in the obvious way.
Clearly the function $\rhohat_i$ defined in \eqref{chihatidef}
descends to a function (which we give the same name)
in $C_{\GrpT{i}}^\infty(\That_i / \mathord \sim)$
and clearly
  \begin{equation}
    \norm{\rhohat_i|_U: C^j(U,\chihat_i)}
    \leq C[\ell,j]
  \end{equation}
for some constant $C[\ell,j]$.
Consequently
  \begin{equation}
  \label{rhosquaredfest}
    \norm{{\rhohat_i}^2 f : C^{0,\alpha}(\Tbar_i,\geuc)}
    \leq
    C[\ell]\norm{f: C^{0,\alpha}(\That_i / \mathord \sim, \chihat_i)}
  \end{equation}
(where we have trivially extended $\rhohat_i^2f$ to a function of the same name
on $\Tbar_i$, simply by requiring it to vanish at the filled in point(s)).

Of course the equation
$\Delta_{\geuc}u=\rhohat_i^2 f$
has a solution on $\Tbar_i$
if and only if the right-hand side has vanishing integral over $\Tbar_i$,
which we do not assume.
Accordingly we would like to permit ourselves the freedom
of adding a constant to the right-hand side.
Soon though (at the end of this section)
we intend to transfer the solution from this model problem
to the initial surface,
so we want to confine any modification of the right-hand side
to the toral region in question,
avoiding any interference on the adjoining catenoidal regions.
Therefore we will we use the cutoff function $T_i^*w_i$ in \eqref{w}
in lieu of the constant function $1$ for the purpose
of adjusting the right-hand side to make it orthogonal to the kernel.

More precisely we note that $T_i^*w_i$ descends smoothly to $\Tbar_i$
and we define
 \begin{equation}
 \label{mudef}
    \mu
    :=
    -\frac{\int_{\Tbar_i} {\rhohat_i}^2 f \, d\xx \, d\yy}
      {\int_{\Tbar_i} \rhohat_i^2 T_i^*w_i \, d\xx \, d\yy},
  \end{equation}
so that $\int_{\Tbar_i} \rhohat_i^2\left(f+\mu T_i^*w_i\right) \, d\xx \, d\yy = 0$.
Consequently there is a unique function $u_0: \Tbar_i \to \R$ solving
  \begin{equation}
  \label{sol0}
    \Delta_{\geuc}u_0=\rhohat_i^2(f+\mu T_i^*w_i)
  \end{equation}
and satisfying $\int_{\Tbar_i} u_0 \, d\xx \, d\yy = 0$;
in particular $u_0$ is necessarily $\GrpT{i}$-invariant.
Note also that by \eqref{chihatidef}
$\rhohat_i \geq 1$ on $\That_i$
and that by \eqref{That} and \eqref{w}
$w=1$ on a region of positive $\geuc$-area
(depending on just $k$ and $\ell$),
while of course $\Tbar_i$ itself has area $\frac{2\pi^2}{k\ell}$;
it then follows from \eqref{rhosquaredfest} and \eqref{mudef}
that
  \begin{equation}
  \label{muest}
    \abs{\mu}
    \leq
    C[k,\ell]\norm{f: C^{0,\alpha}(\That_i/ \mathord \sim,\chihat_i)}
  \end{equation}
for some constant $C[k,\ell]>0$.
The classical global Schauder estimates
applied to \eqref{sol0}
imply in particular that
  \begin{equation}
  \label{sol0est}
    \norm{u_0: C^0(\Tbar_i)}
    \leq
    C[k,\ell]\norm{\rhohat_i^2(f+\mu T_i^*w_i): C^{0,\alpha}(\Tbar_i,\geuc)}
    \leq
    C[k,\ell]\norm{f: C^{0,\alpha}(\That_i/ \mathord \sim,\chihat_i)}
  \end{equation}
(for a possibly larger constant $C[k,\ell]$ than above),
where for the first inequality we have again used the fact
that $\Tbar_i$ is just a flat $\sqrt{2}\pi/k \times \sqrt{2}\pi/\ell$ torus
and for the second we have used \eqref{rhosquaredfest} and \eqref{muest}.

We still need to arrange the rapid decay
of our solution toward the point(s)
in $\Tbar_i$ missing from $\That_i$.
To this end we first observe that,
because both $f$ (by assumption)
and $w$ (by definition \eqref{w})
have support contained in $U=\{1/d_i < 20 \ell\}$,
the solution $u_0$ to \eqref{sol0}
is harmonic on the set
$\Tbar_i \backslash U=\left\{d_i \leq \frac{1}{20 \ell}\right\}$,
where, as we see from \eqref{chihatidef},
$d_i=1/\rhohat_i$.
For $i \in \{1,N\}$
this set has one component---the closed disc of radius
$\frac{1}{20\ell}$ and center $p_0:=(0,0)$---while
otherwise it has two components---the closed discs of radius
$\frac{1}{20\ell}$ and centers
$p_{\pm}:=\left(\frac{\pi}{2\sqrt{2}k},\frac{\pi}{2\sqrt{2}\ell}\right)
 =\pm \left(\frac{mX}{2},\frac{mY}{2}\right)$,
recalling \eqref{RXYdef}.
Now we define $\mubar \in \R$ by
  \begin{equation}
      (-1)^{N-i}\mubar
      :=
      \begin{cases}
        0 \quad \mbox{if} \quad i \in \{1,N\} \\
        \frac{1}{2}[u_0(p_-)-u_0(p_+)] \quad \mbox{if} \quad i \not \in \{1,N\},
      \end{cases}
  \end{equation}
so that by \eqref{sol0est}
  \begin{equation}
  \label{mubarest}
    \abs{\mubar}
    \leq
    C[k,\ell]\norm{f: C^{0,\alpha}(\That_i/ \mathord \sim,\chihat_i)}.
  \end{equation}
Recalling \eqref{vdef}
and noting that $v$ descends smoothly to $\Tbar_i$,
we also define $u: \Tbar_i \to \R$ by
  \begin{equation}
      u
      :=
      \begin{cases}
        u_0-u_0(p_0) \quad \mbox{if} \quad i \in \{1,N\} \\
        u_0-\frac{1}{2}[u_0(p_-)+u_0(p_+)]+(-1)^{N-i}\mubar v_i
         \quad \mbox{if} \quad i \not \in \{1,N\}
      \end{cases}
  \end{equation}
(where we include the alternating signs
because of the one present in definition \eqref{wbar} of $\wbar_i$,
which in turn we included to account for
the alternating direction of the unit normal on the toral regions).
Thus by \eqref{sol0} and \eqref{wbar}
  \begin{equation}
  \label{sol}
    \Delta_{\chihat_i}u_0=f+\mu T_i^*w_i + \mubar T_i^*\wbar_i.
  \end{equation}

Note that $v_i$ is constant on each component of $\Tbar_i \backslash U$,
so, like $u_0$, the function $u$ is harmonic on $\Tbar_i \backslash U$.
By classical harmonic function theory
  \begin{equation}
  \label{harmonicest}
    \norm{u: C^2(\Tbar_i \backslash U,\geuc)}
    \leq
    C[\ell]\norm{u: C^0(\partial U, \geuc)}
  \end{equation}
for some constant $C[\ell]>0$.
On the other hand,
since $v_i(p_\pm)=\pm 1$,
we have $u(p_0)=0$ if $i \in \{1,N\}$
and $u(p_\pm)=0$ otherwise.
Moreover, $u$ is $\GrpT{i}$-invariant
(because $u_0$, $v$, and the constants are),
so, recalling \eqref{GrpTdef},
both first parital derivatives of $u$
also vanish at $p_0$ if $i \in \{1,N\}$
and at both points $p_\pm$ otherwise.
Using Taylor's theorem and \eqref{harmonicest},
we therefore obtain
  \begin{equation}
  \label{tordecayest}
    \norm{u: C^0\left(\Tbar_i \backslash U,\geuc,{d_i}^2\right)}
    \leq
    C[\ell]\norm{u: C^0(\partial U, \geuc)}
    \leq
    C[k,\ell]\norm{f: C^{0,\alpha}(\That_i/ \mathord \sim,\chihat_i)}
  \end{equation}
where we recall \eqref{weightedhoelderdef}
and for the last inequality we use \eqref{sol0est}.
As already observed $\rhohat_i^{-1}=d_i$ on $\Tbar_i \backslash U$,
while on $U$ it is bounded below by $\frac{1}{20\ell}$,
so it now follows from \eqref{weightedhoelderdef}, \eqref{tornorm}
\eqref{sol}, \eqref{tordecayest}, and the standard local Schauder estimates
together with the bounded geometry of $(\That_i / \mathord \sim,\chihat_i)$ that
  \begin{equation}
    \norm{u}_{C^{2,\alpha,2}(\That_i / \mathord \sim)}
    \leq
    C[k,\ell,\alpha]\norm{f: C^{0,\alpha}(\That_i/ \mathord \sim,\chihat_i)},
  \end{equation}
which, along with \eqref{muest}, \eqref{mubarest},
\eqref{sol}, and the already observed $\GrpT{i}$-invariance of $u$,
concludes the proof.
\end{proof}

\subsection*{Exact global solutions}
Now we use Propositions \ref{catsol} and \ref{torsol}
to construct global solutions to the linearized problem
on each initial surface, modulo extended substitute kernel.

\begin{prop}[Global solutions to the linearized problem]
\label{globalsol}
Given a real number $c>0$
and integers
$N \geq 2$ and $\ell \geq k \geq 1$,
there exists $m_0=m_0[N,k,\ell,c]>0$
such that whenever $\zeta,\xi \in [-c,c]^{N-1}$ and $m>m_0$,
there is a linear map
  \begin{equation}
    \Rcal=\Rcal[N,k,\ell,m,\zeta,\xi]
    :C^{0,\alpha,\gamma}_{\Grp[k,\ell,m]}(\Sigma[N,k,\ell,m,\zeta,\xi])
    \to
    C^{2,\alpha,\gamma}_{\Grp[k,\ell,m]}(\Sigma[N,k,\ell,m,\zeta,\xi])
      \times \R^N \times \R^{N-2}
  \end{equation}
(recalling \eqref{globalweightedspace})
defined for all $\alpha,\gamma \in (0,1)$,
and, given $\alpha, \gamma \in (0,1)$,
there is a constant $C=C[N,k,\ell,\alpha,\gamma]>0$
such that whenever
$f \in C^{0,\alpha,\gamma}_{\Grp}(\Sigma)$
and
$(u,(\mu_1,\cdots,\mu_N),(\mubar_2,\cdots,\mubar_{N-1}))=\Rcal f$,
then
  \begin{equation}
    \Lchi u = f + \sum_{i=1}^N \mu_i w_i + \sum_{i=2}^{N-1} \mubar_i \wbar_i
    \qquad \mbox{and} \qquad
    \norm{u}_{2,\alpha,\gamma}
      +\sum_{i=1}^N \abs{\mu_i} + \sum_{i=2}^{N-1} \abs{\mubar_i}
    \leq
    \norm{f}_{0,\alpha,\gamma}
  \end{equation}
(recalling \eqref{globalweightednorm} and \eqref{Lchidef});
moreover, for any fixed $N$, $k$, $\ell$, and $m$, the map
  \begin{equation}
    (f,\zeta,\xi) \mapsto \Rcal[N,k,\ell,m,\zeta,\xi]f
    \mbox{ is continuous}
  \end{equation}
in the sense of Definition \ref{functionalcty}.
\end{prop}

\begin{proof}
Let $c>0$, $\alpha, \gamma \in (0,1)$,
$N \in \Z \cap [2,\infty)$, $\ell \geq k \in \Z \cap [1,\infty)$,
$\zeta,\xi \in [-c,c]^{N-1}$,
and $m \in \Z \cap [m_0,\infty)$,
where $m_0$ is at least as large as the maximum
of the homonymous quantities appearing in Propositions
\ref{catest}, \ref{torest}, \ref{catsol}, and \ref{torsol}
and is subject to an additional lower bound imposed at the end of the proof.
Recalling \eqref{catrdef}, we start by defining,
for $1 \leq i \leq N-1$,
the linear maps
  \begin{equation}
    \Psi_{\catr[i]}: C(\catr[i]) \to C(\Sigma)
  \end{equation}
so that $\Psi_{\catr[i]}f$ is the unique $\Grp$-equivariant function which vanishes
outside $\Grp \catr[i]$ and which satisfies
  \begin{equation}
    \left.\left(\Psi_{\catr[i]}f\right)\right|_{\catr[i]}
      :=
      {\cat_i^{-1}}^*
      \left[
        \left(\cutoff{a_i}{a_i-1/2} \circ \abs{t}\right) \cdot {\cat_i}^*f
      \right]
      =
      \left(\cutoff{a_i}{a_i-1/2} \circ \abs{t \circ \cat_i^{-1}}\right) \cdot f,
  \end{equation}
recalling \eqref{cutoff}, \eqref{adef}, and \eqref{catdef}.
Note that by \eqref{chihatK} the $j\textsuperscript{th}$ $\chihat_K$ covariant derivative of
$\cutoff{a_i}{a_i-1} \circ \abs{t}$ is uniformly $\chihat_K$-bounded
on $\cyl_{a_i}$ by a constant depending on just $j$.
Recalling Proposition \ref{catsol}, we also define
  \begin{equation}
  \label{RtildeK}
    \begin{aligned}
      \RtildeK: C^{0,\alpha,\gamma}_{\Grp}(\Sigma) &\to C^{2,\alpha,\gamma}_{\Grp}(\Sigma)
        \mbox{ by } \\
      \RtildeK f
      :=
      \sum_{i=1}^{N-1} \Psi_{\catr[i]}v_{\catr[i]},
      \qquad &\mbox{with} \qquad
      v_{\catr[i]}
      :=
      {\cat_i^{-1}}^*\left(\RhatK \left({\cat_i}^*\Psi_{\catr[i]}(f|_{\catr[i]})\right)\right)
    \end{aligned}
  \end{equation}
where ${\cat_i}^*\Psi_{\catr[i]}(f|_{\catr[i]})$ is trivially (and smoothly) extended
from $\cyl_{a_i}$ to $\cyl$ so that it vanishes outside $\cyl_{a_i}$,
recalling \eqref{cyldef} and \eqref{cyladef}.
Then
  \begin{equation}
  \label{caterror}
    \Lchi \RtildeK f
    =
    \sum_{i=1}^{N-1}
      \left(
        \left[\Lchi,\Psi_{\catr[i]}\right]v_{\catr[i]}
        +\Psi_{\catr[i]}{\cat_i^{-1}}^*
          \left( {\cat_i}^*\Lchi{\cat_i^{-1}}^*-\LhatK \right){\cat_i}^*v_{\catr[i]}
        +\Psi_{\catr[i]}^2 f|_{\catr[i]}
      \right),
  \end{equation}
where in the first term the brackets indicate the commutator of the operators they enclose,
in the second term we recall \eqref{LhatKdef},
and in the last term we make use of Proposition \ref{catsol}.

We will absorb the ``cutoff error'' in \eqref{caterror},
present in the first and third terms,
into the right-hand side when solving on the toral regions
in the next step.
More precisely, for any given $f \in C^{0,\alpha,\gamma}_{\Grp}(\Sigma)$
we define 
  \begin{equation}
  \label{fT}
    f_T
    :=
    f-\sum_{i=1}^{N-1}\Psi_{\catr[i]}^2 f|_{\catr[i]}
      -\sum_{i=1}^{N-1}\left[\Lchi,\Psi_{\catr[i]}\right]v_{\catr[i]},
  \end{equation}
where each $v_{\catr[i]}$ is defined (for the given $f$) in \eqref{RtildeK}.
Note that $f_T$ is $\Grp$-equivariant and has support contained in
$\Grp \left(\bigcup_{i=1}^N \torr[i]\right)$.
In fact, since
  \begin{equation}
    \tauj{i} \cosh \left(a_i-\frac{1}{2}\right)
    =
    \tauj{i} \cosh a_i \left(\cosh \frac{1}{2} - \tanh a_i \, \sinh \frac{1}{2}\right)
    \geq
    e^{-1/2} \tauj{i} \cosh a_i
    >
    \frac{1}{20\ell m}
  \end{equation}
(using \eqref{adef} for the final inequality),
we have, recalling \eqref{Tdef},
  \begin{equation}
  \label{fTspt}
    \spt {T_i}^*f_T \subset \left\{d_i > \frac{1}{20\ell}\right\}.
  \end{equation}

Next, recalling \eqref{torrdef}, for $1 \leq i \leq N$ we now define the linear maps
  \begin{equation}
    \Psi_{\torr[i]}: C(\torr[i]) \to C(\Sigma)
  \end{equation}
so that $\Psi_{\torr[i]}f$ is the unique $\Grp$-equivariant function on $\Sigma$
vanishing outside $\Grp \torr[i]$ and satisfying
  \begin{equation}
    \begin{aligned}
      \left.\left(\Psi_{\torr[i]}f\right)\right|_{\torr[i]}
      &:=
      {T_i^{-1}}^*
        \left[
          \left(\cutoff{\ln m\tauj{1}^{1/3}}{\ln m\tauj{1}^{1/6}} \circ \ln d_i\right)
            \cdot {T_i}^*f
        \right] \\
      &=
      \left(
        \cutoff{\ln m\tauj{1}^{1/3}}{\ln m\tauj{1}^{1/6}}
          \circ \ln \left(d_i \circ {T_i}^{-1}\right)
      \right)
      \cdot f
    \end{aligned}
  \end{equation}
for which we recall \eqref{Tdef} and \eqref{chihatidef}.
Note that by \eqref{cisummary} and \eqref{taudef}
we have $\sqrt{\tauj{i}}<\tauj{1}^{1/3}<\tauj{1}^{1/6}$
for $1 \leq i \leq N-1$
and that moreover by \eqref{chihatidef}
all $\chihat_i$ covariant derivatives
of $\cutoff{\ln m\tauj{1}^{1/3}}{\ln m\tauj{1}^{1/6}} \circ \ln d_i$
are uniformly $\chihat_i$-bounded on $\T_i$.
Now, recalling Proposition \ref{torsol} and \eqref{fT} just above, we also define
for $1 \leq i \leq N$ the maps
  \begin{equation}
  \label{RtildeT}
    \begin{aligned}
      \RtildeT{i}: C^{0,\alpha,\gamma}_{\Grp}(\Sigma)
      &\to
      C^{2,\alpha,\gamma}_{\Grp}(\Sigma) \times \R \times \R \mbox{ by} \\
      \RtildeT{i}f:=\left(\Psi_{\torr[i]} v_{\torr[i]},\mu,\mubar\right)
      \quad \mbox{with} \quad
      v_{\torr[i]}
      &:=
      {T_i^{-1}}^*\widehat{v}_{\torr[i]}
      \quad \mbox{and} \quad
      \left(\widehat{v}_{\torr[i]},\mu,\mubar\right)
      :=
      \RhatT{i}\left({T_i}^* f_T\right).
    \end{aligned}
  \end{equation}
Here we are implicitly regarding ${T_i}^*f_T$
as a function on $\That_i / \mathord \sim$ (possible because $f_T$ is $\Grp$-equivariant)
after extending it to a function on $\That_i$
which simply vanishes outside $\T_i$,
and moreover we see that \eqref{fTspt}
ensures that this function truly belongs to the domain of $\RhatT{i}$.
It now follows by Proposition \ref{torsol}
that if $(u_{\torr[i]},\mu,\mubar)=\RtildeT{i}f$, then
  \begin{equation}
  \label{torerror}
    \Lchi u_{\torr[i]}
    =
    \left[\Lchi,\Psi_{\torr[i]}\right]v_{\torr[i]}
      +\Psi_{\torr[i]}{T_i^{-1}}^*
        \left(
          {T_i}^*\Lchi{T_i^{-1}}^* - \Delta_{\chihat_i}
        \right){T_i}^*v_{\torr[i]}
      +\Psi_{\torr[i]}f_T|_{\torr[i]}
      +\mu \, w_i + \mubar \, \wbar_i.
  \end{equation}

Next we define the approximate solution operator
  \begin{equation}
    \begin{aligned}
      &\Rtilde: C^{0,\alpha,\gamma}_{\Grp}(\Sigma)
      \to
      C^{2,\alpha,\gamma}_{\Grp}(\Sigma) \times \R^N \times \R^{N-2} \mbox{ by} \\
      &\Rtilde f
      :=
      \left(
        \RtildeK f + \sum_{i=1}^N u_{\torr[i]},
        \left(\mu_1,\ldots,\mu_N\right),
        \left(\mubar_2,\ldots,\mubar_{N-1}\right)
      \right) \mbox{ with } 
      \left(u_{\torr[i]},\mu_i,\mubar_i\right)
      :=
      \RtildeT{i}f,
    \end{aligned}
  \end{equation}
where $\RtildeK$ and $\RtildeT{i}$
are defined in \eqref{RtildeK} and \eqref{RtildeT} above
and where from the output of $\Rtilde$ we are simply omitting
$\mubar_1=\mubar_N=0$, as indicated.
Clearly $\Rtilde$ (from its definition and Propositions \ref{catsol} and \ref{torsol})
is bounded independently of $c$ and $m$.
Moreover, the map $(f,\zeta,\xi) \mapsto \Rtilde[N,k,\ell,m,\zeta,\xi]f$
is manifestly continuous (in the sense of Definition \ref{functionalcty}),
since all the operators (including $\Lchi$ itself) on $\Sigma$
used to construct it clearly enjoy this continuous dependence themselves,
while the maps $\RhatK$ and $\RhatT{i}$ are of course independent of the parameters.
Defining also the operator
  \begin{equation}
    \begin{aligned}
      \Ltilde: C^{2,\alpha,\gamma}_{\Grp}(\Sigma) \times \R^N \times \R^{N-2}
        &\to C^{0,\alpha,\gamma}_{\Grp}(\Sigma) \mbox{ by} \\
      \Ltilde \left(u,\left(\mu_1,\ldots,\mu_N\right),
        \left(\mubar_2,\ldots,\mubar_{N-1}\right)\right)
      &:=
      \Lchi u - \sum_{i=1}^N \mu_i w_i - \sum_{i=2}^{N-1} \mubar_i \wbar_i
    \end{aligned}
  \end{equation}
and using \eqref{caterror}, \eqref{fT}, \eqref{torerror},
and the definitions of $\Psi_{\catr[i]}$ and $\Psi_{\torr[i]}$ above,
we find that for any $f \in C_{\Grp}^{0,\alpha,\gamma}(\Sigma)$
  \begin{equation}
  \label{LR}
    \begin{aligned}
      \Ltilde \Rtilde f-f
      =
      &\sum_{i=1}^{N-1}
          {\cat_i^{-1}}^*
            \left({\cat_i}^*\Lchi{\cat_i^{-1}}^*-\LhatK \right){\cat_i}^*v_{\catr[i]}
        +\sum_{i=1}^N
          {T_i^{-1}}^*
            \left({T_i}^*\Lchi{T_i^{-1}}^* - \Delta_{\chihat_i}\right){T_i}^*v_{\torr[i]} \\
      &+\sum_{i=1}^N \left[\Lchi,\Psi_{\torr[i]}\right]v_{\torr[i]}
    \end{aligned}
  \end{equation}
where $v_{\catr[i]}$ and $v_{\torr[i]}$
are defined in \eqref{RtildeK} and \eqref{RtildeT}.

From \eqref{Lchidef}, items (i) and (v) of Proposition \ref{catest},
items (i) and (iv) of Proposition \ref{torest},
Propositions \ref{catsol} and \ref{torsol},
and the definitions of $\Psi_{\catr[i]}$ and $\Psi_{\torr[i]}$ above
we find that the first two sums in \eqref{LR}
have $C^{0,\alpha,\gamma}$ norm bounded by
$m^{-2}$ times some constant
$C=C[N,k,\ell,\alpha,\gamma]>0$
times $\norm{f}_{0,\alpha,\gamma}$.
As for the commutator terms,
note that each commutator
$\left[\Lchi,\Psi_{\torr[i]}\right]$ itself
has support contained in
$T_i\left(\left\{d_i \leq m\tau^{1/6}\right\}\right)$,
but by Proposition \ref{torsol}
and the definition of $v_{\torr[i]}$
in \eqref{RtildeT}
we know that
  \begin{equation}
    \norm{v_{\torr[i]}: C^{0,\alpha}
      \left(T_i\left(\left\{d_i \leq m\tauj{1}^{1/6}\right\}\right),
      \chi,\frac{m^\gamma}{\rho^\gamma}\right)}
    \leq
    C\norm{f}_{0,\alpha,\gamma}m^{2-\gamma}\tauj{1}^{1/3-\gamma/6}
  \end{equation}
(for a possibly larger $C=C[N,k,\ell,\alpha,\gamma]$ than above).
Thus (making use of line 2 of \eqref{estimatesforforceestimate})
we may take $m$ large enough (in terms of $C$)
so that $\Ltilde \Rtilde$ is a small perturbation
of the identity operator on $C^{0,\alpha,\gamma}_{\Grp}(\Sigma)$
and consequently invertible.
Taking $\Rcal:=\Rtilde \left(\Ltilde \Rtilde \right)^{-1}$
concludes the proof.
\end{proof}

As an immediate application we obtain the first-order correction
of the initial surface toward minimality.

\begin{cor}[The solution to first order]
\label{firstorder}
Given $c>0$, $\alpha,\gamma \in (0,1)$,
and integers $N \geq 2$ and $\ell \geq k \geq 1$,
there exist real numbers $m_0=m_0[N,k,\ell,c,\gamma]>0$
and $C=C[N,k,\ell]>0$ such that whenever
$\zeta,\xi \in [-c,c]^{N-1}$, $m>m_0$, and
  \begin{equation}
    \left(u_1,\left(\lambda_1,\ldots,\lambda_N\right),
      \left(\lambdabar_2,\ldots,\lambdabar_{N-1}\right)\right)
    :=
    -\Rcal \left(\rho^{-2}H-\sum_{i=2}^{N-1} \Dcal_i \wbar_i\right)
  \end{equation}
(recalling \eqref{Ddef}, \eqref{AHdef}, \eqref{rhodef}, \eqref{wbar}, and Proposition \ref{globalsol}),
then
  \begin{equation}
    \norm{u_1}_{2,\alpha,\gamma}
      + \sum_{i=1}^N \abs{\lambda_i} + \sum_{i=2}^{N-1} \abs{\lambdabar_i}
    \leq
    C\tauj{1}
  \end{equation}
(recalling \eqref{globalweightednorm});
moreover, $\lambda_1,\ldots,\lambda_N$
and $\lambdabar_2,\ldots,\lambdabar_{N-2}$
all depend continuously on $(\zeta,\xi)$,
as does $u_1$ (in the sense of Definition \ref{functionalcty}).
\end{cor}

\begin{proof}
All the claims follow directly from
Corollary \ref{Hest} and Proposition \ref{globalsol},
with the obvious supplemental facts that
$\Dcal_i$ is continuous in the parameters
and, in the sense of Definition \ref{functionalcty},
$\wbar_i$ is too.
\end{proof}

\section{The main theorem}
\subsection*{The nonlinear terms}
Recall \eqref{AHdef}, \eqref{Hcaldef} and \eqref{Lcaldef}.
We will need the following estimate for the nonlinear contribution
  \begin{equation}
  \label{Qcal}
    \Qcal[u]=\Qcal[N,k,\ell,m,\zeta,\xi,u]
    :=\Hcal[N,k,\ell,m,\zeta,\xi,u]-H-\Lcal[N,k,\ell,m,\zeta,\xi]u
  \end{equation}
that the deforming function $u$ makes to the mean curvature.
(Of course $H=\Hcal[N,k,\ell,m,\zeta,\xi,0]$.)

\begin{lemma}[The nonlinear terms]
\label{quadest}
Given $C_u,c>0$, $\alpha,\gamma \in (0,1)$,
and integers $N \geq 2$ and $\ell \geq k \geq 1$,
there exists $m_0:=m_0[N,k,\ell,m,C_u,c]>0$
such that
(recalling \eqref{taubardef}, \eqref{rhodef}, \eqref{globalweightednorm}, and \eqref{Qcal})
$\Qcal[N,k,\ell,m,\zeta,\xi,u]$ is well-defined and
  \begin{equation}
    \norm{\rho^{-2}\Qcal[N,k,\ell,m,\zeta,\xi,u]}_{0,\alpha,\gamma}
    \leq
    \taubarj{1}^{1+\gamma/2}
  \end{equation}
whenever $m>m_0$, $\zeta,\xi \in [-c,c]^{N-1}$,
and $u \in C^{2,\alpha}(\Sigma,\chi)$ satisfies
$\norm{u}_{2,\alpha,\gamma} \leq C_u\tauj{1}$;
furthermore, for each fixed $N$, $k$, $\ell$, and $m>m_0$,
the map $(u,\zeta,\xi) \mapsto \Qcal[N,k,\ell,m,\zeta,\xi,u]$
is continuous (in the sense of Definition \ref{functionalcty}).
\end{lemma}

\begin{proof}
That $\Qcal[u]$ is defined at all will be clear
from Lemma \ref{perimm} in conjunction with the estimates below
(which show that $\iota[u]$ \eqref{deformediota} is an immersion and $\Hcal[u]$ is defined).
The continuity follows immediately from
the smooth dependence (Remark \ref{initsurfcty})
of the initial surfaces on the parameters
and from definition \eqref{Hcaldef}.
To make the estimate
we will apply Lemma \ref{perimm}
to the embedding $\iota: \Sigma \to \Sph^3$
of the initial surface into $(\Sph^3,\gsph)$,
as $\Qcal[u]$ can then be read off from item (iv) of the lemma.
First we observe,
recalling \eqref{weightedhoelderdef}, \eqref{gdef}, 
\eqref{AHdef}, and \eqref{rhodef},
that by \eqref{chidef}
and Propositions \ref{catest} and \ref{torest}
  \begin{equation}
  \label{initests}
    \norm{g: C^j\left({T^*\Sigma}^{\otimes 2},\chi,\rho^{-2}\right)}
      +\norm{g^{-1}: C^j\left({T\Sigma}^{\otimes 2},\chi, \rho^2\right)}
      +\norm{A: C^j\left({T^*\Sigma}^{\otimes 2},\chi,\tauj{1}+\rho^{-2}\right)}
    \leq
    C[j].
  \end{equation}
Now, using the notation of Lemma \ref{perimm},
we can apply its system (i)
to estimate $g^t$, $g_t$, and $A^t$.

Actually the estimates become more transparent if we first rescale the system:
we set
  \begin{equation}
  \label{rescaled}
    \widetilde{g}^s:=\rho^2 g^{s/\rho(\cdot)},
    \qquad
    \widetilde{g}_s:=\rho^{-2}g_{s/\rho(\cdot)},
    \qquad \mbox{and} \qquad
    \widetilde{A}^s:=\rho A^{s/\rho(\cdot)},
  \end{equation}
so that by item (i) of Lemma \ref{perimm} and Remark \ref{spherical}
  \begin{equation}
  \label{rescaledsystem}
    \partial_s{\widetilde{g}^s}_{\alpha \beta}
    =
    -2{\widetilde{A}^s}_{\alpha \beta}
    \qquad \mbox{and} \qquad
    \partial_s{\widetilde{A}^s}_{\alpha \beta}
    =
    \rho^{-2}{\widetilde{g}^s}_{\alpha \beta}
      -{\widetilde{g}_s}^{\gamma \delta}
      {\widetilde{A}^s}_{\alpha \gamma}{\widetilde{A}^s}_{\beta \delta}
  \end{equation}
and by \eqref{initests}, \eqref{rhobounds}, item (ii) of Proposition \ref{catest},
and item (ii) of Proposition \ref{torest}
  \begin{equation}
  \label{IC&coeffests}
    \norm{\widetilde{g}^0: C^j\left({T^*\Sigma}^{\otimes 2},\chi\right)}
      +\norm{\widetilde{g}_0: C^j\left({T\Sigma}^{\otimes 2},\chi\right)}
      +\norm{\widetilde{A}^0: C^j\left({T^*\Sigma}^{\otimes 2},\chi\right)}
      +\norm{\rho^{-2}: C^j(\Sigma,\chi)}
    \leq
    C[j].
  \end{equation}
It follows from the system \eqref{rescaledsystem}
and the estimates \eqref{IC&coeffests}
on the initial conditions and coefficients
that there exists some $\epsilon>0$
such that the solution to the system exists
at all points $p \in \Sigma$ whenever $\abs{s}<\epsilon$
and moreover for any nonnegative integers $i$ and $j$
there exists a constant $C[i,j]$ such that
whenever $\abs{s} \leq \epsilon/2$
  \begin{equation}
  \label{rescaledests}
    \norm{{\partial_s}^i \widetilde{g}^s: C^j\left({T^*\Sigma}^{\otimes 2},\chi\right)}
     +\norm{{\partial_s}^i \widetilde{g}_s: C^j\left({T\Sigma}^{\otimes 2},\chi\right)}
     +\norm{{\partial_s}^i \widetilde{A}^s: C^j\left({T^*\Sigma}^{\otimes 2},\chi\right)}
    \leq
    C[i,j].
  \end{equation}

Now let $C_u>0$ and $u \in C^{2,\alpha,\gamma}_{\Grp}(\Sigma)$
with $\norm{u}_{2,\alpha,\gamma} \leq C_u\tauj{1}$.
By \eqref{weightedhoelderdef} and \eqref{globalweightednorm}
  \begin{equation}
    \norm{\rho u: C^{2,\alpha}\left(\Sigma,\chi,m^\gamma \tauj{1} \rho^{1-\gamma}\right)}
    \leq
    C_u,
  \end{equation}
so in particular by \eqref{rhobounds} and line 2 of \eqref{estimatesforforceestimate}
  \begin{equation}
    \norm{\rho u : C^{2,\alpha}(\Sigma,\chi)} \leq m^\gamma \tauj{1}^\gamma \leq \frac{\epsilon}{2},
  \end{equation}
provided $m$ is chosen large enough (in terms of $\epsilon>0$, $c$, and $C_u$).
Consequently we can apply the estimates \eqref{rescaledests}
along with the definitions \eqref{rescaled} to conclude that
for all $t \in [0,1]$
  \begin{equation}
  \label{guAuest}
    \norm{g^{tu}: C^{2,\alpha}\left({T^*\Sigma}^{\otimes 2},\chi, \rho^{-2}\right)}
     +\norm{g^{tu}: C^{2,\alpha}\left({T\Sigma}^{\otimes 2},\chi, \rho^2 \right)}
     +\norm{A^{tu}: C^{2,\alpha}\left({T^*\Sigma}^{\otimes 2},\chi, \rho^{-1}\right)}
    \leq
    C
  \end{equation}
for some constant $C=C[N,k,\ell]>0$ whenever $m>m_0$
for some $m_0=m_0[N,k,\ell,c]>0$.
Thus we also have
  \begin{equation}
  \label{HDeltagest}
    \norm{H^u: C^{2,\alpha}(\Sigma,\chi,\rho)}
      +\norm{g^u-g: C^{2,\alpha}\left({T^*\Sigma}^{\otimes 2},\chi, m^{\gamma}\tauj{1}\rho^{-1-\gamma}\right)}
    \leq
    C[N,k,\ell,C_u],
  \end{equation}
using (i) of Lemma \ref{perimm} to estimate the second norm.

Since $\chi=\rho^2g$,
  \begin{equation}
    \norm{D\left[T^*\Sigma,g\right]-D\left[T^*\Sigma,\chi\right]:
      C^j\left({T^*\Sigma}^{\otimes 2},\chi\right)}
    \leq
    C[j],
  \end{equation}
so, using also $\norm{u}_{2,\alpha,\gamma} \leq C_u\tauj{1}$,
\eqref{guAuest},
and the estimate for the second term of \eqref{HDeltagest},
  \begin{equation}
  \label{remainingests}
    \norm{D[g_u]A^{tu}: C^{1,\alpha}\left({T^*\Sigma}^{\otimes 3},\chi,\rho^{-1}\right)}
      +\norm{D[g^u]^2u: C^{0,\alpha}\left({T^*\Sigma}^{\otimes 2},\chi,\tauj{1}\right)}
    \leq
    C
  \end{equation}
for another constant $C=C[N,k,\ell,C_u]>0$,
whenever $\zeta,\xi \in [-c,c]^{N-1}$ and $m>m_0$
for some $m_0=m_0[N,k,\ell,c,C_u]>0$.
Applying \eqref{guAuest}, \eqref{remainingests}, and $\norm{u}_{2,\alpha,\gamma} \leq C_u\tauj{1}$
(as well as Remark \ref{spherical})
in item (iv) of Lemma \ref{perimm}
and then using \eqref{taudef} and \eqref{rhobounds},
for each $p \in \Sigma$ we obtain
  \begin{equation}
    \frac{\norm{\rho^{-2}\Qcal[u]:C^{0,\alpha}(B[p,1,\chi],\chi)}}{m^\gamma \rho^{-\gamma}}
    \leq
    C\tauj{1}^2m^\gamma \rho(p)^{1-\gamma}
    =
    C\tauj{1}^{1+\gamma/2}m^\gamma \tauj{1}^{1-\gamma/2}\tauj{1}^{\gamma-1}
    \leq
    Ce^{4c}\taubarj{1}^{1+\gamma/2}m^{\gamma}\tauj{1}^{\gamma/2},
  \end{equation}
where $B[p,1,\chi]$ is the $\chi$ metric ball of center $p$ and radius $1$
and $C=C[N,k,\ell,C_u]>0$ is yet another constant,
whenever $\zeta,\xi \in [-c,c]^{N-1}$
and $m$ is sufficiently large in terms of $N$, $k$, $\ell$, and $c$.
The proof is now concluded by invoking line 2 of \eqref{estimatesforforceestimate}.
\end{proof}

\subsection*{Forces through the perturbed surface}
Recall, in addition to \eqref{deformediota} and \eqref{Hcaldef},
the perturbed unit normal $\nu[u]=\nu[N,k,\ell,m,\zeta,\xi,u]$
defined just after \eqref{deformediota} for sufficiently small $u$.
For such $u$ and for each integer $i \in [1,N]$ we define
the force
  \begin{equation}
  \label{Ftildedef}
    \Ftilde_i
    =
    \Ftilde_i[N,k,\ell,m,\zeta,\xi,u]
    :=
    \int_{\Omega_i} \Hcal[u] \, (\gsph \circ \iota[u])(K \circ \iota[u],\nu[u])
      \, \sqrt{\abs{\iota[u]^*\gsph}},
  \end{equation}
the perturbation by $u$ of \eqref{forcedef},
where $\sqrt{\abs{\iota[u]^*\gsph}}$
is the area form induced by $\iota[u]$ and $\gsph$.
We will need the following estimate for $\Ftilde_i$.

\begin{lemma}[Estimates of the perturbations to the forces]
\label{perturbedforce}
Given $C_u,c>0$, $\alpha,\gamma \in (0,1)$,
and integers $N \geq 2$ and $\ell \geq k \geq 1$,
there exist real numbers
$\widetilde{c}:=\widetilde{c}[N,k,\ell,C_u]>0$
and $m_0:=m_0[N,k,\ell,m,C_u,c]>0$
such that (recalling \eqref{forcedef} and \eqref{Ftildedef})
  \begin{equation}
    \abs{\Ftilde_i[N,k,\ell,m,\zeta,\xi,u]-\Fcal_i[N,k,\ell,m,\zeta,\xi]}
    \leq
    \widetilde{c}m^{-2}\tauj{1}
  \end{equation}
whenever $1 \leq i \leq N$,
$m>m_0$, $\zeta,\xi \in [-c,c]^{N-1}$,
and $u \in C^{2,\alpha}(\Sigma,\chi)$ satisfies
$\norm{u}_{2,\alpha,\gamma} \leq C_u\tauj{1}$;
furthermore, for each fixed $i$, $N$, $k$, $\ell$, and $m>m_0$,
the map $(u,\zeta,\xi) \mapsto \Ftilde_i[N,k,\ell,m,\zeta,\xi,u]$
is continuous (in the sense of Definition \ref{functionalcty}).
\end{lemma}

We emphasize that in the statement of Lemma \ref{perturbedforce}
$\widetilde{c}$ does not depend on $c$ or $m$.

\begin{proof}
The continuity is clear from the smooth dependence of the initial surfaces
on the $\zeta,\xi$ parameters
and from definitions \eqref{deformediota} and \eqref{Hcaldef}.
Turning to the estimate, obviously
  \begin{equation}
  \label{Ftelescope}
    \begin{aligned}
      \Ftilde_i-\Fcal_i
      =
      &\int_{\Omega_i} \Hcal[u] \, (\gsph \circ \iota[u])(K \circ \iota[u],\nu[u])
        \, \left[\sqrt{\abs{g[u]}}-\sqrt{\abs{g}}\right] \\
      &+\int_{\Omega_i} \Hcal[u] \,
        \left[
            (\gsph \circ \iota[u])(K \circ \iota[u],\nu[u])-(\gsph \circ \iota)(K \circ \iota, \nu)
          \right] \,
        \sqrt{\abs{g}} \\
      &+\int_{\Omega_i} \left[\Hcal[u]-H\right] \, (\gsph \circ \iota)(K \circ \iota, \nu) \, \sqrt{\abs{g}},
    \end{aligned}
  \end{equation}
using the notation of Lemma \ref{perimm}.

From \eqref{Hcaldef}, \eqref{Lchidef},
Corollary \ref{Hest},
Lemma \ref{quadest},
\eqref{taudef}, \eqref{Ddef}, 
and the assumption that $\norm{u}_{2,\alpha,\gamma} \leq C_u\tauj{1}$
  \begin{equation}
    \norm{\rho^{-2}\Hcal[u]}_{0,\alpha,\gamma}
    =
    \norm{\rho^{-2}H + \Lchi u + \rho^{-2}\Qcal[u]}_{0,\alpha,\gamma}
    \leq
    m\tauj{1}
  \end{equation}
whenever $\zeta,\xi \in [-c,c]^{N-1}$
and $m$ is sufficiently large in terms of $N$, $k$, $\ell$, $C_u$, and $c$.
By \eqref{Killing}
  \begin{equation}
    \abs{(\gsph \circ \iota[u])(K \circ \iota[u],\nu[u])}
      +\abs{(\gsph \circ \iota)(K \circ \iota, \nu)}
    \leq
    2
  \end{equation}
and, using also \eqref{phiunormal},
the proof of Lemma \ref{quadest} (particularly \eqref{guAuest}),
and again the assumption $\norm{u}_{2,\alpha,\gamma} \leq K\tauj{1}$,
  \begin{equation}
    \norm{(\gsph \circ \iota[u])(K \circ \iota[u],\nu[u])
      -(\gsph \circ \iota)(K \circ \iota, \nu):
      C^{1,\alpha}\left(\Sigma,\chi,m^\gamma \rho^{1-\gamma}\right)}
    \leq
    C[N,k,\ell,C_u]\tauj{1}.
  \end{equation}
By item (ii) of Lemma \ref{perimm} and the proof of Lemma \ref{quadest}
(particularly the estimate of the second term in \ref{HDeltagest})
  \begin{equation}
    \norm{\sqrt{\abs{g[u]}}-\sqrt{\abs{g}}: C^{1,\alpha}\left(\Sigma,\chi,m^\gamma \rho^{-1-\gamma}\right)}
    \leq
    C[N,k,\ell,C_u]\tauj{1}.
  \end{equation}
Finally, for the $\chi$ area $\abs{\Omega_i}_\chi$ of $\Omega_i$
we have, whenever $\zeta,\xi \in [-c,c]^{N-1}$
and $m$ is sufficiently large in terms of $N$, $k$, $\ell$, and $c$, the estimate
  \begin{equation}
  \label{Omegaiarea}
    \begin{aligned}
      \abs{\Omega_i}_\chi
      &\leq
      \abs{\torr[i] \backslash (\catr[i-1] \cup \catr[i])}_\chi
        +\abs{\catr[i-1]}_\chi + \abs{\catr[i]}_\chi \\
      &=
      \int_{-\frac{\pi}{\sqrt{2}\ell}}^{\frac{\pi}{\sqrt{2}\ell}}
        \int_{-\frac{\pi}{\sqrt{2}k}}^{\frac{\pi}{\sqrt{2}k}}
        \left(1+Cm^2\tauj{1}\right) \, d\xx \, d\yy
      +
      4 \int_0^{2\pi} \int_0^{m^2+2c} \left(1+Cm^2\tauj{1}\right) \, dt \, d\theta
      \leq
      200m^2,
    \end{aligned}
  \end{equation}
recalling \eqref{Omega}, \eqref{catrdef}, and \eqref{torrdef},
understanding $\catr[0]=\catr[N]=\emptyset$,
and using \eqref{cisummary}, \eqref{adef}
and Propositions \ref{catest} and \ref{torest},
which supply the constant $C=C[N,k,\ell]$.

It now follows from the estimates of the previous paragraph
(and \eqref{rhobounds} and line 2 of \eqref{estimatesforforceestimate})
that, whenever $\zeta, \xi \in [-c,c]^{N-1}$
and $m$ is sufficiently large in terms of $N$, $k$, $\ell$, $c$, and $C_u$,
  \begin{equation}
  \label{Ffirsttwo}
    \begin{aligned}
      &\abs{
        \int_{\Omega_i} \Hcal[u] \, (\gsph \circ \iota[u])(K \circ \iota[u],\nu[u])
        \, \left[\sqrt{\abs{g[u]}}-\sqrt{\abs{g}}\right]
      }
      \leq
      Cm^{3+2\gamma}\tauj{1}^2\norm{\rho^{1-2\gamma}}_{C^0(\Sigma)}
      \leq
      m^{-2}\tauj{1} \mbox{ and} \\
      &\abs{
        \int_{\Omega_i} \Hcal[u] \,
          \left[
            (\gsph \circ \iota[u])(K \circ \iota[u],\nu[u])-(\gsph \circ \iota)(K \circ \iota, \nu)
          \right] \,
        \sqrt{\abs{g}}
      }
      \leq
      Cm^{3+2\gamma}\tauj{1}^2\norm{\rho^{1-2\gamma}}_{C^0(\Sigma)}
      \leq
      m^{-2}\tauj{1}
    \end{aligned}
  \end{equation}
(regardless of the sign of $1-2\gamma$).
Furthermore, using also Lemma \ref{quadest}
(as well as \eqref{taudef}),
  \begin{equation}
    \abs{\int_{\Omega_i} \Qcal[u] \, (\gsph \circ \iota)(K \circ \iota,\nu) \, \sqrt{\abs{g}}}
    \leq
    Cm^2 \taubarj{1}^{1+\gamma/2}
    \leq
    m^{-2}\tauj{1},
  \end{equation}
again for $m$ sufficiently large in terms of $N$, $k$, $\ell$, $c$, and $C_u$.
Therefore
  \begin{equation}
  \label{Fquad}
    \abs{
    \int_{\Omega_i} \left[\Hcal[u]-H\right] \, (\gsph \circ \iota)(K \circ \iota, \nu) \, \sqrt{\abs{g}}
    }
    \leq
    \abs{\int_{\Omega_i} (\gsph \circ \iota)(K \circ \iota, \nu) \, \Lcal u \, \sqrt{\abs{g}}}
      +m^{-2}\tauj{1},
  \end{equation}
recalling \eqref{Lcaldef},
but $K$ is Killing,
so integration by parts (specifically Green's identity) yields
  \begin{equation}
  \label{FGreen}
  \begin{aligned}
    \int_{\Omega_i} (\gsph \circ \iota)(K \circ \iota, \nu) \, \Lcal u \, \sqrt{\abs{g}}
    =
    &-\int_{\Omega_i} u \, (\gsph \circ \iota)(K \circ \iota,\iota_*\nabla_g H) \, \sqrt{\abs{g}} \\
    &+\int_{\partial \Omega_i}
      \left[
        (\gsph \circ \iota)(K \circ \iota, \nu) (\eta u)
        -u \left(\eta \left[(\gsph \circ \iota)(K \circ \iota, \nu)\right]\right)
      \right]
      \, \sqrt{\abs{g}},
  \end{aligned}
  \end{equation}
where $\eta$ is the outward conormal on $\Omega_i$ induced by $g$
(and acts on functions as a derivation).

Using \eqref{Killing}, \eqref{Ddef}, \eqref{wbar}, \eqref{C1Hest}, and \eqref{Omegaiarea}, 
it follows that
  \begin{equation}
  \label{dHest}
    \abs{\int_{\Omega_i} u \, (\gsph \circ \iota)(K \circ \iota,\iota_*\nabla_g H) \, \sqrt{\abs{g}}}
    \leq
    CC_u m^{2+\gamma} \tauj{1}
      \left(\norm{m^2\rho^{-1-\gamma}\tauj{1}+m^2\rho^{1-\gamma}\tauj{1}^2}_0+cm^{1-\gamma}\tauj{1}\right)
    \leq
    m^{-2}\tauj{1}
  \end{equation}
for $m$ sufficiently large in terms of $N$, $k$, $\ell$, $C_u$, and $c$.
Turning to the boundary term, as in the computation following \eqref{balance},
$\partial \Omega_i$ has one or two circular components
(catenoidal waists)
and a single rectangular component.
Suppose $S:=\cat_i(\{t=0\})$ or $S:=\cat_{i-1}(\{t=0\})$ is a circular component
and $T:=\partial \Omega_i \backslash [\cat_{i-1}(\{t=0\}) \cup \cat_i(\{t=0\})]$
is the rectangular component.
By \eqref{Phidef}, \eqref{Omega}, and \eqref{Killing}
  \begin{equation}
  \label{normalprojectionTest}
    m^2\norm{(\gsph \circ \iota)(K \circ \iota, \nu)-1: C^0(T)}
      +m\norm{\eta \left[ (\gsph \circ \iota)(K \circ \iota, \nu)\right]: C^0(T)}
    \leq
    C
  \end{equation}
for some constant $C=C[N,k,\ell]>0$
and obviously $T$ has $g$ length $\abs{T}|_{g|T} \leq Cm^{-1}$
and $\norm{u: C^0(T)} \leq C_u \tauj{1}$,
so
  \begin{equation}
  \label{FTu}
    \abs{\int_T u \left(\eta \left[(\gsph \circ \iota)(K \circ \iota, \nu)\right]\right) \, \sqrt{\abs{g|_T}}}
    \leq
    Cm^{-2}\tauj{1}
  \end{equation}
for another constant $C=C[N,k,\ell,C_u]>0$,
while $\norm{\eta u : C^0(T)} \leq C_u m\tauj{1}$,
so
  \begin{equation}
  \label{FTetau}
    \abs{\int_T \left[(\gsph \circ \iota)(K \circ \iota, \nu)\right] (\eta u) \, \sqrt{\abs{g|_T}}}
    \leq
    Cm^{-2}\tauj{1}
      +\abs{\int_T (\eta u) \, \sqrt{\abs{g|_T}}}
    \leq
    Cm^{-2}\tauj{1},
  \end{equation}
where for the first inequality we have used \eqref{normalprojectionTest}
and for the second we have used the fact that,
because $u$ is $\Grp$-equivariant,
it satisfies periodic boundary conditions on $T$
and accordingly $\int_T (\eta u) \, \sqrt{\abs{g|_T}}=0$.

On the other hand,
on $S$ we have $\norm{u: C^0(S)} \leq C_u m^\gamma \tauj{1}^{1+\gamma}$,
$\norm{\eta u: C^0(S)} \leq C_u m^\gamma \tauj{1}^\gamma$,
and
  \begin{equation}
    \norm{\eta \left[(\gsph \circ \iota)(K \circ \iota, \nu)\right]: C^0(S)}
    \leq
    \sup_S \left(\abs{D[\gsph]K}_{\gsph} + \abs{K}_{\gsph}\abs{A}_g\right)
    \leq
    C\tauj{1}^{-1}
  \end{equation}
for some constant $C=C[N,k,\ell]>0$,
having used item (v) of Proposition \ref{catest}
for the last inequality,
and $S$ has $g$ length $\abs{S}_{g|_S} \leq C\tauj{1}$,
so
  \begin{equation}
  \label{FS}
    \int_S
      \left(
        \abs{u \left(\eta \left[(\gsph \circ \iota)(K \circ \iota, \nu)\right]\right)}
        +\abs{\left[(\gsph \circ \iota)(K \circ \iota, \nu)\right] (\eta u)}
      \right) \, \sqrt{\abs{g|_S}}
    \leq
    C[N,k,\ell,C_u]m^{\gamma}\tauj{1}^{1+\gamma}
    \leq
    m^{-2}\tauj{1},
  \end{equation}
provided $m$ is sufficiently large in terms of $N$, $k$, $\ell$, $c$, and $C_u$,
yet again using line 2 of \eqref{estimatesforforceestimate} for the last inequality.
The proof is now completed by combining
\eqref{Ftelescope}, \eqref{Ffirsttwo}, \eqref{Fquad},
\eqref{FGreen}, \eqref{dHest}, \eqref{FTu}, \eqref{FTetau}, and \eqref{FS}.
\end{proof}

\subsection*{Explicitly defined diffeomorphisms between initial surfaces}
Recall Remark \ref{initsurfcty} and Definition \ref{functionalcty}.
Throughout the construction we have made use of the existence
of maps $I[N,k,\ell,m]$ as in Remark \ref{initsurfcty}
in order to identify function spaces defined on initial surfaces
with identical data $N$, $k$, $\ell$, and $m$
but different $\zeta,\xi$ parameter values.
So far we have made these identifications
merely so as to articulate certain continuity properties,
which do not depend on the choice of $I$.
In the proof of the main theorem, however,
we will need bounds for the norms of the corresponding identification maps
between our normed function spaces,
and so we now explicitly define diffeomorphisms between the initial surfaces.
We define these diffeomorphisms
as compromises between natural identifications on the various standard regions.
More precisely,
recalling \eqref{adef}, \eqref{catdef}, \eqref{catrdef}, \eqref{Tdef}, and \eqref{torrdef},
for any any given data $N$, $k$, $\ell$, $m$, and $\zeta,\xi$
we start by defining
  \begin{equation}
  \label{bardefs}
    \begin{aligned}
      &a_i:=a_i[N,k,\ell,m,\zeta,\xi] \qquad \mbox{and} \qquad
      \abar_i:=a_i[N,k,\ell,m,0,0] \qquad \mbox{for $1 \leq i \leq N-1$,} \\
      &\catr[i]:=\catr[i;N,k,\ell,m,\zeta,\xi] \qquad \mbox{and} \qquad
      \underline{\catr}[i]:=\catr[i;N,k,\ell,m,0,0] \qquad \mbox{for $1 \leq i \leq N-1$,} \\
      &\cat_i:=\cat_i[N,k,\ell,m,\zeta,\xi] \qquad \mbox{and} \qquad
      \underline{\cat}_i:=\cat_i[N,k,\ell,m,0,0] \qquad \mbox{ for $1 \leq i \leq N-1$,} \\
      &\torr[i]:=\torr[i;N,k,\ell,m,\zeta,\xi] \qquad \mbox{and} \qquad
      \underline{\torr}[i]:=\torr[i;N,k,\ell,m,0,0] \qquad \mbox{for $1 \leq i \leq N$, and} \\
      &T_i:=T_i[N,k,\ell,m,\zeta,\xi] \qquad \mbox{and} \qquad
      \underline{T}_i:=T_i[N,k,\ell,m,0,0] \qquad \mbox{for $1 \leq i \leq N$}.
    \end{aligned}
  \end{equation}

We observe that the map
  \begin{equation}
    \left.T_i \circ \underline{T}_i^{-1}\right|_{\underline{\torr}[i] \backslash
      (\underline{\catr}[i-1] \cup \underline{\catr}[i])}:
    \underline{\torr}[i] \backslash (\underline{\catr}[i-1] \cup \underline{\catr}[i])
    \to
    \torr[i] \backslash (\catr[i-1] \cup \catr[i]).
  \end{equation}
(understanding $\catr[0]=\catr[N]=\emptyset$)
is a well-defined diffeomorphism.
We also observe (recalling \eqref{cyladef}) that
whenever $\torr[i] \cap \catr[j] \neq \emptyset$,
the map $T_i \circ \underline{T}_i^{-1} \circ \underline{\cat}_j$
is well-defined on the component of $\cyl_{\abar_j} \backslash \cyl_{\abar_j-1}$
whose image under $\underline{\cat}_j$ lies in $\underline{\torr}[i]$
and that on this set
$T_i \circ \underline{T}_i^{-1} \circ \underline{\cat}_j$
has image contained in $\catr[j]$ and moreover satisfies
  \begin{equation}
    \left(T_i \circ \underline{T}_i^{-1} \circ \underline{\cat}_j\right)(\underline{t},\theta)
    =
    \cat_j
      \left(
        (\sgn \underline{t}) \arcosh \left[\frac{\taubarj{j}}{\tauj{j}} \cosh \underline{t}\right],\theta
      \right),
  \end{equation}
where $\sgn: \R \to \R$ takes the value $1$ when its argument is nonnegative and $-1$ otherwise.
Note that (using the identity \eqref{lnarcosh})
  \begin{equation}
    \arcosh \left(\frac{\taubarj{j}}{\tauj{j}} \cosh \underline{t}\right)
    =
    \abs{\underline{t}}
     +\ln \frac{\taubarj{j}}{2\tauj{j}}
     +\ln \left(1+e^{-2\abs{\underline{t}}}\right)
     +\ln \left(1+\sqrt{1-\tauj{j}^2\taubarj{j}^{-2} \sech^2 \underline{t}}\right)
  \end{equation}

So motivated,
for $1 \leq j \leq N-1$
we define the function
$\widetilde{t}_j: \R \to \R$
by
  \begin{equation}
  \label{ttildedef}
    \widetilde{t}_j(\underline{t})
    :=
    \frac{a_j}{\abar_j}\underline{t} \cdot \cutoff{\abar_j}{\abar_j-1}(\abs{\underline{t}})
      +(\sgn \underline{t})\arcosh \left(\frac{\taubarj{j}}{\tauj{j}} \cosh \underline{t}\right)
        \cdot \cutoff{\abar_j-1}{\abar_j}(\abs{\underline{t}}).
  \end{equation}
Since 
  \begin{equation}
  \label{ttildeests}
    \begin{aligned}
      &\frac{d}{d\underline{t}}\arcosh \left(\frac{\taubarj{j}}{\tauj{j}} \cosh \underline{t}\right)
      =
      \frac{\tanh \underline{t}}{\sqrt{\taubarj{j}^2 \cosh^2 \underline{t}-\tauj{j}^2}}
        \, \taubarj{j} \cosh \underline{t}, \\
      &\arcosh \left(\frac{\taubarj{j}}{\tauj{j}} \cosh \pm \abar_j\right)=\pm a_j,
      \quad \mbox{and by \eqref{adef}} \quad
      \abs{\frac{a_j}{\abar_j}-1} \leq \frac{2c}{m},
    \end{aligned}
  \end{equation}
we see using also \eqref{taudef} and \eqref{adef}
that by taking $m$ sufficiently large in terms of $c$
we can guarantee that
$\widetilde{t}$ takes $[-\abar_j,\abar_j]$
monotonically onto $[-a_j,a_j]$.
Away from the ends of $[-\abar_j,\abar_j]$
this reparametrization is simply multiplication by $a_j/\abar_j \approx 1$,
while close to the ends it almost agrees with the map
$\underline{t} \mapsto \underline{t}+(\sgn \underline{t})(a_j-\abar_j)$.

We can now define the diffeomorphism
  \begin{equation}
    P[\zeta,\xi]=P[N,k,\ell,m,\zeta,\xi]:
    \Sigma[N,k,\ell,m,0,0] \to \Sigma[N,k,\ell,m,\zeta,\xi]
  \end{equation}
by requiring that (i) it commute with the action of $\Grp$ (recalling \eqref{Grpdef}),
  \begin{equation}
  \label{Pdef}
    \begin{aligned}
      &\mbox{(ii)} \quad \mbox{for $1 \leq i \leq N$} \quad
      P[\zeta,\xi]|_{\underline{\torr}[i] \backslash (\underline{\catr}[i-1] \cup \underline{\catr}[i])}
      :=
      T_i \circ \underline{T}_i^{-1} \quad \mbox{and} \\
      &\mbox{(iii)} \quad \mbox{for $1 \leq i \leq N-1$} \quad
      \left(P[\zeta,\xi] \circ \underline{\cat}_i\right)(\underline{t},\theta)
      :=
      \cat_i\left(\widetilde{t}(\underline{t}),\theta\right)
        \mbox{ for all } (\underline{t},\theta) \in \cyl_{\underline{a}_i}
    \end{aligned}
  \end{equation}
(continuing to understand $\catr[0]=\catr[N]=\emptyset$).
We define in turn the map
  \begin{equation}
  \label{Pcal}
    \Pcal=\Pcal[\zeta,\xi]=\Pcal[N,k,\ell,m,\zeta,\xi]:=P[N,k,\ell,m,\zeta,\xi]^*
  \end{equation}
taking functions on $\Sigma[N,k,\ell,m,\zeta,\xi]$ to functions on $\Sigma[N,k,\ell,m,0,0]$.
Clearly the map $I=I[N,k,\ell]: \R^{N-1} \times \R^{N-1} \times \Sigma[N,k,\ell,0,0] \to \Sph^3$
defined by $I(\zeta,\xi,\cdot):=\iota[N,k,\ell,m,\zeta,\xi] \circ P[N,k,\ell,m,\zeta,\xi]$
satisfies the properties specified in Remark \ref{initsurfcty}.
Last we have the following estimate.

\begin{lemma}[Bound for $\Pcal$ and its inverse]
\label{Pcalest}
Given real numbers $\alpha,\gamma \in (0,1)$ and $c>0$
as well as integers $N \geq 2$ and $\ell \geq k \geq 1$,
there exist real numbers $C=C[N,k,\ell,\alpha,\gamma]>0$
and $m_0=m_0[N,k,\ell,c]>0$
such that
whenever $\zeta, \xi \in [-c,c]^{N-1}$
and $m>m_0$
we have (recalling \eqref{globalweightednorm}) the estimates
  \begin{equation}
    \begin{aligned}
      &\norm{\Pcal[\zeta,\xi]u}_{C^{2,\alpha,\gamma}(\Sigma[N,k,\ell,m,0,0])}
        \leq Ce^{2c}\norm{u}_{C^{2,\alpha,\gamma}(\Sigma[N,k,\ell,m,\zeta,\xi])} \mbox{ and} \\
      &\norm{\Pcal[\zeta,\xi]^{-1}v}_{C^{2,\alpha,\gamma}(\Sigma[N,k,\ell,m,\zeta,\xi])}
        \leq Ce^{2c}\norm{v}_{C^{2,\alpha,\gamma}(\Sigma[N,k,\ell,m,0,0])}.
    \end{aligned}
  \end{equation}
\end{lemma}

\begin{proof}
Let $u \in C^{2,\alpha,\gamma}(\Sigma[N,k,\ell,m,\zeta,\xi])$.
By \eqref{Tipre}, \eqref{Tdef}, and \eqref{bardefs}
  \begin{equation}
    T_i^{-1}(\torr[i] \backslash (\catr[i-1] \cup \catr[i]))
    =
    \underline{T}_i^{-1}(\underline{\torr}[i] \backslash (\underline{\catr}[i-1] \cup \underline{\catr}[i])),
  \end{equation}
and by \eqref{Pdef} and \eqref{Pcal}
  \begin{equation}
    \left(T_i^*u\right)(\xx,\yy)=\left(\underline{T}_i^*\left(\Pcal u\right)\right)(\xx,\yy)
      \mbox{ for all } (\xx,\yy) \in T_i^{-1}(\torr[i] \backslash (\catr[i-1] \cup \catr[i])),
  \end{equation}
while by \eqref{ttildedef}
  \begin{equation}
    \left(\cat_i^*u\right)\left((\widetilde{t}(\underline{t}),\theta\right)
    =
    \left(\underline{\cat}_i^*\left(\Pcal u\right)\right)(\underline{t},\theta)
      \mbox{ for all } (\underline{t},\theta) \in \cyl_{\underline{a}_i}
      \mbox{ (equivalently all $\left(\widetilde{t}(\underline{t},\theta\right) \in \cyl_{a_i}$)}.
  \end{equation}
The asserted bounds are now clear from \eqref{globalweightednorm},
using \eqref{taudef}, \eqref{rhodef}, \eqref{ttildeests},
and Propositions \ref{catest} and \ref{torest}.
\end{proof}

\subsection*{The main theorem}
We are ready to prove the main theorem.

\begin{theorem}[The main theorem]
\label{mainthm}
Let $\alpha,\gamma \in (0,1)$.
Given integers $N \geq 2$ and $\ell \geq k \geq 1$,
there are real numbers
$C, \cbar, m_0>0$
such that for every $m>m_0$
there exist parameters $\zeta,\xi \in [-\cbar,\cbar]^{N-1}$
and a function $u \in C_{\Grp}^\infty(\Sigma[N,k,\ell,m,\zeta,\xi])$
(recalling \eqref{Grpdef}, \eqref{initsurfdef}, and Notation \ref{functionspaceappendages})
such that
$\norm{u}_{2,\alpha,\gamma} \leq C\tauj{1}$
(recalling \eqref{taudef} and \eqref{globalweightednorm})
and the image of the normal deformation $\iota[u]: \Sigma \to \Sph^3$
(recalling \eqref{deformediota})
by $u$ of the inclusion
$\iota: \Sigma \to \Sph^3$
is a closed embedded minimal surface
invariant under $\Grp[k,\ell,m]$ and having genus $k\ell m^2(N-1)+1$.
\end{theorem}

\begin{proof}
Fix $\alpha,\gamma \in (0,1)$ and integers $N \geq 2$ and $\ell \geq k \geq 1$.
For each integer $m \geq 1$ set
  \begin{equation}
  \label{functionball}
    B[N,k,\ell,m]
    :=
    \left\{
      v \in C_{\Grp[k,\ell,m]}^{2,\alpha/2}\left(\Sigma[N,k,\ell,m,0,0],\chi \right)
      \; : \;
      \norm{v}_{2,\alpha,\gamma} \leq \taubarj{1}^{1+\gamma/3}
    \right\}
  \end{equation}
 (recalling \eqref{taubardef}).
 Given $\zeta,\xi \in \R^{N-1}$ and assuming $m$ sufficiently large, define also
  \begin{equation}
  \label{1}
    \left(u_1,\left(\lambda_1,\ldots,\lambda_N\right),\left(\lambdabar_2,\ldots,\lambdabar_{N-1}\right)\right)
    :=
    -\Rcal
      \left(
        \rho^{-2}H
        -\sum_{i=2}^{N-1} \Dcal_i \wbar_i
      \right)[N,k,\ell,m,\zeta,\xi],
  \end{equation}
 as in Corollary \ref{firstorder}
(recalling \eqref{Ddef}, \eqref{AHdef}, \eqref{rhodef}, \eqref{wbar}, and Proposition \ref{globalsol}),
and for each $v \in B[N,k,\ell,m]$ define
  \begin{equation}
  \label{2}
    \left(v',\left(\mu_1[v],\ldots,\mu_N[v]\right),\left(\mubar_2[v],\ldots,\mubar_{N-1}[v]\right)\right)
    :=
    -\Rcal \left(\rho^{-2}\Qcal\left[u_1+\Pcal^{-1}v\right]\right)[N,k,\ell,m,\zeta,\xi]
  \end{equation}
(recalling \eqref{Qcal} and \eqref{Pcal}).

Thus, for all $\zeta, \xi \in \R^{N-1}$ and $v \in B[N,k,\ell,m]$,
provided $m$ is sufficiently large in terms of $N$, $k$, $\ell$, and $\zeta,\xi$,
recalling \eqref{Hcaldef}, \eqref{Lcaldef}, \eqref{Lchidef}, \eqref{Qcal}, and Proposition \ref{globalsol},
  \begin{equation}
  \label{motivation}
    \begin{aligned}
      \rho^{-2}\Hcal\left[u_1+v',\zeta,\xi \right]
      &=
      \rho^{-2}H+\Lchi (u_1+v') + \rho^{-2} \Qcal[u_1+v'] \\
      &=
      \left[\Lchi u_1+\left(\rho^{-2}H-\sum_{i=2}^{N-1} \Dcal_i \wbar_i\right)\right]
        + \rho^{-2}\Qcal[u_1+v']+\Lchi v' + \sum_{i=2}^{N-1} \Dcal_i \wbar_i \\
      &=
      \rho^{-2}\Qcal\left[u_1+v'\right]-\rho^{-2}\Qcal\left[u_1+\Pcal^{-1}v\right] \\
      &\;\;\;\;+\sum_{i=1}^N \left(\lambda_i+\mu_i[v]\right)w_i
        +\sum_{i=2}^{N-1} \left(\Dcal_i+\lambdabar_i+\mubar_i[v]\right)\wbar_i.
    \end{aligned}
  \end{equation}
Evidently we want to pick $(v,\zeta,\xi)$ so that
$\Pcal v'=v$ (to make the nonlinear terms cancel),
$\lambda_i+\mu_i[v]=0$ for all $i \in \Z \cap [1,N]$ (to make the $w_i$ terms vanish),
and $\Dcal_i+\lambdabar_i+\mubar_i[v]=0$ for all $i \in \Z \cap [2,N-1]$ (to make the $\wbar_i$ terms vanish).
Recalling \eqref{Killing},
the unit normal $\nu$ for $\Sigma$ specified just above \eqref{AHdef}, and \eqref{w},
we observe that on the support of $w_i|_{\Omega_i}$
the function $(\gsph \circ \iota)(K \circ \iota,\nu)$ has a sign
(namely $(-1)^{N-1}$)
and the function $w_i$ itself is nonnegative.
Consequently, recalling \eqref{Ftildedef},
if $\Pcal v'=v$
and $\Dcal_i+\lambdabar_i+\mubar_i[v]=0$ for all $i \in \Z \cap [2,N-1]$,
then,
for any given $i \in \Z \cap [1,N]$,
$\lambda_i+\mu_i[v]=0$ if and only if $\Ftilde_i=0$.
Accordingly, recalling \eqref{1}, \eqref{2}, and Lemma \ref{FD},
we seek a fixed point for the map
  \begin{equation}
  \label{Jcal}
    \begin{aligned}
    &\Jcal: B[N,k,\ell,m] \times \R^{2N-2}
    \to C^{2,\alpha/2}_{\Grp[k,\ell,m]}(\Sigma[N,k,\ell,m,0,0],\chi) \times \R^{2N-2} \\
    &\Jcal \left(v,\begin{pmatrix}\zeta_1\\\vdots\\\zeta_{N-1}\\\xi_1\\\vdots\\\xi_{N-1}\end{pmatrix}\right)
    \mapsto
    \left(\Pcal v',
      \begin{pmatrix}\zeta_1\\\vdots\\\zeta_{N-1}\\\xi_1\\\vdots\\\xi_{N-1}\end{pmatrix}
      -
      \Theta^{-1}\tauj{1}^{-1}
      \begin{pmatrix}
        m^2\Ftilde_1[N,k,\ell,m,\zeta,\xi,u_1+\Pcal^{-1}v] \\
        \vdots \\
        m^2\Ftilde_N[N,k,\ell,m,\zeta,\xi,u_1+\Pcal^{-1}v] \\
        \Dcal_2[N,k,\ell,m,\zeta,\xi]+\lambdabar_2+\mubar_2[v] \\
        \vdots \\
        \Dcal_{N-1}[N,k,\ell,m,\zeta,\xi]+\lambdabar_{N-1}+\mubar_{N-1}[v]
      \end{pmatrix}
    \right)
    \end{aligned}
  \end{equation}

We will check that the hypotheses of the Schauder fixed-point theorem apply to $\Jcal$,
after restricting its domain as specified below.
It is clear from definition \eqref{Ddef}
and from the continuity assertions made in Proposition \ref{globalsol},
Corollary \ref{firstorder}, Lemma \ref{quadest}, and Lemma \ref{perturbedforce}
that $\Jcal$ is continuous in the sense of Definition \ref{functionalcty},
with the product topology on the domain and target,
the Euclidean topology on the $\R^{2N-2}$ factors,
and the $C^{2,\alpha/2}$ topology on the function-space factors.
Because each initial surface is compact,
the topology of each H\"{o}lder space is independent of the underlying metric
and $C^{2,\alpha}(\Sigma)$ embeds compactly in $C^{2,\alpha/2}(\Sigma)$
(as does the former's $\Grp$-equivariant subspace);
therefore $B[N,k,\ell,m]$ is compact relative to the $C^{2,\alpha/2}$ topology
and is clearly convex.

Now let $C_R$ be the constant $C[N,k,\ell,\alpha,\gamma]$ from Proposition \ref{globalsol},
let $C_1$ be the constant $C[N,k,\ell]$ from Corollary \ref{firstorder},
let $C_P$ be the constant $C[N,k,\ell]$ from Lemma \ref{Pcalest},
let $\widetilde{c}$ be the constant $\widetilde{c}[N,k,\ell,2C_1]$ from Lemma \ref{perturbedforce},
let $C_{\Theta}$ be the constant $C[N,k,\ell]$ from Lemma \ref{FD},
let
  \begin{equation}
  \label{cbar}
    \cbar:=C_\Theta\left(C_\Theta+\sqrt{N} \sqrt{\widetilde{c}^2+4C_1^2}\right),
  \end{equation}
and let $m_1$ be the maximum of the quantities named $m_0[N,k,\ell,\underline{c}]$ from
Proposition \ref{initsurfprops}, Lemma \ref{FD}, and Proposition \ref{globalsol}
as well as the quantity named $m_0[N,k,\ell,\underline{c},\gamma]$ from Corollary \ref{firstorder}
and the quantities named $m_0[N,k,\ell,m,2C_1,\underline{c}]$
from Lemma \ref{quadest} and Lemma \ref{perturbedforce}.

Suppose $m>m_1$, $\zeta,\xi \in [-\cbar,\cbar]$, and $v \in B[N,k,\ell,m]$.
Then by \eqref{functionball}, \eqref{1}, Corollary \ref{firstorder}, and Lemma \ref{Pcalest}
  \begin{equation}
  \label{1est}
    \sum_{i=1}^N \abs{\lambda_i}+\sum_{i=2}^{N-1} \abs{\lambdabar_i} \leq C_1\tauj{1}
    \qquad \mbox{and}\qquad
    \norm{u_1}_{2,\alpha,\gamma}+\norm{\Pcal^{-1}v}_{2,\alpha,\gamma}
    \leq
    C_1\tauj{1}+C_Pe^{2\cbar}\taubarj{1}^{1+\gamma/3}
    \leq
    2C_1\tauj{1},
  \end{equation}
where for the last inequality we use \eqref{taubardef}, \eqref{taudef},
and line 2 of \eqref{estimatesforforceestimate}
and we assume $m>m_2$ for some
$m_2=m_2[N,k,\ell,\gamma] \geq m_1$.
It follows in turn,
using \eqref{2}, Proposition \ref{globalsol}, Lemma \ref{quadest}, and Lemma \ref{Pcalest},
that
  \begin{equation}
  \label{2est}
    \norm{\Pcal v'}_{2,\alpha,\gamma}+\sum_{i=1}^N \abs{\mu_i[v]}+\sum_{i=2}^{N-1} \abs{\mubar_i[v]}
    \leq
    C_Pe^{2c}C_R\taubarj{1}^{1+\gamma/2}
    \leq
    \taubarj{1}^{1+\gamma/3},
  \end{equation}
assuming, for the last inequality, that $m>m_3$ for some $m_3=m_3[N,k,\ell,\gamma] \geq m_2$.
In particular we have verified that
  \begin{equation}
  \label{Bpreservation}
    v \in B[N,k,\ell,m] \Rightarrow \Pcal v' \in B[N,k,\ell,m].
  \end{equation}

Continuing to assume $m>m_3$,
from Lemma \ref{FD}, Lemma \ref{perturbedforce}, \eqref{cbar}, \eqref{1est}, and \eqref{2est}
we find that for any $\zeta,\xi \in [-\cbar,\cbar]^{N-1}$
  \begin{equation}
    \abs{\begin{pmatrix}\zeta_1\\\vdots\\\zeta_N\\\xi_1\\\vdots\\\xi_N\end{pmatrix}
      -
      \Theta^{-1}\tauj{1}^{-1}
      \begin{pmatrix}
        m^2\Ftilde_1[\zeta,\xi,u_1+\Pcal^{-1}v] \\
        \vdots \\
        m^2\Ftilde_N[\zeta,\xi,u_1+\Pcal^{-1}v] \\
        \Dcal_2[\zeta,\xi]+\lambdabar_2+\mubar_2[v] \\
        \vdots \\
        \Dcal_{N-1}[\zeta,\xi]+\lambdabar_{N-1}+\mubar_{N-1}[v]
      \end{pmatrix}
    }
    \leq
    C_\Theta^2 + C_\Theta \sqrt{N\widetilde{c}^2+4(N-2)C_1^2}
    \leq
    \cbar,
  \end{equation}
where the norm $\abs{\cdot}$ is the Euclidean one on $\R^{2N-2}$
and each $\Ftilde_i[N,k,\ell,m,\zeta,\xi,u]$ and $\Dcal_i[N,k,\ell,m,\zeta,\xi]$
have been abbreviated to $\Ftilde_i[\zeta,\xi,u]$ and $\Dcal_i[\zeta,\xi]$ respectively.
In conjunction with \eqref{Bpreservation}
this bound shows that $\Jcal$ (defined in \eqref{Jcal})
maps $B[N,k,\ell,m] \times [-\cbar,\cbar]^{2N-2}$ to itself.
Moreover, it is immediately clear from our observations in the paragraph following \eqref{Jcal}
that $\Jcal$ is continuous and $B[N,k,\ell,m] \times [-\cbar,\cbar]^{2N-2}$ is compact
relative to the topology described there,
and of course $B[N,k,\ell,m] \times [-\cbar,\cbar]^{2N-2}$ is convex.

The Schauder fixed-point theorem therefore applies
to guarantee the existence
of a fixed point $(v,\zeta,\xi)$ for $\Jcal$.
If we set $u:=u_1+\Pcal^{-1}v$,
then, as discussed above in the paragraph containing \eqref{motivation},
we get
  \begin{equation}
  \label{finalu}
    \Hcal[u,\zeta,\xi]=0
    \qquad \mbox{and}\qquad
    \norm{u}_{2,\alpha,\gamma} \leq 2C_1\tauj{1}.
  \end{equation}
That $u$ is actually smooth now follows from the minimality and standard regularity theory.
We have already chosen $m$ sufficiently large that $\iota[u]$ is an immersion.
By taking $m$ possibly even larger, we can guarantee embeddedness as follows.
Recalling \eqref{rhodef},
consider in the initial surface $\Sigma$ the overlapping subsets
$K:=\left\{\rho \geq m^2\right\}$ and $T:=\left\{\rho \leq m^3\right\}$,
so that $K$ has $(N-1)k \ell m^2$ components,
each contained in an isometric copy (under an element of $\Grp$) of some $\catr[i]$,
and $T$ has $N$ components, each a graph over $\T$.
By scaling $\gsph$ it is clear that there exists $\epsilon=\epsilon[N,k,\ell,\cbar]>0$
such that $\iota[u]|_K$ and $\iota[u]|_T$ are embeddings
whenever (given that they are immersions)
$\norm{u|_K: C^0(K)} <\epsilon \tauj{1}$
and $\norm{u|_T: C^0(T)} < \epsilon m^{-3}$.
Both inequalities are ensured by the estimate for $u$ in \eqref{finalu},
assuming $m>m_4$ for some $m_4=m_4[N,k,\ell,\gamma] \geq m_3$
(and, to get the first inequality,
using the decay built into the norm $\norm{\cdot}_{2,\alpha,\gamma}$ \eqref{globalweightednorm}).
Moreover, 
there is a constant $\delta=\delta[N,k,\ell]>0$
so that
the distance between any two components of $K$ is at least $\min\{\delta m^{-1},\delta m^2\tauj{1}\}$,
the distance between any any two components of $T$ is at least $\delta m^2\tauj{1}$,
and the distance between any component of $K \backslash T$ and component of $T \backslash K$
is at least $\delta m^{-2}$.
Of course $2C_1\tauj{1}<m^2\tauj{1}<m^{-2}<m^{-1}$
provided $m>m_0$ for some $m_0=m_0[N,k,\ell,\gamma] \geq m_4$.
Thus $\iota[u]$ is an embedding when $m>m_0$.
In particular its image is diffeomorphic to $\Sigma$,
so by Proposition \ref{initsurfprops} has the stated genus.
This ends the proof.
\end{proof}


\begin{bibdiv}
\begin{biblist}

\bib{BWW}{article}{
  title={Minimal surfaces in the three dimensional sphere with high symmetry}
  author={Bai, S.}
  author={Wang, C.}
  author={Wang, S.}
  journal={preprint}
  eprint={https://arxiv.org/abs/1804.05421}
}

\bib{BreS3}{article}{
  author={Brendle, S.}
  title={Minimal surfaces in $S^3$: a survey of recent results}
  journal={Bulletin of Mathematical Sciences}
  volume={3}
  pages={133--171}
  date={2013}
}


\bib{CS}{article}{
  author={Choe, J.}
  author={Soret, M.}
  title={New minimal surfaces in $S^3$ desingularizing the Clifford tori},
  journal={Mathematische Annalen},
  volume={364},
  number={3-4},
  pages={763--776},
  year={2016},
}

\bib{FPZ}{article}{
  title={Free boundary minimal surfaces in the unit 3-ball},
  author={Folha, A.}
  author={Pacard, F.}
  author={Zolotareva, T.},
  journal={manuscripta mathematica}
  year={2017}
  volume={154}
  number={3}
  pages={359--409}
}


\bib{FS}{article}{
  title={Sharp eigenvalue bounds and minimal surfaces in the ball},
  author={Fraser, A.}
  author={Schoen, R.M.},
  journal={Inventiones Mathematicae},
  volume={203},
  number={3},
  pages={823--890},
  year={2016},
}


\bib{KapDelaunay}{article}{
   author={Kapouleas, N.},
   title={Complete constant mean curvature surfaces in Euclidean
   three-space},
   journal={Annals of Mathematics},
   volume={131},
   date={1990},
   number={2},
   pages={239--330},
}

\bib{KapWenT}{article}{
  title={Constant mean curvature surfaces constructed by fusing Wente tori}
  author={Kapouleas, N.}
  journal={Inventiones Mathematicae}
  volume={119}
  number={3}
  pages={443--518}
  date={1995}
}


\bib{KapClay}{article}{
	title={Constructions of minimal surfaces by gluing minimal immersions}
	author={Kapouleas, N.}
	conference={
	  title={Global Theory of Minimal Surfaces}
	}
	book={
	  series={Clay Mathematics Proceedings}
	  volume={2}
	  publisher={American Mathematical Society}
	  address={Providence}
	  date={2005}
	}
	pages={489--524}
}


\bib{KapRS}{article}{
  title={Doubling and desingularization constructions for minimal surfaces}
  author={Kapouleas, N.}
  pages={281--325}
  book={
        title={Surveys in Geometric Analysis and Relativity celebrating Richard Schoen's
               60th birthday}
         publisher={Higher Education Press and International Press}
         address={Somerville, MA}
         date={2011}
         volume={20}
         series={Advanced Lectures in Mathematics}
       }
}


\bib{KapSphDbl}{article}{
  title={Minimal surfaces in the round three-sphere by doubling
    the equatorial two-sphere, I}
  author={Kapouleas, N.}
  journal={Journal of Differential Geometry}
  volume={106}
  year={2017}
  pages={393--449}
}

\bib{KapMcG}{article}{
  title={Minimal Surfaces in the Round Three-Sphere by Doubling the Equatorial Two-Sphere, II}
  author={Kapouleas, N.}
  author={McGrath, P.}
  journal={preprint}
  eprint={https://arxiv.org/abs/1707.08526}
}

\bib{KWtordesing}{article}{
  title={Minimal surfaces in the three-sphere by desingularizing intersecting Clifford tori}
  author={Kapouleas, N.}
  author={Wiygul, D.}
  journal={preprint}
  eprint={https://arxiv.org/abs/1701.05658}
}

\bib{KY}{article}{
	title={Minimal surfaces in the three-sphere by doubling the Clifford torus}
	author={Kapouleas, N.}
	author={Yang, S.D.}
	journal={American Journal of Mathematics}
	volume={132}
	date={2010}
	pages={257--295}
}


\bib{KPS}{article}{
  title={New minimal surfaces in $S^3$}
  author={Karcher, H.}
  author={Pinkall, U.}
  author={Sterling, I.}
  journal={Journal of Differential Geometry}
  volume={28}
  date={1988}
  pages={169--185}
}

\bib{Ket}{article}{
  title={Equivariant min-max theory}
  author={Ketover, D.}
  journal={preprint}
  eprint={https://arxiv.org/abs/1612.08692}
}

\bib{KMN}{article}{
  title={The catenoid estimate and its geometric applications}
  author={Ketover, D.}
  author={Marques, F.C.}
  author={Neves, A.}
  journal={preprint}
  eprint={https://arxiv.org/abs/1601.04514}
}

\bib{KKS}{article}{
   author={Korevaar, N.},
   author={Kusner, R.},
   author={Solomon, B.},
   title={The structure of complete embedded surfaces with constant mean
   curvature},
   journal={Journal of Differential Geometry},
   volume={30},
   date={1989},
   number={2},
   pages={465--503},
}

\bib{Law}{article}{
  title={Complete minimal surfaces in $S^3$}
  author={Lawson, H.B., Jr.}
  journal={Annals of Mathematics}
  date={1970}
  volume={92}
  pages={335--374}
 }

\bib{MNposric}{article}{
  author={Marques, F.C.}
  author={Neves, A.}
  title={Existence of infinitely many minimal hypersurfaces in positive
    Ricci curvature}
   journal={Inventiones Mathematicae},
   volume={209},
   date={2017},
   number={2},
   pages={577--616},
   issn={0020-9910},
}


\bib{PR}{article}{
  title={Equivariant minimax and minimal surfaces in geometric three-manifolds}
  author={Pitts, J.T.}
  author={Rubinstein, J.H.}
  journal={Bulletin of the American Mathematical Society}
  volume={19}
  number={1}
  date={1988}
  pages={303--309}
}


\bib{SchPC}{article}{
  title={The existence of weak solutions with prescribed singular behavior
         for a conformally invariant scalar equation}
  author={Schoen, R.M.}
  journal={Communications on Pure and Applied Mathematics}
  volume={41}
  number={3}
  date={1988}
  pages={317--392}
}


\end{biblist}
\end{bibdiv}
\end{document}